\documentclass[12pt,twoside]{bristol_thesis2}
\usepackage{cite}
\usepackage{url}
\usepackage{amsmath,amssymb,amsthm,latexsym}
\usepackage{graphicx, epsfig}
\usepackage{mathtools} 
\usepackage{tikz}
\usetikzlibrary{decorations.markings} 

\usepackage{mathrsfs} 

\newcommand{\kommentar}[1]{}

\numberwithin{equation}{section} \numberwithin{lemma}{section}

\begin{document}

\title{Statistics of the zeros of $L$-functions and arithmetic correlations}

\author{Dale Smith}

\titlepage

\chapter*{Abstract}
This thesis determines some of the implications of non-universal and emergent universal statistics on arithmetic correlations and fluctuations of arithmetic functions, in particular correlations amongst prime numbers and the variance of the expected number of prime numbers over short intervals are generalised by associating these concepts to $L$-functions arising from number theoretic objects.

Inspired by work in quantum chaology, which shares the property of displaying emergent universality, in chapter $2$ a heuristic is given to determine the behaviour of a correlation function associated to functions in the Selberg class from the universal form of the $2$-point correlation statistic conjectured for this class. Also in this chapter, the Riemann zeta function is taken as an example of an $L$-function from which the correlations between pairs of prime numbers arise from a non-universal form of the $2$-point correlation statistic for its zeros. Chapter $3$ explores the implications of the $2$-point correlation statistic on an arithmetic variance associated to functions in the Selberg class, generalising the variance of primes in short intervals. Many of the ideas in this thesis are based on the paper, \cite{Buionthevarof2016}, joint with Hung Bui and Jon Keating.

\pagenumbering{roman}

\chapter*{Acknowledgements}
I would like to thank my supervisor, Jon Keating, for providing me with valuable opportunities during my time in Bristol. I'm also very grateful to Hung Bui and Ghaith Hiary for many helpful discussions, to Brian Conrey and Chris Hughes for taking the time to conduct my viva and provide various useful comments, and to Aman Arora and Adam Moas for sharing their programming expertise with me. Finally, I wish to express my gratitude to my family for their continued support. 

\chapter*{Author's declaration} 
I declare that the work in this dissertation was carried out in accordance with the requirements of 
the University's Regulations and Code of Practice for Research Degree 
Programmes and that it 
has not been submitted for any other academic award. Except where indicated by specific 
reference in the text, the work is the candidate's own work. Work done in collaboration with, or with 
the assistance of, others, is indicated as 
such. Any views expressed in the dissertation are those of 
the author. \\

Dale Smith, \today.

\listoftables

\listoffigures

\tableofcontents

\chapter*{Notation} \thispagestyle{empty}
A list of commonly used notation is given here, where used otherwise this is indicated in the text.
\section*{Symbols}
\begin{align*}
&m \textrm{ and } n &\textrm{natural numbers} \\
&p \textrm{ and } q &\textrm{prime numbers} \\
&x \textrm{ and } y &\textrm{real variables} \\
&s &\textrm{complex variable}
\end{align*}
\section*{Operations}
\begin{align*}
&* \textrm{ and } cc &\textrm{complex conjugate} \\
&\Im &\ \textrm{imaginary part} \\
&\Re & \textrm{real part} \\
&' & \textrm{derivative} \\
&\int & \textrm{integral} \\
&\mathscr{F} &\textrm{Fourier transform}
\end{align*}
\section*{Conventions}
\begin{itemize}
\item The integral, $\int,$ is to be taken over the whole domain of the variables of integration unless specified otherwise by a subscript. For integration over a domain, $D,$ this is denoted by $\int_D.$
\item A function, $f,$ with the Fourier transform, $\mathscr{F}[f](k),$ is given by the convention $\mathscr{F}[f](k) = \int f(x) \exp(- i k x) \, dx,$ so that $f(x) = \frac{1}{2 \pi} \int \mathscr{F}[f](k) \exp(i k x) \, dk.$
\end{itemize}
\section*{Functions and distributions} \thispagestyle{empty}
\begin{align*}
&f(s) \textrm{ and } g(s) &\textrm{generic functions} \\
&1_A &\textrm{indicator function for condition } A \\
&\Gamma(s) &\textrm{gamma function} \\
&F(s) &\textrm{function in the Selberg class} \\
&m(F) &\textrm{polar order of } F \\
&\textrm{deg}(F) &\textrm{degree of } F \\
&\mathfrak{q}(F) &\textrm{conductor of } F
\end{align*}
\begin{align*}
&\Lambda(n; F) &\textrm{von Mangoldt function associated to } F \\
&\psi(x, h; F) &\sum_{x < n \leq x + h} \Lambda(n; F) \\
&Var^{mul}(X, \delta; F) &\frac{1}{X} \int_{1 \leq x \leq X} |\psi(x, \delta x; F) - m(F) \delta x|^2 \, dx \\
&Var^{fix}(X, h; F) &\frac{1}{X} \int_{1 \leq x \leq X} |\psi(x, h; F) - m(F) h|^2 \, dx \\
&\textrm{Si}(x) &\textrm{sine integral } \\
&\textrm{erf}(x) &\textrm{error function } \\
&\delta(x) &\textrm{Dirac delta distribution}
\end{align*}
\section*{Relations}
In some domain, $f$ is said to be Big O of $g$ or less than less than $g$ if 
$$|f| < c |g|,$$ \thispagestyle{empty}
where $c$ is a constant. If $c$ depends on $x$ then this is written
$$f = O_x(g) \textrm{ or } f \ll_x g$$ respectively, with the subscript dropped otherwise.

\chapter{Introduction} \pagenumbering{arabic} 
\nocite{MR1077246} \nocite{MR2477804} \nocite{MR2175035} \nocite{MR2128842} \nocite{MR2251029} \nocite{MR1710182} \nocite{MR1946551} \nocite{MR2800717} \nocite{MR2145172} \nocite{berry1986riemann} \nocite{MR2360010} \nocite{MR1853612} \nocite{zbMATH03783102} This thesis begins by introducing concepts in $2$ apparently disparate fields of study, namely quantum chaology and analytic number theory. However, as will emerge, there are similarities that provide intrigue for further investigation.
\section{Quantum chaology}
Quantum chaology is the study of quantum systems with classically chaotic counterparts. A remarkable feature of such systems is that the behaviour of their energy spectra may be deduced from the underlying classical dynamics, in particular via the periodic orbits in the system.
\subsection{Spectral behaviour}
Gutzwiller's trace formula  provides the leading order $\hbar \to 0$ asymptotic for the spectral energy density,
\begin{equation} \label{Spectral density}
d(e; D) = \sum_{n \geq 1} \delta(e - e(n; D)),
\end{equation} 
where $D$ corresponds to a system with classically chaotic dynamics with energy spectrum $e(n; D) \textrm{ for } n = 1, 2, \ldots.$
\begin{equation}
d(e; D) \sim d^{mean}(e; D) + d^{osc}(e; D),
\end{equation}
with mean density,
\begin{equation} \label{Mean spectral density}
d^{mean}(e; D) = \Omega(e(D)) / (2 \pi \hbar)^f,
\end{equation}
where $f$ is the number of degrees of freedom of $D$ and $\Omega(e(D))$ is the volume of the energy shell.
The oscillating part of the spectral density, $d^{osc}(e; D),$ is given by
\begin{equation} \label{Oscillating part of spectral density}
\frac{1}{2 \pi \hbar}\sum_{\substack{p \\ k \geq 1}} A(p, k; D) T(p; D) \exp(i k / \hbar (S(p; D) - \pi \mu(p; D) / 2)) + cc,
\end{equation}
where $A(p, k; D), T(p; D), S(p; D) \textrm{ and } \mu(p; D)$ are respectively the stability amplitude, the period, the action and the Maslov index associated to primitive orbits, $p,$ and their repetitions, $k,$ in $D.$ In particular the stability amplitude for the $k$th repetition of $p$ is given by
\begin{equation}
A(p, k; D) = |\textrm{det}(M^k(p; D) - I(D))|^{- 1 / 2},
\end{equation}
where $M(p; D)$ is the monodromy matrix associated to $D$ and $I(D)$ is the identity matrix, with its dimension equal to $M(p; D).$

The $2$-point correlation statistic,
\begin{equation} \label{2-point correlation for D}
R_2(x, E; D) = \left<d(e + x_1; D) d(e + x_2; D)\right>_E,
\end{equation}
provides information on the distribution of energies averaged over an interval around $E,$ specifically on pairs of energies with the separation $x = x_2 - x_1.$ Substituting \eqref{Mean spectral density} and \eqref{Oscillating part of spectral density} into \eqref{2-point correlation for D} gives a sum over the primitive orbits $p \textrm{ and } q,$ and their repetitions $k \textrm{ and } l,$ resulting in
\begin{equation} \label{Terms in 2-point correlation for D}
R_2(x, E; D) = (d^{mean}(E; D))^2 + R_2^d(x, E; D) + R_2^o(x, E; D),
\end{equation}
with $R_2^d(x, E; D)$ corresponding to the diagonal terms, with $k S(p; D) = l S(q; D),$ and $R_2^o(x, E; D)$ to the off-diagonal terms, with $k S(p; D) \neq l S(q; D).$ 

Bogomolny and Keating calculated this statistic in \cite{bogomolny1996gutzwiller}. As described in \cite{keating1999periodic}, there is an exact formula for the diagonal contribution,
\begin{equation}
R_2^d(x, E; D) = \frac{1}{2 \pi^2 \hbar^2} \Re\left[\frac{d^2}{dw^2} \log Z(w; D) - a(w; D)\right]_{w = i x / \hbar},
\end{equation}
where
\begin{equation}
1 / Z(w; D) = \prod_p Z(w, p; D)
\end{equation}
with
\begin{equation}
Z(w, p; D) = \prod_{n \geq 0} \left(1 - \frac{\exp(w T(p; D))}{|\Xi(p; D)| \Xi^n(p; D)}\right)^{(n + 1) g(p; D)}
\end{equation}
and
\begin{equation}
a(w; D) = \sum_{\substack{p \\ n \geq 0}} \frac{(n + 1) g(p; D) T^2(p; D)}{(|\Xi(p; D)| \Xi^n(p; D) \exp(- w T(p; D)) - 1)^2},
\end{equation}
where $g(p; D)$ is the multiplicity of the $p$th primitive orbit and $\Xi(p; D)$ the largest eigenvalue of $M(p; D).$ 

Furthermore, there is a non-zero off-diagonal contribution, indicating the presence of correlations between the actions of distinct orbits. For the case of ergodic systems without time-reversal symmetry, $Z(w; D)$ has a pole at $w = 0$ with residue $\gamma,$ and
\begin{align}
R_2^o(x, E; D) = \left|Z(i x / \hbar; D) / \sqrt{2} \pi \hbar \gamma\right|^2 \nonumber \\
\Re\left(\exp(2 \pi i d^{mean}(E; D) x) \prod_p \chi(p, x / \hbar; D)\right), \label{Off-diagonal contribution}
\end{align}
where
\begin{align}
\chi(p, x; D) = |Z(0, p; D) / Z(i x, p; D)|^2 \nonumber \\
_2\phi_1(a(p; D), b(p; D); c(p; D); q(p; D), z(p; D)),
\end{align}
with $a(p; D) = b(p; D) = \exp(- i x T(p; D)),$ $c(p; D) = q(p; D) = \Xi^{- 1}(p; D),$ $z(p; D) = |\Xi(p; D)|^{- 1} \exp(i x T(p; D)),$ and
\begin{equation}
_2\phi_1(a, b; c; q, z) = \sum_{n \geq 0} \frac{(a; q)_n (b; q)_n}{(c; q)_n} \frac{z^n}{(q; q)_n}
\end{equation}
is the $q$-hypergeometric series, with the $q$-Pochhammer symbol given by
\begin{equation}
(d; q)_n = \begin{dcases*} 1 & if $n = 0$ \\ \prod_{0 \leq i \leq n - 1} (1 - d q^i) & if $n > 0.$ \end{dcases*}
\end{equation}

Correlations between off-diagonal terms are however unknown in general, as such \eqref{Off-diagonal contribution} was obtained indirectly, from the diagonal terms.

For a generic non-time-reversal symmetric system $g(p; D)$ is replaced by $1,$ rescaling the $2$-point correlation then gives the system independent expression
\begin{equation} \label{Universal diagonal contribution}
\tilde{R}_2^d(x; D) = \lim_{E \to \infty} \frac{R_2^d(x / d^{mean}(E; D), E; D)}{(d^{mean}(E; D))^2} = - \frac{1}{2 \pi^2 x^2}.
\end{equation}
Rescaling the off-diagonal contribution in the same way, the analytic structure of $Z$ around the pole gives
\begin{equation} \label{Universal off-diagonal contribution}
\tilde{R}_2^o(x; D) = \frac{\cos 2 \pi x}{2 \pi^2 x^2},
\end{equation}
together with the diagonal contribution giving a universal form for the rescaled $2$-point correlation;
\begin{equation}
\tilde{R}_2(x; D) = 1 - \left(\frac{\sin \pi x}{\pi x}\right)^2.
\end{equation}
This emergent universality coincides with a celebrated conjecture of Bohigas, Giannoni and Schmit. In \cite{MR730191} it was conjectured that the form of the universal spectral statistics are purely a consequence of the symmetries present in the system under consideration, for further information on the BGS conjecture several chapters of \cite{MR2920518} are recommended.
\subsection{Correlations in the actions of \\ periodic orbits}
Assuming the BGS conjecture, in \cite{MR1253066} correlations between the actions of distinct orbits were quantified and compared with numerics for specific systems. Specialising to systems with $2$ degrees of freedom and even Maslov indicies, the weighted correlation statistic,
\begin{align} 
P(x, T; D) = \sum_{\substack{p, q \\ k, l \geq 1 \\ (p, k) \neq (q, l)}} (- 1)^{(\mu(p; D) - \mu(q; D)) / 2} A(p, k; D) A(q, l; D) \nonumber \\ 
\delta(x - (S(p, k; D) - S(q, l; D)) \Delta T f(T - T(p; D)) f(T - T(q; D)) \label{Action correlations}
\end{align}
was studied, where
\begin{equation}
f(y) = \sqrt{2} / (\Delta T) \exp(- 2 \pi / (\Delta T)^2 y^2)
\end{equation}
for the time interval $\Delta T,$ with the pairs $(p, k) \textrm{ and } (q, l)$ denoting the $k$th repetition of the $p$th periodic orbit and the $l$th repetition of the $q$th periodic orbit respectively. Rescaling \eqref{Action correlations} to
\begin{equation} \label{Rescaled action correlations}
\hat{P}(z) = \frac{P\left(\Omega(E; D) / (2 \pi T) z, T; D\right)}{(2 \pi T^2) / \Omega(E; D)}
\end{equation}
for their numerics, this was determined by inverting a general expression for the spectral statistics. Using the off-diagonal spectral correlations given by the universal form \eqref{Universal off-diagonal contribution} for non-time-reversal symmetric systems, \eqref{Rescaled action correlations} was shown to be well approximated by
\begin{equation} \label{Rescaled averaged action correlations}
- 2 / \pi \left(\frac{\sin z / 2}{z}\right)^2.
\end{equation}

Making the diagonal approximation, that is neglecting any correlations between periodic orbits, in \cite{MR805089} Berry found the universal behaviour given by \eqref{Universal diagonal contribution}. The correlations indicated by \eqref{Rescaled averaged action correlations} were subsequently found to be related to crossings of periodic orbits, the underlying mechanism for the leading order off-diagonal contributions, \eqref{Universal off-diagonal contribution}, now known as the Sieber-Richter pair due to the work in \cite{Sieber}. 

More precise action correlation information may be recovered from the added structure available in \eqref{Off-diagonal contribution}. Analogous correlations will be deduced in the next chapter in the simpler case of a class of model systems coming from analytic number theory, where the analogues of \eqref{Rescaled averaged action correlations} and its more precise counterpart are of an arithmetic nature.
\section{The Selberg class}
Selberg's eponymous class was introduced in \cite{MR1220477} as an attempt to provide an axiomatic foundation for $L$-functions arising in analytic number theory. A function, $F,$ in the Selberg class, is given by
\begin{equation}
F(s) = \sum_{n \geq 1} \frac{a(n; F)}{n^s} = \exp \left(\sum_{n \geq 1} \frac{b(n; F)}{n^{s}}\right)
\end{equation}
for $\Re(s) > 1,$ where for every $\epsilon > 0$ the complex multiplicative function
\begin{equation}
a(n; F) = O(n^{\epsilon}),
\end{equation}
with $a(1; F) = 1.$ The complex function
\begin{equation}
b(n; F) = O(n^{\theta})
\end{equation}
for some $\theta < 1 / 2,$ with $b(n; F) = 0$ unless $n$ is a power of a prime number. A function in the Selberg class is therefore specified by the values of its coefficients at prime powers,
\begin{equation}
F(s) = \prod_p \left(\sum_{l \geq 0} \frac{a(p^l; F)}{p^{l s}}\right) = \prod_p \exp\left(\sum_{k \geq 1} \frac{b(p^k; F)}{p^{k s}}\right). 
\end{equation}
The coefficients may be related at prime powers by equating the derivatives,
\begin{equation}
- \log p \sum_{l \geq 0} \frac{l a(p^l; F)}{p^{l s}} = - \log p \sum_{m \geq 0} \frac{a(p^m; F)}{p^{m s}} \sum_{k \geq 1} \frac{k b(p^k; F)}{p^{k s}},
\end{equation}
so that
\begin{equation}
\sum_{n \geq 0} \frac{n a(p^n; F)}{p^{n s}} = \sum_{n \geq 0} \frac{1}{p^{n s}} \sum_{\substack{m \geq 0, k \geq 1 \\ n = m + k}} k a(p^m; F) b(p^k; F)
\end{equation}
and
\begin{equation}
a(p^n; F) = \frac{1}{n} \sum_{1 \leq k \leq n} k a(p^{n - k}; F) b(p^k; F).
\end{equation}

Functions in the Selberg class are also endowed with a meromorphic continuation to the rest of the complex plane, with a pole of order $m(F) \geq 0 \textrm{ at } s = 1$ and a functional equation, 
\begin{equation}
F(s) = X(s; F) \overline{F}(1 - s),
\end{equation}
where
\begin{equation}
X(s; F) = \epsilon(F) Q^{1 - 2 s} \prod_{1 \leq j \leq r} \frac{\Gamma(\lambda_j (1 - s) + \mu_j^*)}{\Gamma(\lambda_j s + \mu_j)}
\end{equation}
for some $|\epsilon(F)| = 1, \, Q > 0, \, \lambda_j > 0 \textrm{ and } \Re(\mu_j) \geq 0 \textrm{ for } 1 \leq j \leq r \textrm{ and } r \geq 0, \textrm{ where } r = 0 \textrm{ gives an empty product of unity.}$ Given a function, $G,$ the overline notation is used for the operation
\begin{equation}
\overline{G}(s) = G^*(s^*).
\end{equation}
The form of the functional equation is not uniquely determined by $F,$ for example applying Legendre's duplication formula,
\begin{equation}
\Gamma(s) = \pi^{- 1 / 2} 2^{s - 1} \Gamma((s + 1) / 2) \Gamma(s / 2),
\end{equation}
to the gamma functions appearing in $X(s; F),$ doubles the value of $r.$ However the degree of $F,$
\begin{equation}
\textrm{deg}(F) = 2 \sum_{1 \leq j \leq r} \lambda_j,
\end{equation}
is invariant under changes to the form of the functional equation. Another invariant is the conductor of $F,$ 
\begin{equation}
\mathfrak{q}(F) = (2 \pi)^{\textrm{deg}(F)} Q^2 \prod_{1 \leq j \leq r} \lambda_j^{2 \lambda_j}.
\end{equation}

By the exponential form of $F$ there are no zeros in the region $\Re(s) > 1.$ Furthermore, by the functional equation $F$ is also non-zero for $\Re(s) < 0$ except for the points $- \mu_j / \lambda_j \textrm{ for } 1 \leq j \leq r, \textrm{ if } \mu_j = 0, \textrm{ and } - (m + \mu_j) / \lambda_j \textrm{ for } 1 \leq j \leq r \textrm{ and } m = 1, 2, \ldots$ otherwise, corresponding to the poles of the gamma function.

The leading order $|t| \to \infty$ asymptotics for the analogue of \eqref{Spectral density},
\begin{equation}
d(t; F) = \sum_{n \geq 1} \delta(t - t(n; F)),
\end{equation}
where the $n$th non-trivial zero of $F$ is denoted by $\sigma(n; F) + i t(n; F),$ have a similar form to \eqref{Mean spectral density} and \eqref{Oscillating part of spectral density}. For $r >0$ the mean density is given by
\begin{equation} \label{Mean density}
d^{mean}(t; F) = \frac{1}{2 \pi} \log \left(\mathfrak{q}(F) \left(\frac{|t|}{2 \pi}\right)^{\textrm{deg}(F)}\right)
\end{equation}
and the oscillating contribution is given formally by,
\begin{equation} \label{Oscillating part of density}
d^{osc}(t; F) = - \frac{1}{2 \pi} \sum_{\substack{p \\ k \geq 1}}  b(p^k; F) \log p^k p^{- k / 2} \exp(- i k |t| \log p) + cc.
\end{equation}
As in the case of quantum chaology, it is not feasible to directly calculate the 2-point correlation, since in general the inherent correlations are unknown, in particular for the function
\begin{equation}
\Lambda(n; F) = b(n; F) \log n,
\end{equation}
the negative values of which serve as coefficients for the logarithmic derivative of $F$,
\begin{equation}
\frac{F'}{F}(s) = - \sum_{n \geq 1} \frac{\Lambda(n; F)}{n^s} \textrm{ for } \Re(s) > 1.
\end{equation}

In \cite{MR1981179} it was shown that as $x \to \infty,$ for any function in the Selberg class 
\begin{equation} \label{Prime number theorem for the Selberg class}
\sum_{n \leq x} \Lambda(n; F) \sim m(F) x,
\end{equation}
assuming a conjecture made by Selberg in \cite{MR1220477}, namely Selberg's orthogonality conjecture, which says that
\begin{equation}
\sum_{p \leq x} \frac{a(p; F) a^*(p; G)}{p} = 1_{F = G} \log \log x + O(1),
\end{equation}
where $F \textrm{ and } G$ are primitive functions in the Selberg class, that is, they cannot be expressed as a product of other functions in the Selberg class in a non-trivial way. Assuming instead, as conjectured by Selberg in \cite{MR1220477}, that all non-trivial zeros of any function in the class have a real part of $1 / 2,$ there is the explicit formula
\begin{align}
& \sum_{n \leq x} \Lambda(n; F) = m(F) x - \sum_{n \geq 1} \frac{x^{1 / 2 + i t(n; F)}}{1 / 2 + i t(n; F)} 1_{|t(n; F)| \leq Z} \\ \nonumber
& + O(\log x \min\{1, x / (Z ||x||)\}) + O(x (\log x Z)^2 / Z),
\end{align}
where $||x||$ is the distance from $x$ to the closest integer, which follows as in chapter 17 of \cite{MR1790423}.

Under the same assumption on the horizontal positions in the complex plane of the non-trivial zeros, as well as Selberg's orthogonality conjecture, in \cite{MR1692847} Murty and Perelli studied
\begin{equation}
\mathcal{F}(X, T; F, G) = \sum_{m, n \geq 1} X^{i (t(m; F) - t(n; G))} w(t(m; F) - t(n; G)) 1_{|t(m; F)|, |t(n; F)| < T},
\end{equation}
where 
\begin{equation}
w(u) = \frac{4}{4 + u^2}.
\end{equation}
Fixing the primitive function $F,$ and allowing the other primitive function, $G,$ to vary, they made the conjecture that
\begin{equation}
\mathcal{F}(X, T; F, G) = \frac{1_{F = G} T \log X}{\pi} + \frac{\textrm{deg}(F) \textrm{deg}(G) T (\log T)^2}{\pi X^2} (1 + o(1)) + o(1)
\end{equation}
uniformly for $T^{A_1} \leq X \leq T^{A_2} \textrm{ for any fixed } A_1 \textrm{ and } A_2 \textrm{ satisfying } 0 < A_1 \leq A_2 \leq \textrm{deg}(F)$ and
\begin{equation}
\mathcal{F}(X, T; F, G) = \frac{\textrm{deg}(F) T \log T}{\pi} (1_{F = G} + o(1))
\end{equation}
uniformly for $T^{A_1} \leq X \leq T^{A_2} \textrm{ for any fixed } A_1 \textrm{ and } A_2 \textrm{ satisfying } \textrm{deg}(F) \leq A_1 \leq A_2.$

Murty and Perelli showed that unique factorisation into primitive functions follows from their conjecture for $\mathcal{F}(X, T; F, G),$ a result shown in \cite{MR1253620} to be a consequence of Selberg's orthogonality conjecture. Their conjecture also has implications for the fluctuations about the average behaviour given in \eqref{Prime number theorem for the Selberg class} over short intervals. In chapter $3$ the variances
\begin{equation}
Var^{mul}(X, \delta; F) = \frac{1}{X} \int_{1 \leq x \leq X} |\psi(x, \delta x; F) - m(F) \delta x|^2 \, dx
\end{equation}
and
\begin{equation}
Var^{fix}(X, h; F) = \frac{1}{X} \int_{1 \leq x \leq X} |\psi(x, h; F) - m(F) h|^2 \, dx,
\end{equation}
are studied over short multiplicative and fixed intervals respectively, where
\begin{equation}
\psi(x, h; F) = \sum_{x < n \leq x + h} \Lambda(n; F),
\end{equation}
employing a more precise form of the conjecture made by Murty and Perelli with $F = G.$ In doing so, functions are combined via the tensor product
\begin{equation}
(F \otimes \overline{G})(s) = \prod_{p} \textrm{exp} \left(\sum_{l \geq 1} \frac{l b(p^l; F) b^*(p^l; G)}{p^{l s}}\right).
\end{equation}
Its logarithmic derivative is given by
\begin{equation}
\frac{(F \otimes \overline{G})'}{(F \otimes \overline{G})}(s) = - \sum_{\substack{p \\ l \geq 1}} \frac{l^2 b(p^l; F) b^*(p^l; G) \log p}{p^{l s}},
\end{equation}
so for $\Re(s) > 1$
\begin{equation}
\frac{(F \otimes \overline{G})'}{(F \otimes \overline{G})}(s) = - \sum_p \frac{a(p; F) a^*(p; G) \log p}{p^s} + O(1).
\end{equation}
Applying partial summation and making use of Selberg's orthogonality conjecture gives
\begin{align}
\sum_{p \leq x} \frac{a(p; F) a^*(p; G)}{p^{1 + \sigma}} = & \sigma \int_{1 \leq t \leq x} \frac{1_{F = G} \log \log t + O(1)}{t^{1 + \sigma}} \, dt \nonumber \\
& + O\left(\frac{1_{F = G} \log \log x}{x^{\Re(\sigma)}}\right),
\end{align}
so that
\begin{align}
\sum_p \frac{a(p; F) a^*(p; F)}{p^{1 + \sigma}} = & \sigma \int_{t \geq 1} \frac{\log \log t}{t^{1 + \sigma}} \, dt + O(1) \nonumber \\
= & \log \frac{1}{\sigma} + O(1),
\end{align}
showing that
\begin{equation}
\frac{(F \otimes \overline{F})'}{(F \otimes \overline{F})}(1 + \sigma) \sim - \frac{1}{\sigma},
\end{equation}
and the presence of a simple pole at $s = 1$ for $(F \otimes \overline{G})(s)$ when $F = G.$
\subsection{Examples of $L$-functions}
An alternative to the Selberg class formulation of $L$-functions is to consider them as arising from automorphic representations of the general linear group over a number field, with primitive $L$-functions of degree $d$ attached to cuspidal automorphic representations of $GL(d) / \mathbb{Q}.$ Historically, the first $L$-function studied was the Riemann zeta function.
\subsubsection{The Riemann zeta function}
Originally studied by Euler as a function of a real variable, Riemann first treated\footnote{Riemann's influencial paper can be found in the appendix of \cite{MR1854455}}
\begin{equation}
\zeta(s) = \sum_{n \geq 1} \frac{1}{n^s}
\end{equation}
as a function of a complex variable, for $\Re(s) > 1,$ and then continued to the rest of the complex plane, except for $s = 1$ where there is a pole of order $1.$ The zeros occur in complex conjugate pairs since $\zeta(s^*) = (\zeta(s))^*,$ also by the Euler product representation for $\Re(s) > 1,$
\begin{equation}
\zeta(s) = \prod_p (1 - p^{- s})^{- 1},
\end{equation}
there are no zeros in this region of the complex plane. Riemann showed that there is a functional equation,
\begin{equation}
\zeta(s) = \pi^{s - 1 / 2} \frac{\Gamma(\frac{1 - s}{2})}{\Gamma(\frac{s}{2})} \zeta(1 - s),
\end{equation}
so that the only zeros with $\Re(s) < 0$ occur at $s = - 2 n \textrm{ for } n = 1, 2, \ldots.$ The functional equation also shows that each zero generates another by reflection about the critical line, $\Re(s) = 1 / 2,$ where Riemann's famous hypothesis asserts all the non-trivial zeros lie. Upon stating what is now known as the Riemann hypothesis, Riemann sketched a proof of the prime number theorem,
\begin{equation}
\sum_{n \leq x} \Lambda(n) \sim x \textrm{ as } x \to \infty,
\end{equation}
where 
\begin{equation}
\Lambda(n) = \begin{dcases*} \log p & if $n$ is a power of a prime, $p$ \\ 0 & otherwise \end{dcases*}
\end{equation}
is the von Mangoldt function. The prime number theorem was subsequently proved independently by De la Val{\'e}e Poussin and Hadamard by showing that there are no zeros with real part equal to $1,$ for further details one of the many texts covering the theory of the Riemann zeta function may be consulted, for example \cite{MR0434929}, \cite{MR882550}, \cite{MR1854455} or \cite{MR933558}.

The statistic $\mathcal{F}(X, T; F, G)$ was first studied in \cite{MR0337821}. Montgomery considered
\begin{align}
\mathcal{F}(X, T) = & \sum_{m, n \geq 1} X^{i (\gamma(m) - \gamma(n))} w(\gamma(m) - \gamma(n)) 1_{0 < \gamma(m), \gamma(n) < T} \nonumber \\
= & \frac{1}{2} \mathcal{F}(X, T; \zeta, \zeta),
\end{align}
where assuming the Riemann hypothesis the $m$th zero is denoted $1 / 2 + i \gamma(m).$ Montgomery proved that
\begin{equation}
\mathcal{F}(X, T) = \frac{T \log X}{2 \pi} + \frac{T (\log T)^2}{2 \pi X^2} (1 + o(1)) + o(1)
\end{equation}
uniformly for $T^{A_1} \leq X \leq T^{A_2} \textrm{ for any fixed } A_1 \textrm{ and } A_2 \textrm{ satisfying } 0 < A_1 \leq A_2 < 1$ and conjectured that 
\begin{equation}
\mathcal{F}(X, T) = \frac{T \log T}{2 \pi} (1 + o(1))
\end{equation}
uniformly for $T^{A_1} \leq X \leq T^{A_2} \textrm{ for any fixed } A_1 \textrm{ and } A_2 \textrm{ satisfying } 1 \leq A_1 \leq A_2.$ His conjecture was motivated by the work of Hardy and Littlewood on the distribution of pairs of primes, who in \cite{MR1555183} conjectured that
\begin{equation} \label{Hardy-Littlewood conjecture}
\lim_{N \to \infty} \frac{1}{N} \sum_{n \leq N} \Lambda(n + r) \Lambda(n) = \begin{dcases*} 2 \prod_p \left(1 - \frac{1}{(p - 1)^2}\right) \prod_{q | r} \frac{q - 1}{q - 2} & if $r$ is even \\ 0 & otherwise, \end{dcases*}
\end{equation}
where $p \textrm{ and } q \textrm{ are both primes greater than } 2.$ Extensive numerical investigations support Montgomery's conjecture, an informative survey of results is given in \cite{MR1868473}.

In \cite{MR1018376} an equivalence was established, relating fluctuations about the leading term given by the prime number theorem over short intervals to Montgomery's conjecture. Assuming the Riemann hypothesis it was shown that
\begin{equation}
\mathcal{F}(X, T) \sim \frac{T \log T}{2 \pi} 
\end{equation}
uniformly for $T \leq X \leq T^A \textrm{ for any fixed } A \geq 1,$
\begin{equation}
\int_{1 \leq x \leq X} \left(\sum_{x < n \leq (1 + \delta) x} \Lambda(n) - \delta x \right)^2 \, dx \sim \frac{1}{2} \delta X^2 \log \frac{1}{\delta}
\end{equation}
uniformly for $X^{- 1} \leq \delta \leq X^{- \epsilon} \textrm{ for every fixed } \epsilon > 0,$ and 
\begin{equation} \label{Goldston and Montgomery with h fixed}
\int_{0 \leq x \leq X} \left(\sum_{x < n \leq x + h} \Lambda(n) - h \right)^2 \, dx \sim h X \log \frac{X}{h}
\end{equation}
uniformly for $1 \leq h \leq X^{1 - \epsilon} \textrm{ for every fixed } \epsilon > 0,$ are all equivalent to one another. Making an assumption on the form of the error term in the Hardy-Littlewood conjecture for finite $N,$ in \cite{MR2104891} \eqref{Goldston and Montgomery with h fixed} was shown to follow with the lower order terms $- h X (\gamma + \log 2 \pi) + O\left(h^{15/16} X (\log X)^{17/16} + h^2X^{1 / 2 + \epsilon}\right)$ uniformly for $\log X \leq h \leq X^{1 / 2}$ and some $\epsilon > 0,$ where $\gamma$ is Euler's constant, with no assumptions made on the distribution of the zeros of the Riemann zeta function. In fact, assuming a $k$-tuple form of the Hardy-Littlewood conjecture with the same assumption on the error term it was shown that
\begin{align}
& \int_{0 \leq x \leq X} \left(\sum_{x < n \leq x + h} \Lambda(n) - h \right)^k \, dx \nonumber \\ 
= & \mu_k h^{k / 2} \int_{1 \leq x \leq X} (\log x / h + 1 - \gamma - \log 2 \pi)^{k / 2} \, dx \nonumber \\
& + O(X (h \log X)^{k / 2} (h / \log X)^{- 1 / (8 k)} + h^k X^{1 / 2 + \epsilon})
\end{align}
uniformly for $\log X \leq h \leq X^{1 / k} \textrm{ and some } \epsilon > 0,$ where
\begin{equation}
\mu_k = \begin{dcases*} \prod_{0 \leq j \leq (k - 2) / 2} (2 j + 1) & if $k$ is even \\ 0 & otherwise, \end{dcases*}
\end{equation}
for $k = 1, 2, \ldots,$ are the centralised moments of a standard normal distribution.

Assuming the lower order terms given in \cite{MR2104891} to \eqref{Goldston and Montgomery with h fixed}, in \cite{MR2009438} the arguments in \cite{MR1018376} were refined, leading to a more precise expression for the conjectural behaviour of $\mathcal{F}(X, T),$
\begin{equation}
\mathcal{F}(X, T) = \frac{T}{2 \pi} \log \frac{T}{2 \pi} - \frac{T}{2 \pi} + o(T)
\end{equation}
uniformly for $T^{1 + \epsilon} \leq X \leq T^A \textrm{ for every fixed } \epsilon > 0 \textrm{ and } A \geq 1 + \epsilon.$ A more precise form of Montgomery's conjecture was first given in \cite{bogomolny1996gutzwiller}, expressed via the 2-point correlation for the zeros of the Riemann zeta function,
\begin{equation} \label{2-point correlation for zeta}
R_2(x, T) = \left<d(t + x_1) d(t + x_2)\right>_T,
\end{equation}
where the mean and oscillating contribution to the density of Riemann zeros are given by
\begin{equation} \label{Mean density for zeta}
d^{mean}(t) = \frac{1}{2 \pi} \log \frac{|t|}{2 \pi}
\end{equation}
and
\begin{equation} \label{Oscillating part of density for zeta}
d^{osc}(t) = - \frac{1}{\pi} \sum_{n \geq 1} \frac{\Lambda(n)}{n^{1 / 2}} \cos(t \log n)
\end{equation}
respectively. As exposited in \cite{keating1999periodic}, the diagonal contribution is controlled by the behaviour of the Riemann zeta function on the line $1$-line, $s = 1 + i w \textrm{ with } w \in \mathbb{R}.$ Substituting \eqref{Mean density for zeta} and \eqref{Oscillating part of density for zeta} into \eqref{2-point correlation for zeta} gives the analogue of \eqref{Terms in 2-point correlation for D} for the Riemann zeta function,
\begin{equation}
R_2(x, T) = (d^{mean}(T))^2 + R_2^d(x, T) + R_2^o(x, T),
\end{equation}
where
\begin{equation}
R_2^d(x, T) = \frac{1}{2 \pi^2} \sum_{n \geq 1} \frac{(\Lambda(n))^2}{n} \cos (x \log n)
\end{equation}
and
\begin{equation}
R_2^o(x, T) = \frac{1}{2 \pi^2} \sum_{n \geq 1, \, r \neq 0} \frac{\Lambda(n + r) \Lambda(n)}{n} \cos(x \log n - T r / n).
\end{equation}
The diagonal contribution may be expressed as
\begin{align}
R_2^d(x, T) & = \frac{\Re}{2 \pi^2} \sum_{\substack{p \\ k \geq 1}} \frac{(\log p)^2}{p^{(1 + i x) k}} \nonumber \\
& = \frac{\Re}{2 \pi^2} \left(\sum_{\substack{p \\ k \geq 1}} \frac{k (\log p)^2}{p^{(1 + i x) k}} - \sum_{\substack{p \\ k \geq 1}} \frac{(k - 1) (\log p)^2}{p^{(1 + i x) k}}\right).
\end{align}
Evaluating the logarithm of the conditionally convergent Euler product on the $1$-line gives 
\begin{equation}
\sum_{\substack{p \\ k \geq 1}} \frac{1}{k p^k} \exp(- i k w \log p),
\end{equation}
resulting in
\begin{equation}
R_2^d(x, T) = - \frac{\Re}{2 \pi^2} \left(\frac{d^2}{dw^2} \log \zeta(1 + i w) + \sum_{\substack{p \\ k \geq 0}} (\log p)^2 k p^{- (1 + i x) (k + 1)}\right)_{w = x}.
\end{equation}
The off-diagonal contribution may be determined by making use of the Hardy-Littlewood conjecture, and was shown by Bogomolny and Keating to be given by
\begin{equation}
R_2^o(x, T) = \frac{1}{2 \pi^2} |\zeta(1 + i x)|^2 \Re\left(\exp(2 \pi i x d^{mean}(T)) \prod_p \left(1 - \left(\frac{1 - p^{i x}}{p - 1}\right)^2\right)\right).
\end{equation}

An alternative route to the 2-point correlation was given in \cite{MR2325314}. Assuming the Riemann hypothesis Conrey and Snaith made the conjecture that
\begin{align}
& \sum_{m, n \geq 1} f(\gamma(m) - \gamma(n)) 1_{0 < \gamma(m), \gamma(n) < T} \nonumber \\
= & \int_{0 \leq t \leq T, - T \leq x \leq T} \left(d^{mean}(t) \delta(x) + R_2(x, t)\right) f(x) \, dx dt \nonumber \\
& + O\left(T^{1 / 2 + \epsilon}\right)
\end{align}
treating the integral as principal value near $x = 0,$ for some $\epsilon > 0$ and holomorphic function in the strip $|\Im(s)| < 2,$ where $f(y) \in \mathbb{R} \textrm{ for } y \in \mathbb{R} \textrm{ with } f(- y) = f(y)$ and $f(y) \ll (1 + y^2)^{- 1} \textrm{ as } y \to \infty.$
The diagonal and off-diagonal contributions to the 2-point correlation are then given by
\begin{equation}
R_2^d(x, T) = \frac{1}{2 \pi^2} \left(\left(\frac{\zeta'}{\zeta}\right)'(1 + i x) - \sum_p \left(\frac{\log p}{p^{1 + i x} - 1}\right)^2\right)
\end{equation}
and
\begin{equation}
R_2^o(x, T) = \frac{1}{2 \pi^2} \textrm{exp}(- 2 \pi i x d^{mean}(T)) |\zeta(1 + i x)|^2 \prod_p \frac{\left(1 - \frac{1}{p^{1 + i x}}\right) \left(1 - \frac{2}{p} + \frac{1}{p^{1 + i x}}\right)}{\left(1 - \frac{1}{p}\right)^2}.
\end{equation}
This was found by first forming an expression for
\begin{equation} \label{Ratios of zeta functions with shifts}
 \int_{0 \leq t \leq T} \frac{\zeta\left(\frac{1}{2} + \alpha + i t\right) \zeta\left(\frac{1}{2} + \beta - i t\right)}{\zeta\left(\frac{1}{2} + \gamma + i t\right) \zeta\left(\frac{1}{2} + \delta - i t\right)} \, dt
\end{equation}
and then differentiating with respect to $2$ of the complex shifts along the critical line so as to utilise the residue theorem. Integrals over ratios of zeta functions with shifts, \eqref{Ratios of zeta functions with shifts}, were first considered by Farmer in \cite{zbMATH00436504} and subsequently evaluated by Conrey, Farmer and Zirnbauer in \cite{MR2482944} by replacing the functions on the numerator by
\begin{equation}
\sum_{n \geq 1} \left(\frac{1}{n^s} + X(s) \frac{1}{n^{1 - s}}\right)
\end{equation}
and those on the denominator by
\begin{equation}
\left(\sum_{n \geq 1} \frac{\mu(n)}{n^s}\right)^{- 1},
\end{equation}
where $s$ is the argument of the function and $\mu$ is the M$\ddot{\textrm{o}}$bius function, a multiplicative function given by
\begin{equation}
\mu(n) = \begin{dcases*} (- 1)^k & if $n \textrm{ has } k \geq 0 \textrm{ distinct prime factors}$ \\ 0 & otherwise. \end{dcases*}
\end{equation}
Dropping the highly oscillatory terms when taking the integral, namely terms containing just $1$ factor of $X(s)$ and the off-diagonal terms, the integral over the ratios of Riemann zeta functions was then shown to yield
\begin{align}
& \int_{0 \leq t \leq T} \frac{\zeta(1 + \alpha + \beta) \zeta(1 + \gamma + \delta)}{\zeta(1 + \alpha + \delta) \zeta(1 + \beta + \gamma)} A(\alpha, \beta, \gamma, \delta) \, dt \nonumber \\
& + \int_{0 \leq t \leq T} \left(\frac{t}{2 \pi}\right)^{- \alpha - \beta} \frac{\zeta(1 - \alpha - \beta) \zeta(1 + \gamma + \delta)}{\zeta(1 - \beta + \delta) \zeta(1 - \alpha + \gamma)} A(- \beta, - \alpha, \gamma, \delta) \, dt
\end{align}
after applying this recipe, where 
\begin{equation}
A(\alpha, \beta, \gamma, \delta) = \prod_p \frac{\left(1 - \frac{1}{p^{1 + \gamma + \delta}}\right) \left(1 - \frac{1}{p^{1 + \beta + \gamma}} - \frac{1}{p^{1 + \alpha + \delta}} + \frac{1}{p^{1 + \gamma + \delta}}\right)}{\left(1 - \frac{1}{p^{1 + \beta + \gamma}}\right) \left(1 - \frac{1}{p^{1 + \alpha + \delta}}\right)},
\end{equation}
this procedure is followed in the next chapter for the Selberg class.

In analogy with \eqref{Universal diagonal contribution} and \eqref{Universal off-diagonal contribution}, rescaling gives
\begin{equation} \label{Universal diagonal contribution for zeta}
\tilde{R}_2^d(x) = \lim_{t \to \infty} \frac{R_2^d(x / d^{mean}(t), t)}{(d^{mean}(t))^2} = - \frac{1}{2 \pi^2 x^2}
\end{equation}
and similarly
\begin{equation} \label{Universal off-diagonal contribution for zeta}
\tilde{R}_2^o(x) = \frac{\cos 2 \pi x}{2 \pi^2 x^2},
\end{equation}
together giving the rescaled $2$-point correlation as
\begin{equation} \label{Montgomery}
\tilde{R}_2(x) = 1 - \left(\frac{\sin \pi x}{\pi x}\right)^2.
\end{equation}
In \cite{MR1246830} this universal expression was shown by Keating to follow from knowledge of the asymptotic behaviour of the function in the Hardy-Littlewood conjecture, \eqref{Hardy-Littlewood conjecture}, also deduced in \cite{MR1246830}. Asymptotically as $|r| \to \infty$ the typical behaviour is given by $1 - \frac{1}{2 |r|},$ which is to be interpreted as
\begin{equation}
\frac{1}{2 |r|} \sum_{|R| \leq |r|} \lim_{N \to \infty} \frac{1}{N} \sum_{n \leq N} \Lambda(n + R) \Lambda(n) \sim 1 - \frac{\log |r|}{2 |r|} \textrm{ as } |r| \to \infty.
\end{equation}
Considering analogues of the action correlation statistics, \eqref{Action correlations} and \eqref{Rescaled action correlations} for the Riemann zeta function the $- \frac{1}{2 |r|}$ term in the averaged Hardy-Littlewood conjecture was recovered in \cite{MR1253066} by assuming Montgomery's conjecture, \eqref{Montgomery}, specifically by inverting the off-diagonal contribution given by \eqref{Universal off-diagonal contribution for zeta}. In the next chapter \eqref{Montgomery} is found to be the rescaled $2$-point correlation for the Selberg class. This universal behaviour of the $2$-point correlation is then used to generalise the inversion calculation performed in \cite{MR1253066}.
\subsubsection{Dirichlet $L$-functions}
Named in honour of Dirichlet, the series defining these $L$-functions generalises that of the Riemann zeta function,
\begin{equation} \label{Dirichlet series}
L(s, \chi_q) = \sum_{n \geq 1} \frac{\chi_q(n)}{n^s} \textrm{ for } \Re(s) > 1,
\end{equation}
where the coefficients in the numerator are multiplicative characters,
\begin{equation}
\chi_q(m n) = \chi_q(m) \chi_q(n) \, \forall m \textrm{ and } n \in \mathbb{Z},
\end{equation}
with periodicity $q,$
\begin{equation}
\chi_q(n + q) = \chi_q(n) \, \forall n \in \mathbb{Z}.
\end{equation}
It is stipulated that $\chi(1) = 1$ so as to avoid characters from being identically equal to zero, however, if the greatest common divisor of $n$ and $q,$ $(n, q),$ is not $1$ then the Dirichlet character is equal to zero, and conversely $\chi_q(n) \neq 0 \Rightarrow (n, q) = 1.$ 

The Riemann zeta function is generated from the principal character,
\begin{equation}
\chi_q^p(n) = \begin{dcases*} 1 & if $(n, q) = 1$ \\ 0 & otherwise, \end{dcases*}
\end{equation}
where
\begin{equation}
L(s; \chi_q^p) = \prod_{p | q} \left(1 - \frac{1}{p^s}\right) \zeta(s),
\end{equation}
which follows from the Euler product representation
\begin{equation}
L(s; \chi_q) = \prod_p \left(1 - \frac{\chi_q(p)}{p^s}\right)^{- 1} \textrm{ for } \Re(s) > 1.
\end{equation}

Given a Dirichlet character of period $q | Q,$ the multiplicative function, 
$$\begin{dcases*} \chi_q(n) & if $(n, Q) = 1$ \\ 0 & otherwise \end{dcases*},$$
is an induced Dirichlet character. If a Dirichlet character is primitive, that is to say it is not induced from another Dirichlet character, and if it is also non-principal, then \eqref{Dirichlet series} may be continued to an analytic function on the whole complex plane.

Furthermore, there is a functional equation,
\begin{equation}
L(s, \chi_q) = \frac{i^{a(\chi_q)} q^{1 / 2}}{\tau(\chi_q^*)} \left(\frac{\pi}{q}\right)^{s - 1 / 2} \frac{\Gamma((1 - s + a(\chi_q)) / 2)}{\Gamma((s + a(\chi_q)) / 2)} L(1 - s, \chi_q^*),
\end{equation}
where $\tau(\chi_q)$ is the Gauss sum,
\begin{equation}
\sum_{1 \leq n \leq q} \chi_q(n) \exp(2 \pi i n / q),
\end{equation}
and
\begin{equation}
a(\chi_q) = \begin{dcases*} 0 & if $\chi_q(- 1) = 1$ \\ 1 & otherwise. \end{dcases*}
\end{equation}

From the Euler product together with the multiplicative property of the characters,
\begin{equation}
\frac{L'}{L}(s, \chi_q) = - \sum_{n \geq 1} \frac{\chi_q(n) \Lambda(n)}{n^s} \textrm{ for } \Re(s) > 1,
\end{equation}
giving the generalised von Mangoldt function associated to a Dirichlet $L$-function,
\begin{equation}
\Lambda(n, \chi_q) = \chi_q(n) \Lambda(n).
\end{equation}
The mean and oscillating contribution to the density of zeros of Dirichlet $L$-functions are then given by
\begin{equation} \label{Mean density for Dirichlet}
d^{mean}(t; \chi_q) = \frac{1}{2 \pi} \log \frac{q |t|}{2 \pi}
\end{equation}
and
\begin{equation} \label{Oscillating part of density for Dirichlet}
d^{osc}(t; \chi_q) = - \frac{1}{2 \pi} \sum_{n \geq 1} \frac{\Lambda(n, \chi_q)}{n^{1 / 2}} \exp(i t \log n) + cc
\end{equation}
respectively.

In \cite{MR3030176} Bogomolny and Keating used the Hardy-Littlewood conjecture applied to primes in arithmetic progressions to deduce the $2$-point correlation,
\begin{equation} \label{2-point correlation for Dirichlet}
R_2(x, T; \chi_q) = \left<d(t + x_1; \chi_q) d(t + x_2; \chi_q)\right>_T,
\end{equation}
showing that
\begin{equation}
R_2^d(x, T; \chi_q) = - \frac{1}{4 \pi^2} \frac{\partial^2}{\partial x^2} \log(|\zeta(1 + i x)|^2 \Phi^d(x) \Psi(x, q)),
\end{equation}
where
\begin{equation}
\Phi^d(x) = \prod_p \exp\left(2 \sum_{n \geq 1} \frac{1 - n}{n^2 p^n} \cos(n x \log p)\right)
\end{equation}
and
\begin{equation}
\Psi(x, q) = \exp\left(- \sum_{p | q, n \geq 1} \frac{1}{n^2 p^n} \exp(i n x \log p) + cc \right).
\end{equation}
Also
\begin{equation}
R_2^o(x, T; \chi_q) = \frac{1}{4 \pi^2} \exp(i x \log(T q / (2 \pi)) |\zeta(1 + i x)|^2 \Phi^o(x) \Psi^o(x, q) + cc,
\end{equation}
where
\begin{equation}
\Phi^o(x) = \prod_p \left(1 - \frac{(p^{i x} - 1)^2}{(p - 1)^2}\right)
\end{equation}
and
\begin{equation}
\Psi^o(x, q) = \prod_{p | q} \left(1 + \frac{(p^{i x / 2} - p^{- i x / 2})^2}{p - p^{- i x}}\right)^{- 1}.
\end{equation}
In the rescaled limit this recovers the universal form \eqref{Montgomery}, which was deduced earlier by Bogomolny and Leb{\oe}uf in \cite{MR1284685} by using the method in \cite{MR1246830}, applying the averaged form of the Hardy-Littlewood conjecture for primes in arithmetic progressions.

From the functional equation it may be seen that Dirichlet $L$-functions are degree $1$ $L$-functions, examples of degree $2$ $L$-functions will now be given.
\subsubsection{Degree $2:$ $L$-functions associated to the Ramanujan tau function and elliptic curves}
The Ramanujan tau $L$-function is constructed from a modular form of weight 12, denoted by $\Delta,$
\begin{equation}
\Delta\left(\frac{a s + b}{c s + d}\right) = (c s + d)^{12} \Delta(s)
\end{equation}
for all integers $a, b, c \textrm{ and } d \textrm{ with } a d - b c = 1.$ The Ramanujan tau values are defined implicitly through
\begin{equation}
x \prod_{n \geq 1} (1 - x^n)^{24} = \sum_{n \geq 1} \tau(n) x^n
\end{equation}
and are the Fourier coefficients of
\begin{equation}
\Delta(z) = \sum_{n \geq 1} \tau(n) \exp(2 \pi i n z).
\end{equation}
The $L$-function associated to the Ramanujan tau function, for $\Re(s) > 1,$ is given by
\begin{equation}
L(s; \Delta) = \sum_{n \geq 1} \frac{\lambda(n; \Delta)}{n^s},
\end{equation}
where $\lambda(n; \Delta)$ is the $n$th normalised coefficient, 
\begin{equation}
\lambda(n; \Delta) = \frac{\tau(n)}{n^{11 / 2}}.
\end{equation}
For $\Re(s) > 1$ there is also the Euler product
\begin{equation}
L(s; \Delta) = \prod_p \left(1 - \frac{\lambda(p; \Delta)}{p^s} + \frac{1}{p^{2 s}}\right)^{- 1},
\end{equation}
which may be written in the form
\begin{equation}
\prod_p \left(1 - \frac{\alpha(p; \Delta)}{p^s}\right)^{- 1} \left(1 - \frac{\alpha^*(p; \Delta)}{p^s}\right)^{- 1}
\end{equation}
with $\alpha(p; \Delta) = \exp(i \theta(p; \Delta)) \textrm{ for real } \theta(p; \Delta).$ There is a functional equation,
\begin{equation}
L(s, \Delta) = \pi^{2 s - 1} \frac{\Gamma((1 - s) / 2 + 11 / 4) \Gamma((1 - s) / 2 + 13 / 4)}{\Gamma(s / 2 + 11 / 4) \Gamma(s / 2 + 13 / 4)} L(1 - s, \Delta),
\end{equation}
which shows that this is a degree $2$ $L$-function.

Another example of a degree $2$ $L$-function is given by associating an $L$-function to an elliptic curve,
\begin{equation}
E: y^2 + a_1 x y + a_3 y = x^3 + a_2 x^2 + a_4 x + a_6,
\end{equation}
which may be specified more succinctly by specifying its Weierstrass coefficients, $(a_1, a_2, a_3, a_4, a_6).$ Associated to such an elliptic curve is a modular form, the normalised coefficients of which, $\lambda(n; E) \textrm{ for } n = 1, 2, \ldots,$ generate the $L$-function,
\begin{equation}
L(s; E) = \sum_{n \geq 1} \frac{\lambda(n; E)}{n^s},
\end{equation}
for $\Re(s) > 1.$
As for the Ramanujan tau case there is an Euler product which may be written in the form
\begin{equation}
L(s; E) = \prod_p \left(1 - \frac{\alpha(p; E)}{p^s}\right)^{- 1} \left(1 - \frac{\beta(p; E)}{p^s}\right)^{- 1},
\end{equation}
where
\begin{equation}
\alpha(p; E) = \begin{dcases*} \exp(i \theta(p; E)) & if $p \nmid N$ \\ \lambda(p; E) & otherwise, \end{dcases*}
\end{equation}
and
\begin{equation}
\beta(p; E) = \begin{dcases*} \exp(- i \theta(p; E)) & if $p \nmid N$ \\ 0 & otherwise, \end{dcases*}
\end{equation}
with
\begin{equation}
\lambda(p; E) = 2 \cos(\theta(p; E)),
\end{equation}
where $N$ is the conductor of the elliptic curve and $\theta(p; E)$ a real parameter. The functional equation for these degree $2$ $L$-functions is given by 
\begin{equation}
L(s, E) = \left(\frac{\pi}{N^{1 / 2}}\right)^{2 s - 1} \frac{\Gamma((1 - s) / 2 + 1 / 4) \Gamma((1 - s) / 2 + 3 / 4)}{\Gamma(s / 2 + 1 / 4) \Gamma(s / 2 + 3 / 4)} L(s, E).
\end{equation}
\section{Outline of results in this thesis}
Chapter $2$ contains a heuristic study of the correlations of non-trivial zeros of functions in the Selberg class and the consequences for correlations of generalised von Mangoldt functions in the limiting universal case, and specifically for the von Mangoldt function in the non-universal case for the Riemann zeta function:
\subsection{Heuristic $1$}
The $2$-point correlation for non-trivial zeros of functions in the Selberg class,
$$
R_2(x, T; F) = \left<d(t + x_1; F) d(t + x_2; F)\right>_T,
$$
is given by
$$
R_2(x, T; F) = R_2^d(x, T; F) + R_2^o(x, T; F),
$$
with
$$
R_2^d(x, T; F) = \frac{1}{2 \pi^2} \left(\frac{(F \otimes \overline{F})'}{(F \otimes \overline{F})}\right)'(1 + i x) - B(i x; F)
$$
and
\begin{align}
R_2^o(x, T; F) = & \frac{1}{2 \pi^2 (r(F \otimes \overline{F})^2} (F \otimes \overline{F})(1 + i x) (F \otimes \overline{F})(1 - i x) \nonumber \\
& \exp\left(- 2 \pi i x d^{mean}(T; F)\right) A(i x; F), \nonumber
\end{align}
where $r(F \otimes \overline{F})$ is the residue of the pole of $F \otimes \overline{F}$ at $s = 1,$
\begin{align}
B(i x; F) = & \sum_p (\log p)^2 \left(- \sum_{h + m = k + n} \frac{a(p^m; F) a^*(p^n; F) \mu(p^h; F) \mu^*(p^k; F) m n}{p^{(n + k) (1 + i x)}} \right. \nonumber \\
& \left. \vphantom{- \sum_{h + m = k + n} \frac{a(p^m; F) a^*(p^n; F) \mu(p^h; F) \mu^*(p^k; F) m n}{p^{(n + k) (1 + i x)}}} + \sum_{l \geq 1} \frac{l^3 |b(p^l; F)|^2}{p^{l (1 + i x)}}\right) \nonumber
\end{align}
and
\begin{align}
A(i x; F) = & \prod_p \left(\sum_{h + m = k + n} \frac{a(p^m; F) a^*(p^n; F) \mu(p^h; F) \mu^*(p^k; F)}{p^{- r m + n + (1 + i x) k}}\right) \nonumber \\
& \exp\left(\sum_{l \geq 1} l |b(p^l; F)|^2 \left(\frac{2}{p^l} - \frac{1}{p^{l (1 - i x)}} - \frac{1}{p^{l (1 + i x)}}\right)\right). \nonumber
\end{align}
\subsection{Heuristic $2$}
The smoothed correlation statistic,
$$
\tilde{C}(r; F) = \int (d^{mean}(T))^2\tilde{R}_2^o(x d^{mean}(T); F) \exp(i r T) \, dx dT,
$$
where
$$
\tilde{R}_2^o(x; F) = \lim_{T \to \infty} \frac{R_2^o(x / d^{mean}(T; F), T; F)}{(d^{mean}(T; F))^2},
$$
is given asymptotically as $|r| \to \infty$ by
$$
- \frac{\textrm{deg}(F)}{\pi r} \textrm{Si}(r).
$$
\subsection{Heuristic $3$}
The correlation statistic for the Riemann zeta function case,
$$
C(r) = \int R_2^o(x, T) \exp(i r T) \, dx dT,
$$
is given by the Hardy-Littlewood conjecture,
$$
\begin{dcases*} 2 \prod_p \left(1 - \frac{1}{(p - 1)^2}\right) \prod_{q | r} \frac{q - 1}{q - 2} & if $r$ is even \\ 0 & otherwise, \end{dcases*}
$$
where $p \textrm{ and } q \textrm{ are both primes greater than } 2.$

In chapter $3$ the implications of the distribution of non-trivial zeros of functions in the Selberg class on the variance of generalised summatory von Mangoldt functions over short intervals are given:
\subsection{Theorem $1$}
If $1 < A_1 < A_2 < \textrm{deg}(F)$ and
$$
\mathcal{F}(X, T; F) = \frac{T \log X}{\pi} + O \left( T^{1 - c} \right)
$$
uniformly for $T^{A_1} \ll X \ll T^{A_2}$ for some $c > 0,$ then for any fixed $1 / A_2 < B_1 \leq B_2 < 1 / A_1$
\begin{align}
& Var^{mul}(X, \delta; F) = \frac{1}{6} \delta X (3 \log X - 4 \log 2) + O(\delta^{1 + c / 2} X) \nonumber \\
& + O_{\epsilon}(\delta^{1 - \epsilon} X (\delta X^{1 / A_1})^{- 2 A_1 / (4 A_1 + 1)}) + O_{\epsilon}(\delta^{1 - \epsilon} X (\delta X^{1 / A_2})^{1 / 2}) \nonumber
\end{align}
uniformly for $X^{- B_2} \ll \delta \ll X^{- B_1}.$
\subsection{Theorem $2$}
If $\textrm{deg}(F) < A_1 < A_2$ and
$$
\mathcal{F}(X, T; F) \sim \frac{T}{\pi} \left(\textrm{deg}(F) \log \frac{T}{2 \pi} + \log \mathfrak{q}(F) - \textrm{deg}(F)\right) + O\left(T^{1 - c}\right)
$$
uniformly for $T^{A_1} \ll X \ll T^{A_2}$ for some $c > 0,$ then for any fixed $1 / A_2 < B_1 \leq B_2 < 1 / A_1$
\begin{align}
& Var^{mul}(X, \delta; F) \nonumber \\
= & \frac{1}{2} \delta (\textrm{deg}(F) \log \frac{1}{\delta} + \log \mathfrak{q}(F) + (1 - \gamma - \log 2 \pi) \textrm{deg}(F)) \nonumber \\
& + O(\delta^{1 + c / 2} X) + O_{\epsilon}(\delta^{1 - \epsilon} X (\delta X^{1 / A_2})^{1 / 2}) + O_{\epsilon}(\delta^{1 - \epsilon} X (\delta X^{1 / A_2})^{- 2 A_1 / (4 A_1 + 1)}) \nonumber
\end{align}
uniformly for $X^{- B_2} \ll \delta \ll X^{- B_1}.$
\subsection{Theorem $3$}
If $1 / \textrm{deg}(F) < B_1 < B_2 \leq B_3 < 1$ and
$$
Var^{mul}(X, \delta; F) = \frac{1}{6} \delta X (3 \log X - 4 \log 2) + O(\delta^{1 + c} X)
$$
uniformly for $X^{- B_3} \ll \delta \ll X^{- B_1}$ for some $c > 0,$ then
\begin{align}
& Var^{fix}(X, h; F) = \frac{1}{6} h (6 \log X - (3 + 8 \log 2)) \nonumber \\
& + O_{\epsilon}(h X^{\epsilon} (h / X)^{c / 3}) + O_{\epsilon}(h X^{\epsilon} (h X^{- (1 - B_1)})^{1 / (3 (1 - B_1))}) \nonumber
\end{align}
uniformly for $X^{1 - B_3} \ll h \ll X^{1 - B_2}.$
\subsection{Theorem $4$}
If $0 < B_1 < B_2 \leq B_3 < 1 / \textrm{deg}(F)$ and
\begin{align}
Var^{mul}(X, \delta; F) = & \frac{1}{2} \delta X \left(\textrm{deg}(F) \log \frac{1}{\delta} - \textrm{deg}(F) (\log 2 \pi + \gamma - 1) + \log \mathfrak{q}(F)\right) \nonumber \\
& + O \left(\delta^{1+c} X \right) \nonumber
\end{align}
uniformly for $X^{- B_3} \ll \delta \ll X^{- B_1}$ for some $c > 0,$ then
\begin{align}
Var^{fix}(X, h ; F) = & h \left(\textrm{deg}(F) \log \frac{X}{h} - \textrm{deg}(F) (\log 2 \pi + \gamma) + \log \mathfrak{q}(F)\right) \nonumber \\ 
& + O_\epsilon\left(h X^{\epsilon} (h / X)^{c /3} \right) + O_\epsilon \left( h X^{1 - (B_2 - B_1)/3 (1 - B_1) + \epsilon}\right) \nonumber
\end{align}
uniformly for $X^{1 - B_3} \ll h \ll X^{1 - B_2}.$
\subsection{Remarks}
Theorem $1$ and Theorem $3$ are valid only when the degree of the function in the Selberg class considered is greater than $1.$ When the degree is equal to $1$ the behaviour of the variance first studied by Goldston and Montgomery in \cite{MR1018376} and later refined by Chan in \cite{MR2009438} is recovered, as given by Theorem $2$ and Theorem $4.$ The degree $1$ case is therefore a special case of the more general variance which does not exhibit the change in behaviour for $h$ of order $X^{1 - 1 / \textrm{deg}(F)}.$ Data exhibiting this phenomenom is plotted at the end of chapter $3$ and the relation to Heuristic $2$ and Heuristic $3$ is discussed in the concluding chapter.

\chapter{Zero statistics and arithmetic correlations associated to $L$-functions}
In this chapter the $2$-point correlation for the Selberg class is calculated by using the method described by Conrey and Snaith in \cite{MR2325314}. The limiting form of this statistic is then used to generalise the work in \cite{MR1253066} from the case of the Riemann zeta function to the whole Selberg class and then specifically for the Riemann zeta function case the heuristics are applied to recover the Hardy-Littlewood conjecture by inverting the $2$-point correlation for this case.
\section[The $2$-point correlation statistic for the Selberg class]{The $2$-point correlation statistic \\ for the Selberg class}  
Substituting the density of Selberg class zeros, \eqref{Mean density} and \eqref{Oscillating part of density}, into the $2$-point correlation,
\begin{equation}
R_2(x, T; F) = \left<d(t + x_1; F) d(t + x_2; F)\right>_T,
\end{equation}
gives the analogue of the expression for the $2$-point correlation for the system $D,$ \eqref{Terms in 2-point correlation for D}, now applied to the Selberg class for a generic function $F$ contained therein,
\begin{equation}
R_2(x, T; F) = (d^{mean}(T; F))^2 + R_2^d(x, T; F) + R_2^o(x, T; F),
\end{equation}
where
\begin{equation}
R_2^d(x, T; F) = \frac{1}{2 \pi^2} \sum_{n \geq 1} \frac{|\Lambda(n; F)|^2}{n} \cos (x \log n)
\end{equation}
and
\begin{equation}
R_2^o(x, T; F) = \frac{1}{4 \pi^2} \sum_{n \geq 1, \, r \neq 0} \frac{\Lambda(n + r; F) \Lambda^*(n; F)}{n} \exp(i (x \log n - T r / n)) + cc.
\end{equation}
The diagonal and off-diagonal contributions are now indirectly calculated by following the recipe described in the introduction for the Riemann zeta function case. 

The integral over the ratios of functions in the Selberg class,
\begin{equation}
\int_{|t| \leq T} \frac{F\left(\frac{1}{2} + \alpha + i t\right) \overline{F}\left(\frac{1}{2} + \beta - i t\right)}{F\left(\frac{1}{2} + \gamma + i t\right) \overline{F}\left(\frac{1}{2} + \delta - i t\right)} \; dt,
\end{equation}
is evaluated by replacing functions in the numerator by
\begin{equation}
\sum_{n \geq 1} \left(\frac{a(n; F)}{n^s} + X(s; F) \frac{a^*(n; F)}{n^{1 - s}}\right)
\end{equation}
and those in the denominator by
\begin{equation}
\left(\sum_{n \geq 1} \frac{\mu(n; F)}{n^s}\right)^{- 1},
\end{equation}
where $s$ is the argument of the function. The values of the multiplicative function $\mu(n; F)$ are defined implicitly as the coefficients for the reciprocal of $F(s),$
\begin{equation}
\frac{1}{F}(s) = \sum_{n \geq 1} \frac{\mu(n; F)}{n^s}.
\end{equation}

Making the substitutions results in
\begin{align}
& \int_{|t| \leq T} \frac{(F \otimes \overline{F}) (1 + \alpha + \beta) (F \otimes \overline{F}) (1 + \gamma + \delta)}{(F \otimes \overline{F}) (1 + \alpha + \delta) (F \otimes \overline{F}) (1 + \beta + \gamma)} A(\alpha, \beta, \gamma, \delta; F) \, dt \nonumber \\
& + \int_{|t| \leq T} \left(\mathfrak{q}(F) \left(\frac{|t|}{2 \pi}\right)^{\textrm{deg}(F)}\right)^{- (\alpha + \beta)} \frac{(F \otimes \overline{F}) (1 - \beta - \alpha) (F \otimes \overline{F}) (1 + \gamma + \delta)}{(F \otimes \overline{F}) (1 -\beta + \delta) (F \otimes \overline{F}) (1 - \alpha + \gamma)} \nonumber \\
& A(- \beta, - \alpha, \gamma, \delta; F) \, dt, \label{Ratios conjecture}
\end{align}
where
\begin{align}
& A(\alpha, \beta, \gamma, \delta; F) = \prod_p \left(\sum_{h + m = k + n} \frac{a(p^m; F) a^*(p^n; F) \mu(p^h; F) \mu^*(p^k; F)}{p^{(1 / 2 + \alpha) m + (1 / 2 + \beta) n + (1 / 2 + \gamma) h + (1 / 2 + \delta) k}}\right) \nonumber \\
& \exp\left(\sum_{l \geq 1} l |b(p^l; F)|^2 \left(\frac{1}{p^{l (1 + \alpha + \delta)}} + \frac{1}{p^{l (1 + \beta + \gamma)}} - \frac{1}{p^{l (1 + \alpha + \beta)}} - \frac{1}{p^{l (1 + \gamma + \delta)}}\right)\right).
\end{align}
The factor appearing in the second term of the integral follows from the identity
\begin{equation}
X(1 / 2 + i t + \alpha; F) \overline{X} (1 / 2 - i t + \beta; F) = \left(\mathfrak{q}(F) \left(\frac{|t|}{2 \pi}\right)^{\textrm{deg}(F)}\right)^{- (\alpha + \beta)} + O\left(\frac{1}{|t|}\right).
\end{equation}
Differentiating \eqref{Ratios conjecture} with respect to $\alpha$ and $\beta$ and then setting $\gamma = \alpha$ and $\delta = \beta$ gives an expression for a product of shifted logarithmic derivatives,
\begin{align}
& \int_{t \leq |T|} \frac{F'}{F}(1 / 2 + i t + \alpha) \frac{\overline{F}'}{\overline{F}}(1 / 2 - i t + \beta) \, dt \nonumber \\
\sim & \int_{t \leq |T|} \left(\frac{(F \otimes \overline{F})'}{(F \otimes \overline{F})}\right)'(1 + \alpha + \beta) \, dt \nonumber \\
& + \int_{t \leq |T|} \frac{1}{(r(F \otimes \overline{F}))^2} \left(\mathfrak{q}(F) \left(\frac{|t|}{2 \pi}\right)^{\textrm{deg}(F)}\right)^{- (\alpha + \beta)} A(- \beta, - \alpha, \alpha, \beta; F) \nonumber \\
& (F \otimes \overline{F})(1 - \alpha - \beta) (F \otimes \overline{F})(1 + \alpha + \beta) \, dt \nonumber \\
& + \int_{t \leq |T|} \frac{\partial^2}{\partial \alpha \partial \beta} A(\alpha, \beta, \gamma, \delta)|_{\gamma = \alpha, \delta = \beta} \, dt, \label{Derivative of ratios conjecture}
\end{align}
where $r(F \otimes \overline{F})$ is the residue of the simple pole of $F \otimes \overline{F}$ at $s = 1,$ which enters upon applying the identity
\begin{equation}
\frac{\partial}{\partial \alpha} \frac{g(\alpha, \gamma)}{(F \otimes \overline{F})(1 - \alpha + \gamma)}|_{\gamma = \alpha} = \frac{g(\alpha, \alpha)}{r(F \otimes \overline{F})}
\end{equation}
for a function, $g,$ analytic at $(\alpha, \alpha).$

As in the introduction, a function $f$ is taken which is holomorphic in the strip $|\Im(s)| < 2,$ where $f(y) \in \mathbb{R} \textrm{ for } y \in \mathbb{R} \textrm{ with } f(- y) = f(y)$ and $f(y) \ll (1 + y^2)^{- 1} \textrm{ as } y \to \infty.$ Assuming the generalisation of Riemann's hypothesis for the Selberg class the sum
\begin{equation}
\sum_{m, n \geq 1} f(t(m) - t(n)) 1_{|t(m)|, |t(n)| \leq T}
\end{equation}
is now evaluated by means of the residue theorem,
\begin{equation}
\sum_{m, n \geq 1} f(t(m) - t(n)) 1_{|t(m)|, |t(n)| \leq T} = \frac{1}{(2 \pi i)^2} \int_{\mathcal{C}, \mathcal{C}} \frac{F'}{F}(u) \frac{F'}{F}(v) f(- i (u - v)) \, du dv,
\end{equation}
where $1 / 2 < a <1$ and $\mathcal{C}$ is a positively oriented rectangle with vertices at $1 - a - i T$, $a -i T$, $a + i T$ and $1 - a + i T.$

The integrals with contours passing through $1 - a$ may be reflected by use of the functional equation, in particular taking the logarithmic derivative,
\begin{equation}
\frac{F'}{F}(s) = \frac{X'}{X}(s; F) - \frac{\overline{F}'}{\overline{F}}(1 - s),
\end{equation}
where
\begin{align}
\frac{X'}{X}(s; F) = & - 2 \log Q - \sum_{1 \leq j \leq r} \lambda_j \left(\frac{\Gamma'}{\Gamma}(\lambda_j s + \mu_j) + \frac{\Gamma'}{\Gamma}(\lambda_j (1 - s) + \mu_j^*)\right) \nonumber \\
= & - 2 \pi d^{mean}(T; F) + O\left(\frac{1}{|t|}\right).
\end{align}

The integral with the contour passing through $a$ and $1 - a$ is given by
\begin{align}
& - \frac{1}{(2\pi i)^2} \int_{a - i T}^{a + i T} \int_{1 - a - i T}^{1 - a + i T} \frac{F'}{F}(u) \frac{F'}{F}(v) f(- i (u - v)) \, du dv \nonumber \\
= & \frac{i}{(2 \pi )^2} \int_{- 2 T - i (1 - 2 a)}^{2 T - i (1 - 2 a)} f(\eta) \int_{a - i T_1}^{a + i T_2} \frac{F'}{F}(v) \frac{F'}{F}(v + i \eta) \, dv d\eta,
\end{align}
where $T_1 = \textrm{min} \{T, T + \Re(\eta)\} \textrm{ and } T_2 = \textrm{min} \{T, T - \Re(\eta)\}.$
Making use of the functional equation then gives
\begin{equation}
- \frac{i}{(2 \pi)^2} \int_{- 2 T - i (1 - 2 a)}^{2 T - i (1 - 2 a)} f(\eta) \int_{a - i T_1}^{a + i T_2} \frac{F'}{F}(v) \frac{\overline{F}'}{\overline{F}}(1 - v - i \eta) \, dv d\eta
\end{equation}
with a contribution of $O_{\epsilon}(T^{\epsilon})$ coming from the $\frac{X'}{X}$ term after moving the contour to the right of $1.$ The main contribution is therefore
\begin{equation}
\int_{- 2 T - i (1 - 2 a)}^{2 T - i (1 - 2 a)} f(\eta) \int_{- T_1 \leq t \leq T_2} \frac{F'}{F} (a + i t) \frac{\overline{F}'}{\overline{F}}(1 - a - i t - i \eta) \, dt d\eta.
\end{equation}
Using \eqref{Derivative of ratios conjecture}, the inner integrand then yields
\begin{align}
\left(\frac{(F \otimes \overline{F})'}{(F \otimes \overline{F})}\right)' (1 + \alpha + \beta) & \nonumber \\
+ \frac{1}{(r(F \otimes \overline{F}))^2} \left(\mathfrak{q}(F) \left(\frac{|t|}{2 \pi}\right)^{\textrm{deg}(F)}\right)^{- (\alpha + \beta)} A(- \beta, - \alpha, \alpha, \beta; F) & \nonumber \\
(F \otimes \overline{F})(1 - \alpha - \beta) (F \otimes \overline{F})(1 + \alpha + \beta), & \nonumber \\
+ \frac{\partial^2}{\partial \alpha \partial \beta} A(\alpha, \beta, \gamma, \delta; F)|_{\gamma = \alpha, \delta = \beta},
\end{align}
where $\alpha = a - 1 / 2$ and $\beta = 1 / 2 - a - i \eta.$ Setting these values then gives
\begin{equation}
A(\alpha, \beta, \gamma, \delta; F) = A(i \eta; F),
\end{equation}
where
\begin{align}
A(s; F) = & \prod_p \left(\sum_{h + m = k + n} \frac{a(p^m; F) a^*(p^n; F) \mu(p^h; F) \mu^*(p^k; F)}{p^{- s m + n + (1 + s) k}}\right) \nonumber \\
& \exp\left(\sum_{l \geq 1} l |b(p^l; F)|^2 \left(\frac{2}{p^l} - \frac{1}{p^{l (1 - s)}} - \frac{1}{p^{l (1 + s)}}\right)\right)
\end{align}
and
\begin{equation}
\frac{\partial^2}{\partial \alpha \partial \beta} A(\alpha, \beta, \gamma, \delta; F)|_{\gamma = \alpha, \delta = \beta} = - B(i \eta; F),
\end{equation}
with
\begin{align}
B(s; F) = & \sum_p (\log p)^2 \left(- \sum_{h + m = k + n} \frac{a(p^m; F) a^*(p^n; F) \mu(p^h; F) \mu^*(p^k; F) m n}{p^{(n + k) (1 + s)}} \right. \nonumber \\
& \left. \vphantom{- \sum_{h + m = k + n} \frac{a(p^m; F) a^*(p^n; F) \mu(p^h; F) \mu^*(p^k; F) m n}{p^{(n + k) (1 + s)}}} + \sum_{l \geq 1} \frac{l^3 |b(p^l; F)|^2}{p^{l (1 + s)}}\right).
\end{align}
Extending the inner region of integration results in
\begin{equation}
\frac{1}{(2 \pi)^2} \int_{|t| \leq T} \int_{- 2 T - i (1 - 2 a)}^{2 T - i (1 - 2 a)} f(\eta) g(\eta, t) \, d\eta dt,
\end{equation}
incurring an error of $O(T^{\epsilon}),$ where
\begin{align}
g(\eta, t) = & \left(\frac{(F \otimes \overline{F})'}{(F \otimes \overline{F})}\right)' (1 - i \eta) + \frac{1}{(r(F \otimes \overline{F}))^2} \left(\mathfrak{q}(F) \left(\frac{|t|}{2 \pi}\right)^{\textrm{deg}(F)}\right)^{i \eta} \nonumber \\
& A(- i \eta; F) (F \otimes \overline{F})(1 - i \eta) (F \otimes \overline{F})(1 + i \eta) - B(- i \eta; F).
\end{align}
As $\eta \to 0$
\begin{equation}
g(\eta, t) = \frac{2 \pi i}{\eta} d^{mean}(t; F) + O(1),
\end{equation}
moving the path of integration of the inner integral to the real line then passes through $0,$ giving the principal value
\begin{equation}
\frac{f(0)}{2} \int_{|t| \leq T} d^{mean}(t; F) \, dt + \frac{1}{(2 \pi)^2} \int_{|t| \leq T, |\eta| \leq 2 T} f(\eta) g(- \eta, t) \, d\eta dt.
\end{equation}

The integral with both contours passing through $1 - a$ proceeds similarly except the functional equation is used for both of the logarithmic derivatives appearing in the integrand. The main contribution is given by
\begin{equation}
\frac{1}{(2 \pi)^2} \int_{|v| \leq T, |u| \leq T} \frac{X'}{X}(1 / 2 + i u) \frac{X'}{X}(1 / 2 + i v) f(u - v) \, du dv,
\end{equation}
where an error of $O(T^{\epsilon})$ is incurred by moving contours involving the logarithmic derivative of $F$ to the right of $1$ and those involving the logarithmic derivative of $X$ to the $1 / 2$ line. Using the even property of $f$ then gives
\begin{align}
& 2 \int_{|v| \leq T, v \leq u \leq T} d^{mean}(u; F) d^{mean}(v; F) f(u - v) \, du dv \nonumber \\
= & 2 \int_{0 \leq \eta \leq 2 T, - T \leq t \leq T - \eta} d^{mean}(t + \eta; F) d^{mean}(t; F) \, dt d\eta.
\end{align}
Incurring an error term of $(\log T)^3$ the integral is given by
\begin{align}
& 2 \int_{0 \leq \eta \leq T, |t| \leq T} f(\eta)( d^{mean}(t; F))^2 \, dt d\eta \nonumber \\
= & \int_{|t| \leq T, |\eta| \leq 2 T} f(\eta) (d^{mean}(t; F))^2 \,  d\eta dt.
\end{align}
The remaining integral gives a contribution of $O(T^{\epsilon}),$ therefore in total the recipe gives the $2$-point correlation, $R_2(x, T; F) = R_2^d(x, T; F) + R_2^o(x, T; F),$ with
\begin{equation}
R_2^d(x, T; F) = \frac{1}{2 \pi^2} \left(\frac{(F \otimes \overline{F})'}{(F \otimes \overline{F})}\right)'(1 + i x) - B(i x; F)
\end{equation}
and
\begin{align}
R_2^o(x, T; F) = & \frac{1}{2 \pi^2 (r(F \otimes \overline{F})^2} (F \otimes \overline{F})(1 + i x) (F \otimes \overline{F})(1 - i x) \nonumber \\
& \exp\left(- 2 \pi i x d^{mean}(T; F)\right) A(i x; F).
\end{align}
For $x \to 0$
\begin{equation}
\left(\frac{(F\otimes \overline{F})'}{(F\otimes \overline{F})}\right)'(1 + i x) = - \frac{1}{x^2} + O(1),
\end{equation}
\begin{equation}
B(i x; F) = O(1),
\end{equation}
\begin{equation}
(F \otimes \overline{F})(1 + i x) (F \otimes \overline{F})(1 + i x) = \frac{(r(F \otimes \overline{F}))^2}{x^2} + O(1)
\end{equation}
and
\begin{equation}
A(i x; F) = 1 + O(x^2).
\end{equation}
Rescaling the $2$-point correlation then once more gives the universal expressions
\begin{equation}
\tilde{R}_2^d(x; F) = \lim_{T \to \infty} \frac{R_2^d(x / d^{mean}(T; F), T; F)}{(d^{mean}(T; F))^2} = - \frac{1}{2 \pi^2 x^2}
\end{equation}
and
\begin{equation}
\tilde{R}_2^o(x; F) = \frac{\cos 2 \pi x}{2 \pi^2 x^2}.
\end{equation}
\section{Inverting the limiting form of the $2$-point correlation statistic for the Selberg class}
Assuming for any $F$ there is a function $C(r; F), \textrm{ with } C(- r; F) = C^*(r; F),$ that captures the correlations, that is allows the off-diagonal contribution to be expressed as
\begin{equation}
R_2^o(x, T; F) = \frac{1}{4 \pi^2} \int_{n \geq 1} \frac{C(r; F)}{n} \exp(i (x \log n - T r / n)) \, dr dn + cc,
\end{equation}
then an expression for $C(r; F)$ may be determined by Fourier inversion.

\begin{align}
\int R_2^o(x, T; F) \, dx & = \frac{1}{2 \pi} \int_{n \geq 1} \frac{C(r; F)}{n} \delta(\log n) \exp(-  i T r / n) \, dr dn + cc \\ 
& = \frac{1}{2 \pi} \int_{y \geq 0} C(r; F) \delta(y) \exp(- i T r \exp(- y)) \, dr dy + cc \\
& = \frac{1}{4 \pi} \int C(r; F) \exp(- i T r) \, dr + cc \\
& = \frac{1}{2\pi} \int C(r; F) \exp(- i T r) \, dr.
\end{align}
Giving the Fourier transform 
\begin{equation}
\mathscr{F}[C(r; F)] = 2 \pi \int R_2^o(x, T; F) \, dx,
\end{equation}
and by Fourier inversion,
\begin{equation}
C(r; F) = \int R_2^o(x, T; F) \exp(i r T) \, dx dT. \label{Arithmetic correlation function}
\end{equation}
The mean density of zeros then necessarily enters the smoothed correlation function
\begin{equation}
\tilde{C}(r; F) = \int (d^{mean}(T))^2\tilde{R}_2^o(x d^{mean}(T); F) \exp(i r T) \, dx dT.
\end{equation}
The integral over $x$ is
\begin{align}
\frac{1}{2 \pi^2} \int \frac{1}{x^2} \cos(2 \pi x d^{mean}(T)) \, dx = & \frac{1}{\pi} \int \frac{1}{x^2} \cos(x d^{mean}(T)) \, dx \nonumber \\
= & \frac{1}{\pi} \mathscr{F}\left[\frac{1}{x^2}\right](d^{mean}(T)).
\end{align}
This may be deduced by considering the Fourier transform of the triangle function
\begin{equation}
T(x) = \begin{dcases*} 1 - |x| & if $|x| \leq 1$ \\ 0 & otherwise, \end{dcases*}
\end{equation}
noting that the triangle function may be expressed in terms of the sign function,
\begin{equation}
\textrm{sgn}(x) = \begin{dcases*} - 1 & if $x \leq 0$ \\ 1 & otherwise, \end{dcases*}
\end{equation}
by the linear relation
\begin{equation}
T(x) = \frac{1}{2} (1 - x) \textrm{sgn}(1 - x) + \frac{1}{2} (1 + x) \textrm{sgn}(1 + x) - x \textrm{sgn}(x).
\end{equation}
Then 
\begin{equation}
\mathscr{F}[T](k) = 2 \int_{0 \leq x \leq 1} (1 - x) \cos(k x) \, dx = \left(\textrm{sinc}\left(\frac{k}{2}\right)\right)^2,
\end{equation}
by integrating by parts, giving
\begin{align}
\int \left(\textrm{sinc}\left(\frac{k}{2}\right)\right)^2 \textrm{exp}(i k x) \, dk = & \pi (1 - x) \textrm{sgn}(1 - x) \nonumber \\
& +  \pi (1 + x) \textrm{sgn}(1 + x) \nonumber \\
& - 2 \pi x \textrm{sgn}(x)
\end{align}
by Fourier inversion. Expanding the sinc function in terms of exponential functions then gives
\begin{align}
\int \left(\textrm{sinc}\left(\frac{k}{2}\right)\right)^2 \textrm{exp}(i k x) \, dk = & - \mathscr{F}\left[\frac{1}{k^2}\right](1 - x) \nonumber \\
& - \mathscr{F}\left[\frac{1}{k^2}\right](1 + x) \nonumber \\
& + 2 \mathscr{F}\left[\frac{1}{k^2}\right](x),
\end{align}
giving the identification
\begin{equation}
\mathscr{F}\left[\frac{1}{x^2}\right](k) = - \pi k \textrm{sgn}(k).
\end{equation}
The smoothed correlation function is then given by another Fourier transform,
\begin{align}
\tilde{C}(r; F) & = - \mathscr{F} \left[d^{mean}(T) \textrm{sgn}(d^{mean}(T))\right](r) \nonumber \\
& \sim \frac{\textrm{deg}(F)}{\pi} \int_{0 \leq T \leq 1} \log T \cos (r T) \, dT \textrm{ as } |r| \to \infty \nonumber \\
& \sim - \frac{\textrm{deg}(F)}{\pi r} \textrm{Si}(r) \textrm{ as } |r| \to \infty.
\end{align}
\section{Inverting the $2$-point correlation statistic for the Riemann zeta function}
The Hardy-Littlewood conjecture is now shown to follow from the expression for the off-diagonal contribution to the $2$-point correlation for the case of the Riemann zeta function. 

The off-diagonal contribution to the $2$-point correlation is given by
\begin{align}
R_2^o(x, T) & = \frac{1}{2 \pi^2} \exp(- 2 \pi i x d^{mean}(T)) |\zeta(1 + i x)|^2 \prod_p \frac{\left(1 - \frac{1}{p^{1 + i x}}\right) \left(1 - \frac{2}{p} + \frac{1}{p^{1 + i x}}\right)}{\left(1 - \frac{1}{p}\right)^2} \nonumber \\
& = \frac{1}{2 \pi^2} \exp(- 2 \pi i x d^{mean}(T)) |\zeta(1 + i x)|^2 \prod_p \frac{\left(1 - \frac{1}{p^{1 + i x}}\right) (p^2 - 2 p + p^{1 - i x})}{(p - 1)^2}.
\end{align}
Using the Euler product on the $1$-line then gives
\begin{align}
& & \frac{1}{2 \pi^2} \exp(- 2 \pi i x d^{mean}(T)) \zeta(1 - i x) \prod_p \left(1 - \frac{1}{(p - 1)^2} + \frac{p^{1 - i x}}{(p - 1)^2}\right) \nonumber \\
& = & \frac{1}{2 \pi^2} \exp(- 2 \pi i x d^{mean}(T)) \zeta(1 - i x) \prod_p \left(1 + \frac{p^{1 - i x} \left(1 - \frac{1}{p^{1 - i x}}\right)}{(p - 1)^2}\right).
\end{align}
Introducing Euler's totient function,
\begin{equation}
\phi(n) = n \prod_{p | n} \left(1 - \frac{1}{p}\right),
\end{equation}
and the M$\ddot{\textrm{o}}$bius function, this can then be expressed as
\begin{equation}
\frac{1}{2 \pi^2} \exp(- 2 \pi i x d^{mean}(T)) \zeta(1 - i x) \prod_p \left(1 + \left(\frac{\mu(p)}{\phi(p)}\right)^2 p^{1 - i x} \left(1 + \frac{\mu(p)}{p^{1 - i x}}\right)\right),
\end{equation}
which, by multiplicativity of $\phi$ and $\mu$ yields
\begin{equation}
\frac{1}{2 \pi^2} \exp(- 2 \pi i x d^{mean}(T)) \zeta(1 - i x) \sum_{n \geq 1} \left(\frac{\mu(n)}{\phi(n)}\right)^2 n^{1 - i x} \prod_{p | n} \left(1 + \frac{\mu(p)}{p^{1 - i x}}\right),
\end{equation}
where the sum is over squarefree values of $n$ due to the presence of $\mu(n).$ Making use of the multiplicativity of $\mu$ once more gives
\begin{equation}
\frac{1}{2 \pi^2} \exp(- 2 \pi i x d^{mean}(T)) \zeta(1 - i x) \sum_{n \geq 1} \left(\frac{\mu(n)}{\phi(n)}\right)^2 n^{1 - i x} \sum_{d | n} \frac{\mu(d)}{d^{1 - i x}}.
\end{equation}
Using now the divergent series representation of the Riemann zeta function on the $1$-line results in
\begin{equation}
\frac{1}{2 \pi^2} \exp(- 2 \pi i x d^{mean}(T)) \sum_{n \geq 1} \left(\frac{\mu(n)}{\phi(n)}\right)^2 n^{1 - i x} \sum_{m \geq 1, d | n} \frac{\mu(d)}{(md)^{1 - i x}}.
\end{equation}
Replacing $m d$ by $l,$ the sum over $m$ and $d | n$ becomes a sum over $l,$ $d | n$ and $d | l,$ 
\begin{equation}
\sum_{l \geq 1, d | n, d | l} \frac{\mu(d)}{l^{1 - i x}} = \sum_{l \geq 1, d | (l, n)} \frac{\mu(d)}{l^{1 - i x}} = \sum_{l \geq 1, (l, n) = 1} \frac{1}{l^{1 - i x}},
\end{equation}
giving
\begin{equation}
\frac{1}{2 \pi^2} \exp(- 2 \pi i x d^{mean}(T)) \sum_{l, n \geq 1, (l, n) = 1} \left(\frac{\mu(n)}{\phi(n)}\right)^2 \left(\frac{n}{l}\right)^{1 - i x}.
\end{equation}
This may be approximated by
\begin{equation}
\frac{1}{2 \pi^2} \exp(- 2 \pi i x d^{mean}(T))  \sum_{l, n \geq 1, (l, n) = 1} \left(\frac{\mu(n)}{\phi(n)}\right)^2 \int_{0 \leq y \leq T / 2 \pi} \frac{1}{y^{1 - i x}} \delta(y - l / n) \, dy,
\end{equation}
where the resulting error is reduced by taking the value of $T$ to be large. Substituting the expression for the mean density and taking the exponential inside the integral then gives
\begin{equation} \label{Off-diagonal contribution for inversion}
\frac{1}{2 \pi^2} \sum_{n \geq 1, l \geq 1, (l, n) = 1} \left(\frac{\mu(n)}{\phi(n)}\right)^2 \int_{0 \leq z \leq 1} \frac{1}{z} \exp(i x \log z) \delta(t / (2 \pi) z- l / n) \, dz,
\end{equation}
after making a change of variables in the integral.

Substituting \eqref{Off-diagonal contribution for inversion} for $R_2^o(x, T)$ into \eqref{Arithmetic correlation function} then gives
\begin{align}
& \frac{1}{2 \pi^2} \sum_{l, n \geq 1, (l, n) = 1} \left(\frac{\mu(n)}{\phi(n)}\right)^2 \int \int_{0 \leq z \leq 1} \frac{1}{z} \exp\left(i x \log z\right) \delta(T / (2 \pi) z- l / n) \, dz \nonumber \\
& \exp(i r T) \, dx dT.
\end{align}
Integrating over $x$ then introduces a delta function,
\begin{align}
& \frac{1}{\pi} \sum_{l, n \geq 1, (l, n) = 1} \left(\frac{\mu(n)}{\phi(n)}\right)^2 \int \int_{0 \leq z \leq 1} \frac{1}{z} \delta(T / (2 \pi) z - l / n) \delta(\log z) \nonumber \\
& \exp(i r t) \, dz dT \nonumber \\
= & \frac{1}{\pi} \sum_{n \geq 1, l \geq 1, (l, n) = 1} \left(\frac{\mu(n)}{\phi(n)}\right)^2 \int \int_{w \leq 0} \delta(w) \delta(T / (2 \pi) \exp(w) - l / n) \nonumber \\
& \exp(i r T) \, dw dT.
\end{align}
Integrating over $w \leq 0$ captures half the mass of the delta function, giving
\begin{align}
& \frac{1}{2 \pi} \sum_{l, n \geq 1, (l, n) = 1} \left(\frac{\mu(n)}{\phi(n)}\right)^2 \int \delta(T / (2 \pi) - l / n) \exp(i r T) \, dT \nonumber \\
= & \sum_{l, n \geq 1, (l, n) = 1} \left(\frac{\mu(n)}{\phi(n)}\right)^2 \int \delta(u - l / n) \exp(2 \pi i h u) \, du \nonumber \\
= & \sum_{l, n \geq 1, (l, n) = 1} \left(\frac{\mu(n)}{\phi(n)}\right)^2 \exp((2 \pi i l r) / n),
\end{align}
which may be approximated by
\begin{equation}
\sum_{n \geq 1} \left(\frac{\mu(n)}{\phi(n)}\right)^2 c(r, n),
\end{equation}
where
\begin{equation}
c(r, n) = \sum_{1 \leq l \leq n, (l, n) = 1} \exp((2 \pi i l r) / n)
\end{equation}
is Ramanujan's sum.

Ramanujan sums are multiplicative,
\begin{equation}
c(l, m) c(l, n) = c(l, mn)
\end{equation}
if $(m, n) = 1.$ They also satisfy
\begin{equation}
c(n, p) = \begin{dcases*} p - 1 & if $p | n$ \\ - 1 & otherwise \end{dcases*}
\end{equation}
for all primes, $p.$ Using these properties then gives
\begin{align}
& \sum_{n \geq 1} \left(\frac{\mu(n)}{\phi(n)}\right)^2 c(r, n) \nonumber \\
= & \prod_p \left(1 + \left(\frac{\mu(p)}{\phi(p)}\right)^2 c(r, p)\right) \nonumber \\
= & \begin{dcases*} 2 \prod_{p > 2} \left(1 + \left(\frac{\mu(p)}{\phi(p)}\right)^2 c(r, p)\right) & if $r$ is even \\ 0 & otherwise, \end{dcases*}
\end{align}
where 
\begin{align}
& \prod_{p > 2} \left(1 + \left(\frac{\mu(p)}{\phi(p)}\right)^2 c(r, p)\right) \nonumber \\
= & \prod_{p \nmid r, p > 2} \left(1 - \frac{1}{(p - 1)^2}\right) \prod_{q | r, q > 2} \left(1 + \frac{1}{q - 1}\right) \nonumber \\
= & \prod_{p > 2} \left(1 - \frac{1}{(p - 1)^2}\right) \prod_{q | r, q > 2} \left(\frac{q - 1}{q - 2}\right).
\end{align}

\chapter{Generalised primes in short intervals}
In this chapter the $2$-point correlation derived in the previous chapter is related to fluctuations of the sum of the generalised von Mangoldt function over short intervals. Throughout this chapter the non-trivial zeros of all functions in the Selberg class are assumed to have a real part of $1 / 2.$ To begin, the expression is put into a form similar to that used by Murty and Perelli in \cite{MR1692847}. In particular the test function, $f(z),$ is first given the exponential weight $\cos(z \log X) \exp(- z^2).$ This allows a form similar to Murty and Perelli's function to be deduced,
\begin{equation}
\tilde{\mathcal{F}}(X, T; F) = \sum_{m, n \geq 1} X^{i ((t(m; F) - t(n; F))} \exp(- (t(m; F) - t(n; F))^2) 1_{|t(m; F)|, |t(n; F)| \leq T}.
\end{equation}
It is then conjectured that $\tilde{\mathcal{F}}(X, T; F)$ is asymptotically equivalent to $\mathcal{F}(X, T; F)$ so that the analysis first performed in \cite{MR1018376} and then refined in \cite{MR2009438} may be followed. A modification to the arguments in \cite{MR1018376} and \cite{MR2009438} is then considered in Appendix C, making use of $\tilde{\mathcal{F}}(X, T; F)$ as a pose to $\mathcal{F}(X, T; F).$ In agreement with Murty and Perelli's conjecture made in \cite{MR1692847}, it is conjectured that
\begin{equation}
\tilde{\mathcal{F}}(X, T; F) \sim \mathcal{F}(X, T; F) \sim \frac{T \log X}{\pi}
\end{equation}
uniformly for $X \leq T^{\textrm{deg}(F)}$ and
\begin{equation}
\tilde{\mathcal{F}}(X, T; F) \sim \mathcal{F}(X, T; F) \sim \frac{T}{\pi} \left(\textrm{deg}(F) \log \frac{T}{2 \pi} + \log \mathfrak{q}(F) - \textrm{deg}(F)\right)
\end{equation}
uniformly for $X \geq T^{\textrm{deg}(F)}.$ 

It is shown that if $1 < A_1 < A_2 < \textrm{deg}(F)$ and
\begin{equation}
\mathcal{F}(X, T; F) = \frac{T \log X}{\pi} + O \left( T^{1 - c} \right)
\end{equation}
uniformly for $T^{A_1} \ll X \ll T^{A_2}$ for some $c > 0,$ then for any fixed $1 / A_2 < B_1 \leq B_2 < 1 / A_1$
\begin{align}
& Var^{mul}(X, \delta; F) = \frac{1}{6} \delta X (3 \log X - 4 \log 2) + O(\delta^{1 + c / 2} X) \nonumber \\
& + O_{\epsilon}(\delta^{1 - \epsilon} X (\delta X^{1 / A_1})^{- 2 A_1 / (4 A_1 + 1)}) + O_{\epsilon}(\delta^{1 - \epsilon} X (\delta X^{1 / A_2})^{1 / 2})
\end{align}
uniformly for $X^{- B_2} \ll \delta \ll X^{- B_1}.$ While for the regime $\textrm{deg}(F) < A_1 < A_2$ if
\begin{equation}
\mathcal{F}(X, T; F) \sim \frac{T}{\pi} \left(\textrm{deg}(F) \log \frac{T}{2 \pi} + \log \mathfrak{q}(F) - \textrm{deg}(F)\right) + O\left(T^{1 - c}\right)
\end{equation}
uniformly for $T^{A_1} \ll X \ll T^{A_2}$ for some $c > 0,$ then for any fixed $1 / A_2 < B_1 \leq B_2 < 1 / A_1$
\begin{align}
& Var^{mul}(X, \delta; F) \nonumber \\
= & \frac{1}{2} \delta (\textrm{deg}(F) \log \frac{1}{\delta} + \log \mathfrak{q}(F) + (1 - \gamma - \log 2 \pi) \textrm{deg}(F)) \nonumber \\
& + O(\delta^{1 + c / 2} X) + O_{\epsilon}(\delta^{1 - \epsilon} X (\delta X^{1 / A_2})^{1 / 2}) + O_{\epsilon}(\delta^{1 - \epsilon} X (\delta X^{1 / A_2})^{- 2 A_1 / (4 A_1 + 1)})
\end{align}
uniformly for $X^{- B_2} \ll \delta \ll X^{- B_1}.$

From these relations the variances for fixed differences are then found. It is shown that if $1 / \textrm{deg}(F) < B_1 < B_2 \leq B_3 < 1$ and
\begin{equation}
Var^{mul}(X, \delta; F) = \frac{1}{6} \delta X (3 \log X - 4 \log 2) + O(\delta^{1 + c} X)
\end{equation}
uniformly for $X^{- B_3} \ll \delta \ll X^{- B_1}$ for some $c > 0,$ then
\begin{align}
& Var^{fix}(X, h; F) = \frac{1}{6} h (6 \log X - (3 + 8 \log 2)) \nonumber \\
& + O_{\epsilon}(h X^{\epsilon} (h / X)^{c / 3}) + O_{\epsilon}(h X^{\epsilon} (h X^{- (1 - B_1)})^{1 / (3 (1 - B_1))}),
\end{align}
uniformly for $X^{1 - B_3} \ll h \ll X^{1 - B_2}.$ For the regime $0 < B_1 < B_2 \leq B_3 < 1 / \textrm{deg}(F)$ if
\begin{align}
Var^{mul}(X, \delta; F) = & \frac{1}{2} \delta X \left(\textrm{deg}(F) \log \frac{1}{\delta} - \textrm{deg}(F) (\log 2 \pi + \gamma - 1) + \log \mathfrak{q}(F)\right) \nonumber \\
& + O \left(\delta^{1+c} X \right)
\end{align}
uniformly for $X^{- B_3} \ll \delta \ll X^{- B_1}$ for some $c > 0,$ then
\begin{align}
Var^{fix}(X, h ; F) = & h \left(\textrm{deg}(F) \log \frac{X}{h} - \textrm{deg}(F) (\log 2 \pi + \gamma) + \log \mathfrak{q}(F)\right) \nonumber \\ 
& + O_\epsilon\left(h X^{\epsilon} (h / X)^{c /3} \right) + O_\epsilon \left( h X^{1 - (B_2 - B_1)/3 (1 - B_1) + \epsilon}\right)
\end{align}
uniformly for $X^{1 - B_3} \ll h \ll X^{1 - B_2}.$ Further relations between $\mathcal{F}(X, T; F),$ $Var^{mul}(X, \delta; F)$ and $Var^{fix}(X, h; F)$ may be found in \cite{Buionthevarof2016}.

\section{A reduced form of the $2$-point correlation statistic for the Selberg class}
The $2$-point correlation deduced in the previous chapter is now reduced by specifying the test function. Firstly
\begin{equation}
\tilde{\mathcal{F}}(X, T; F) = \sum_{m, n \geq 1} f(t(m; F) - t(n; F)) 1_{|t(m; F)|, |t(n; F)| \leq T}= J_1 + J_2 + O_\epsilon(T^{1 / 2 + \epsilon}),
\end{equation}
where
\begin{equation}
J_1 = \int_{|t| \leq T} d^{mean}(t; F) \, dt
\end{equation}
and
\begin{equation}
J_2 = \int_{|t| \leq T, |\eta| \leq 2 T} \cos(\eta \log X) \exp(- \eta^2) R_2(X, T; F) \, d\eta dt,
\end{equation}
with
\begin{align}
\tilde{\mathcal{F}}(X, T; F) = & \sum_{m, n \geq 1} f(t(m; F) - t(n; F)) 1_{|t(m; F)|, |t(n; F)| \leq T} \nonumber \\
= & \sum_{m, n \geq 1} X^{i ((t(m; F) - t(n; F))} \exp(- (t(m; F) - t(n; F))^2) 1_{|t(m; F)|, |t(n; F)| \leq T}.
\end{align}

Since $f = \cos(z \log X) \exp(- z^2)$ is even
\begin{equation}
\int_{|t| \leq 2 T} \eta^{2 k - 1} f(\eta) \, d\eta  = 0,
\end{equation}
and
\begin{align}
& \int_{|\eta| \leq 2 T} \eta^{2 k} f(\eta) \, d\eta = \sum_{j \geq 0} \frac{(- 1)^j (\log X)^{2 j}}{(2 j)!} \int_{|\eta| \leq 2 T} \eta^{2 (k + j)} \exp(- \eta^2) \, d\eta \nonumber \\
= & \sqrt{\pi}\sum_{j \geq 0} \frac{(- 1)^j (2 k + 2 j)!}{2^{2 k + 2 j} (2 j)! (k + j)!} (\log X)^{2 j} + O((2T)^{2 k - 1} \exp(2 (\log X) T - 4 T^2))
\end{align}
for any $k \in \mathbb{Z}$. In particular,
\begin{align}
& \int_{|\eta| \leq 2 T} \eta^{2 k} f(\eta) \, d\eta \ll (\log X / 2)^{2 k} \exp(- (\log X)^2 / 4) \nonumber \\
& + (2 T)^{2 k - 1} \exp(2 (\log X) T - 4 T^2)
\end{align}
for any $k \geq 0$.

Moreover
\begin{equation}
\int_{|\eta| \leq 2 T} \eta^{2 k - 1} \cos(2 \pi d^{mean}(t; F) \eta) f(\eta) \, d\eta = 0,
\end{equation}
and
\begin{align}
& \int_{|\eta| \leq 2 T} \eta^{2 k} \cos(2 \pi d^{mean}(t; F) \eta) f(\eta) \, d\eta \nonumber \\
= & \sum_{i, j \geq 0} \frac{(- 1)^{i + j} (2 \pi d^{mean}(t; F))^{2 i} (\log X)^{2 j}}{(2 i)! (2 j)!} \int_{|\eta| \leq 2 T} \eta^{2 (k + i + j)} \exp(- \eta^2) \, d\eta \nonumber \\
= & \sqrt{\pi} \sum_{i, j \geq 0} \frac{(- 1)^{i + j} (2 k + 2 i + 2 j)!}{2^{2 k + 2 i + 2 j} (2 i)! (2 j)! (k + i + j)!} (2 \pi (d^{mean}(t; F))^{2 i} (\log X)^{2 j} \nonumber \\
& + O((2 T)^{2 k - 1} \exp(2 (\log X + 2 \pi d^{mean}(T; F)) T - 4 T^2))
\end{align}
for any $k \in \mathbb{Z}$. 

In particular,
\begin{align}
\int_{|\eta| \leq 2 T} \eta^{2 k} \cos(2 \pi d^{mean}(t; F) \eta) f(\eta) \, d\eta \ll ((\log X + 2 \pi d^{mean}(t; F)) / 2)^{2 k} & \nonumber \\
\exp(- (\log X - 2 \pi d^{mean}(t; F))^2 / 4) & \nonumber \\
+ (2 T)^{2 k - 1} \exp(2 (\log X + 2 \pi d^{mean}(T; F)) T - 4 T^2) &
\end{align}
for any $k \geq 0,$ and hence for $\alpha < \textrm{deg}(F)$
\begin{align}
& \int_{|t| \leq T, |\eta| \leq 2 T} \eta^{2 k} \cos(2 \pi d^{mean}(t; F) \eta) f(\eta) \, d\eta dt \ll_{\epsilon} T^{\alpha / \textrm{deg}(F) + \epsilon} (d^{mean}(T; F))^{2 k} \nonumber \\
& + (2 T)^{2 k} \exp(2 (\log X + 2 \pi d^{mean}(T; F)) T - 4 T^2),
\end{align}
whereas for $\alpha > \textrm{deg}(F)$
\begin{align}
& \int_{|t| \leq T, |\eta| \leq 2 T} \eta^{2 k} \cos(2 \pi d^{mean}(t; F) \eta) f(\eta) \, d\eta dt \ll_{\epsilon} T (\log X)^{2 k} \exp(- c (\log X)^2) \nonumber \\
& + (2 T)^{2 k} \exp(2 (\log X + 2 \pi d^{mean}(T; F)) T - 4 T^2)
\end{align}
with some absolute constant $c  > 0$, for any $k \geq 0$. Similarly,
\begin{equation}
\int_{|\eta| \leq 2 T} \eta^{2 k} \sin(2 \pi d^{mean}(t; F) \eta) f(\eta) \, d\eta = 0,
\end{equation}
and for $\alpha < \textrm{deg}(F)$
\begin{align}
& \int_{|t| \leq T, |\eta| \leq 2 T} \eta^{2 k + 1} \sin(2 \pi d^{mean}(t; F) \eta) f(\eta) \, d\eta dt \ll_{\epsilon} T^{\alpha / \textrm{deg}(F) + \epsilon} (d^{mean}(T; F))^{2 k + 1} \nonumber \\
& + (2 T)^{2 k} \exp(2 (\log X + 2 \pi d^{mean}(T; F)) T - 4 T^2),
\end{align}
whereas for $\alpha > \textrm{deg}(F)$
\begin{align}
& \int_{|t| \leq T, |\eta| \leq 2 T} \eta^{2 k + 1} \sin(2 \pi d^{mean}(t; F) \eta) f(\eta) \, d\eta dt \ll_{\epsilon} T (\log X)^{2 k + 1} \exp(- c (\log X)^2) \nonumber \\
& + (2 T)^{2 k + 1} \exp(2 (\log X + 2 \pi d^{mean}(T; F)) T - 4 T^2)
\end{align}
with some absolute constant $c  > 0$, for any $k \geq 0$.

Expanding terms then gives
\begin{align}
J_2 = & \int_{|t| \leq T, |\eta| \leq 2 T} (d^{mean}(t; F))^2 f(\eta) \, d\eta dt - \int_{|t| \leq T, |\eta| \leq 2 T} 2 \eta^{- 2} f(\eta) \, d\eta dt \nonumber \\
& + \int_{|t| \leq T, |\eta| \leq 2 T} 2 \eta^{- 2} \cos(2 \pi d^{mean}(t; F) \eta) f(\eta) \, d\eta dt + E \nonumber \\
= & 2 \pi^{- 3 / 2} \int_{|t| \leq T} \sum_{i \geq 2, j \geq 0} \frac{(- 1)^{i + j} \pi^{2 i} (2 i + 2 j - 2)!}{2^{2 j} (2 i)! (2 j)! (i + j - 1)!} (d^{mean}(t; F))^{2 i} (\log X)^{2 j} \, dt \nonumber \\
& + E,
\end{align}
where for every $A > 0$
\begin{equation}
E \ll_{\epsilon, A} \begin{dcases*} T^{\alpha / \textrm{deg}(F) + \epsilon} & if $\alpha < \textrm{deg}(F)$ \\ T^{- A} & otherwise. \end{dcases*}
\end{equation}
Inserting the sum into the symbolic manipulation software, Mathematica, gives
\begin{align}
& - \frac{\sqrt{\pi}}{8} |\log X - 2 \pi d^{mean}(t; F)| \textrm{erf}\left(\frac{|\log X - 2 \pi d^{mean}(t; F)|}{2}\right) \nonumber \\
& - \frac{\sqrt{\pi}}{8} (\log X + 2 \pi d^{mean}(t; F)) \textrm{erf}\left(\frac{\log X + 2 \pi d^{mean}(t; F)}{2}\right) \nonumber \\
& + \frac{\sqrt{\pi}}{4} \log X \textrm{erf}\left(\frac{\log X}{2}\right) + O\left(\exp(- (\log X - 2 \pi d^{mean}(t; F))^2 / 4)\right) \nonumber \\
& + O\left((d^{mean}(T; F))^2 \exp(- (\log X)^2 / 4)\right) \nonumber \\
= & - \frac{\sqrt{\pi}}{4} (\max\{\log X, 2 \pi d^{mean}(t; F)\} - \log X) \nonumber \\
& + O\left((d^{mean}(T; F))^2 \exp(- (\log X - 2 \pi d^{mean}(t; F))^2 / 4)\right) \nonumber \\
& + O\left((d^{mean}(T; F))^2 \exp(- (\log X)^2 / 4)\right).
\end{align}

Giving
\begin{equation}
J_2 = - \frac{1}{2 \pi} \int_{|t| \leq T} (\max\{\log X, 2 \pi d^{mean}(t; F)\} - \log X) \, dt + E.
\end{equation}

Together with $J_1$ this then gives the expression
\begin{equation}
\tilde{\mathcal{F}}_F(X,T) = \frac{T \log X}{\pi} + O_\epsilon(X^{1 / \textrm{deg}(F) + \epsilon}) + O_\epsilon(T^{1 / 2 + \epsilon})
\end{equation}
unifomly for $X \leq T^{\textrm{deg}(F) - \epsilon}$ and
\begin{equation}
\tilde{\mathcal{F}}(X,T; F) = \frac{T}{\pi} \left(\textrm{deg}(F) \log \frac{T}{2 \pi} + \log \mathfrak{q}(F) - \textrm{deg}(F)\right) + O_\epsilon(T^{1 / 2 + \epsilon})
\end{equation}
uniformly for $X \geq T^{\textrm{deg}(F) + \epsilon}$.
\section{$Var^{mul}(X, \delta; F)$}
Assuming the same expression holds for $\mathcal{F}(X, T; F)$ consideration is now given to
\begin{align}
I(X,T) = & \int_{|t| \leq T} \left|\sum_{n \geq 1} \frac{X^{i t(n; F)}}{1 + (t - t(n; F))^2} 1_{|t(n; F)| \leq Z} \right|^2 \, dt \nonumber \\
= & \sum_{m, n \geq 1} X^{i (t(m; F) - t(n; F))} 1_{|t(m; F)|, |t(n; F)| \leq Z} \nonumber \\
& \int_{|t| \leq T} \frac{1}{(1 + (t - t(m; F))^2)(1 + (t - t(n; F))^2)} \, dt.
\end{align}
Given that the number of zeros between $t \textrm{ and } t + 1 \textrm{ is } \ll \log |t|,$ the summation can be restricted to zeros with $|t(m; F)| \leq T \leq Z \textrm{ and } |t(n; F)| \leq T \leq Z$ incurring an error of size $\ll (\log T)^2.$ Similarly the domain of integration can be extended to the whole real line, incurring an error of size $\ll (\log T)^3.$ Therefore
\begin{align}
I(X, T) = & \sum_{m, n \geq 1} X^{i (t(m; F) - t(n; F))} 1_{|t(m; F)|, |t(n; F)| \leq T} \nonumber \\
& \int \frac{1}{(1 + (t - t(m; F))^2)(1 + (t - t(n; F))^2)} \, dt + O \left( \left( \log T \right)^3 \right) \nonumber \\
= & \frac{\pi}{2} \mathcal{F}(X, T; F) + O \left( \left( \log T \right)^3 \right).
\end{align}
In the regime $1 < A_1 < A_2 < \textrm{deg}(F)$ it then follows that
\begin{equation}
I(X,T) = \frac{T \log X}{2} + O \left( T^{1-c} \right)
\end{equation}
uniformly for $X^{1 / A_2} \ll T \ll X^{1 / A_1},$ whereas for the regime $\textrm{deg}(F) < A_1 < A_2$
\begin{equation}
I(X,T) = \frac{T}{2} \left(\textrm{deg}(F) \log \frac{T}{2 \pi} + \log \mathfrak{q}(F) - \textrm{deg}(F)\right) T + O \left( T^{1-c} \right)
\end{equation}
uniformly for $X^{1 / A_2} \ll T \ll X^{1 / A_1}$.

Taking
\begin{equation}
a(s) = \frac{(1 + \delta)^s - 1}{s}
\end{equation}
gives
\begin{equation}
\left|a(it)\right|^2 = 4 \left(\frac{\sin \kappa t}{t}\right)^2,
\end{equation}
where $\kappa = \frac{\log(1 + \delta)}{2}$ and
\begin{equation}
a(s) \ll \min\left\{\delta, 1 / |s| \right\} \textrm{ and } a'(s) \ll \min\left\{\delta^2, \delta / |s|\right\},
\end{equation}
for $\Re(s) \leq 1.$ To proceed a lemma is now demonstrated.
\subsection{Lemma A}
Taking $f$ to be a non-negative function satisfying the bound $f(t) = O_{\epsilon}(|t|^{\epsilon}),$ if 
\begin{equation}
R(u) = \int_{|t| \leq u} f(t) \, dt - u = O(u^{1 - c})
\end{equation}
uniformly for $U_1 \leq u \leq U_2,$ for some $0 < c < 1,$ where $U_1 = \kappa^{- (1 - c_1)} \textrm{ and } U_2 = \kappa^{- (1 + c_2)},$ for some $0 < c_1 < 1 \textrm{ and } 0 < c_2 < 1,$ then
\begin{equation}
\int \left(\frac{\sin \kappa u}{u}\right)^2 f(u) \, du = \frac{\pi}{2} \kappa + O(\kappa^{1 + c}) + O_{\epsilon}(\kappa^{1 + c_1 - \epsilon}) + O_{\epsilon}(\kappa^{1 + c_2 - \epsilon})
\end{equation}
as $\kappa \to 0^{+}.$
\subsubsection{Proof}
Firstly
\begin{align}
& \int \left(\frac{\sin \kappa u}{u}\right)^2 f(u) \, du = \int_{U_1 \leq u \leq U_2} \left(\frac{\sin \kappa u}{u}\right)^2 (f(u) + f(- u)) \, du \nonumber \\
+ & \int_{u \geq U_2} \left(\frac{\sin \kappa u}{u}\right)^2 (f(u) + f(- u)) \, du +  \int_{|u| \leq U_1} \left(\frac{\sin \kappa u}{u}\right)^2 f(u) \, du.
\end{align}
The first term is expressible as
\begin{equation}
\int_{U_1 \leq u \leq U_2} \left(\frac{\sin \kappa u}{u}\right)^2 \, dR(u) + \int_{U_1 \leq u \leq U_2} \left(\frac{\sin \kappa u}{u}\right)^2 \, du.
\end{equation}
Now
\begin{align}
& \int_{U_1 \leq u \leq U_2} \left(\frac{\sin \kappa u}{u}\right)^2 \, du = \int_{u \geq 0} \left(\frac{\sin \kappa u}{u}\right)^2 \, du \nonumber \\
- & \int_{0 \leq u \leq U_1} \left(\frac{\sin \kappa u}{u}\right)^2 \, du -  \int_{u \geq U_2} \left(\frac{\sin \kappa u}{u}\right)^2 \, du,
\end{align}
where
\begin{equation}
\int_{u \geq 0} \left(\frac{\sin \kappa u}{u}\right)^2 \, du = \kappa \int_{v \geq 0} \left(\frac{\sin v}{v}\right)^2 \, dv.
\end{equation}
Using the definite integral
\begin{equation}
\int_{v \geq 0} \left(\frac{\sin v}{v}\right)^2 \, dv = \frac{\pi}{2}
\end{equation}
then gives
\begin{equation}
\int_{u \geq 0} \left(\frac{\sin \kappa u}{u}\right)^2 \, du = \frac{\pi}{2} \kappa.
\end{equation}
The remaining terms contribute to the errors. Integrating by parts gives
\begin{align}
\int_{U_1 \leq u \leq U_2} \left(\frac{\sin \kappa u}{u}\right)^2 \, dR(u) & \ll \int_{U_1 \leq u \leq U_2} \left(\left|\frac{\kappa \sin 2 \kappa u}{u^2}\right| + \left|\frac{(\sin \kappa u)^2}{u^3}\right|\right) |R(u)| \, du \nonumber \\
+ U_2^{- 2} R(U_2) + \kappa^2 R(U_1),
\end{align}
where
\begin{equation}
\int_{U_1 \leq u \leq U_2} \left|\frac{\kappa \sin 2 \kappa u}{u^2}\right| |R(u)| \, du \ll \kappa^2 \int_{U_1 \leq u \leq U_2} u^{- c} \, du \ll \kappa^2 U_1^{1 - c}, 
\end{equation}
with
\begin{equation}
\kappa^2 U_1^{1 - c} = \kappa^2 \kappa^{- (1 - c_1) (1 - c)} = \kappa^{1 + c + c_1 (1 - c)} \ll \kappa^{1 + c},
\end{equation}
similarly
\begin{equation}
\int_{U_1 \leq u \leq U_2} \left|\frac{(\sin \kappa u)^2}{u^3}\right| |R(u)| \, du \ll \kappa^2 \int_{U_1 \leq u \leq U_2} u^{- c} \, du \ll \kappa^{1 + c}.
\end{equation}
Also
\begin{equation}
U_2^{- 2} R(U_2) \ll U_2^{- 1 - c},
\end{equation}
where
\begin{equation}
U_2^{- 1 - c} = \kappa^{- (1 + c_2) (- 1 - c)} = \kappa^{1 + c + c_2 (1 + c)} \ll \kappa^{1 + c},
\end{equation}
and
\begin{equation}
\kappa^2 R(U_1) \ll \kappa^2 U_1^{1 - c},
\end{equation}
where
\begin{equation}
\kappa^2 U_1^{1 - c} = \kappa^{2 - (1 - c_1) (1 - c)} = \kappa^{1 + c + c_1 (1 - c)} \ll \kappa^{1 + c},
\end{equation}
giving 
\begin{equation}
\int_{U_1 \leq u \leq U_2} \left(\frac{\sin \kappa u}{u}\right)^2 \, dR(u) \ll \kappa^{1 + c}.
\end{equation}
Also
\begin{equation}
\int_{0 \leq u \leq U_1} \left(\frac{\sin \kappa u}{u}\right)^2 \, du \ll \kappa^2 \int_{0 \leq u \leq U_1} \, du \ll \kappa^{1 + c_1},
\end{equation}
and
\begin{equation}
\int_{u \geq U_2} \left(\frac{\sin \kappa u}{u}\right)^2 \, du \ll \int_{u \geq U_2} u^{- 2} \, du \ll \kappa^{1 + c_2}.
\end{equation}
Similarly
\begin{equation}
\int_{u \geq U_2} \left(\frac{\sin \kappa u}{u}\right)^2 (f(u) + f(- u)) \, du \ll_{\epsilon} \int_{u \geq U_2} u^{- 2 + \epsilon} \, du \ll_{\epsilon} U_2^{- 1 + \epsilon} \ll_{\epsilon} \kappa^{1 + c_2 - \epsilon}
\end{equation}
and
\begin{equation}
\int_{|u| \leq U_1} \left(\frac{\sin \kappa u}{u}\right)^2 f(u) \, du \ll_{\epsilon} \kappa^2 \int_{0 < u \leq U_1} u^{\epsilon} \, du \ll_{\epsilon} \kappa^2 U_1^{1 + \epsilon} \ll_{\epsilon} \kappa^{1 + c_1 - \epsilon}. 
\end{equation}
$\Box$

Applying Lemma A in the range
\begin{equation}
\delta^{- (1 - c_1)} \ll T \ll \delta^{- (1 + c_2)},
\end{equation}
gives
\begin{align}
& \int \left|a(it)\right|^2 \left|\sum_{n \geq 1} \frac{X^{i t(n; F)}}{1 + (t - t(n; F))^2} 1_{|t(n; F)| \leq Z}\right|^2 \, dt \nonumber \\
= & \pi \kappa \log X + O\left(\kappa^{1 + c}\right) + O_{\epsilon}\left(\kappa^{1 + c_1 - \epsilon}\right) + O_{\epsilon}\left(\kappa^{1 + c_2 - \epsilon}\right) \nonumber \\
= & \frac{\pi}{2} \delta \log X + O\left(\delta^{1 + c}\right) + O_{\epsilon}\left(\delta^{1 + c_1 - \epsilon}\right) + O_{\epsilon}\left(\delta^{1 + c_2 - \epsilon}\right).
\end{align}
Denoting the integral above by $J$ and the integral with $a(i t)$ replaced by $a(1 / 2 + i t(n; F))$ by $K,$ with $J = \int|A|^2$ and $K =\int|B|^2,$ since the number of zeros between $t \textrm{ and } t + 1 \textrm{ is } \ll \log(|t|),$
\begin{equation}
A, B \ll \min \left\{ \delta, 1/|t| \right\} \log|t|
\end{equation}
and
\begin{equation}
a(i t) - a(1 / 2 + i t(n; F)) \ll \left(1 + |t - t(n; F)|\right) \min\left\{\delta^2, \delta / |t| \right\}.
\end{equation}
Therefore
\begin{equation}
A - B \ll \min \left\{ \delta^2, \delta/|t| \right\} \left( \log|t| \right)^2,
\end{equation}
giving
\begin{equation}
|A|^2 - |B|^2 \ll \min \left\{ \delta^3, \delta/|t|^2 \right\} \left( \log|t| \right)^3
\end{equation}
and
\begin{equation}
J - K \ll \delta^2 \left( \log \frac{1}{\delta} \right)^3.
\end{equation}
Using this then gives
\begin{align}
& \int \left|\sum_{n \geq 1} \frac{a(1 / 2 + i t(n; F)) X^{i t(n; F)}}{1 + (t - t(n; F))^2} 1_{|t(n; F)| \leq Z}\right|^2 \, dt \nonumber \\ 
= & \frac{\pi}{2} \delta \log X + O\left(\delta^{1 + c}\right) + O_\epsilon\left(\delta^{1 + c_1 - \epsilon}\right) + O_\epsilon\left(\delta^{1 + c_2 - \epsilon}\right).
\end{align}

Making use of the Fourier transform,
\begin{equation}
\int \frac{1}{1 + (t - u)^2} \exp(- 2 \pi i t v) \, dt = \pi \exp(- 2 \pi i u v) \exp(- 2 \pi |v|),
\end{equation}
an application of Plancherel's formula gives
\begin{align}
& \int \left|\sum_{n \geq 1} \frac{a(1 / 2 + i t(n; F)) X^{i t(n; F)}}{1 + (t - t(n; F))^2} 1_{|t(n; F)| \leq Z}\right|^2 \, dt \nonumber \\
= & \frac{\pi}{2} \int \left|\sum_{n \geq 1} a(1 / 2 + i t(n; F)) \exp(i t(n; F) (Y + y)) 1_{|t(n; F)| \leq Z}\right|^2 \exp(- 2 |y|) \, dy,
\end{align}
using the substitution $Y = \log X$ and $y = - 2 \pi v.$ A special case of lemma 1 of \cite{MR2927496} is now used to remove the weight.
\subsection{Lemma B}
Taking $f$ to be a non-negative function satisfying
\begin{equation}
\int f(T + y) \exp(- 2 |y|) \, dy = 1 + O(\exp(- c Y))
\end{equation}
for $Y \leq T \leq Y + \log 2$ for some $c > 0,$ then
\begin{equation}
\int_{0 \leq y \leq \log 2} f(Y + y) \exp(2 y) \, dy = \frac{3}{2} + O(\exp(- c Y / 2)).
\end{equation}

Applying Lemma B now gives
\begin{align}
& \int_{X \leq x \leq 2 X} \left|\sum_n a(1 / 2 + i t(n; F)) x^{1 / 2 + i t(n; F)} 1_{|t(n; F)| \leq Z}\right|^2 \, dx \nonumber \\
= & \frac{3}{2} \delta X^2 \log X + O(\delta^{1 + c / 2} X^2) + O_{\epsilon}(\delta^{1 + c_1 / 2 - \epsilon} X^2) + O_{\epsilon}(\delta^{1 + c_2 / 2 - \epsilon} X^2)
\end{align}
after the change of variable $x = \exp(Y + y),$ provided that
\begin{equation}
X^{1 / A_2} \ll \delta^{- (1 - c_1)} < \delta^{- (1 + c_2)} \ll X^{1 / A_1}.
\end{equation}
Taking $Z = X^2$ in the explicit formula then gives
\begin{align}
& \int_{X \leq x \leq 2 X} |\psi(x, \delta x; F) - m(F) \delta x|^2 \, dx \nonumber \\
= & \frac{3}{2} \delta X^2 \log X + O(\delta^{1 + c / 2} X^2) + O_{\epsilon}(\delta^{1 + c_1 / 2 - \epsilon} X^2) + O_{\epsilon}(\delta^{1 + c_2 / 2 - \epsilon} X^2)
\end{align}
Summing over the dyadic intervals $[2^{- j} X, 2^{- j + 1} X],$ with $1 \leq j \leq J$ and
\begin{equation}
2^J = \delta^{(1 + c_2) A_1} X
\end{equation}
and noting that
\begin{align}
&3 \sum_{1 \leq j \leq J} j 4^{- j} = 4 \sum_{1 \leq j \leq J} j (4^{- j} - 4^{- (j + 1)}) = 1 + 4 \sum_{2 \leq j \leq J} 4^{- j} - J 4^{- J} \nonumber \\
=&4 \sum_{1 \leq j \leq J} 4^{- j} - J 4^{- J} = \frac{4 - 4^{- J + 1}}{3} - J 4^{- J} = \frac{4 - (3 J + 4) 4^{- J}}{3}
\end{align}
then gives
\begin{align}
& \int_{2^{- J} X \leq x \leq X} |\psi(x, \delta x; F) - m(F) \delta x|^2 \, dx \nonumber \\
= & \frac{1 - 4^{- J}}{2} \delta X^2 \log X - \frac{4 - (3 J + 4) 4^{- J}}{6} \delta X^2 \log 2 \nonumber \\
& + O(\delta^{1 + c / 2} X^2) + O_{\epsilon}(\delta^{1 + c_1 / 2 - \epsilon} X^2) + O_{\epsilon}(\delta^{1 + c_2 / 2 - \epsilon} X^2).
\end{align}
To bound the contribution coming from the range $2^{- J} X \leq x \leq X$ the steps from lemma $5$ of the work of Saffari and Vaughan in \cite{MR0480388} are now followed, adapted to the Selberg class.
\subsection{Lemma C}
\begin{equation}
Var^{mul}(X, \delta; F) \ll \delta X^2 \left(\log \frac{2}{\delta}\right)^2
\end{equation}
for $0 < \delta \leq 1.$
\subsubsection{Proof}
Taking $Z \to \infty$ in the explicit formula gives
\begin{align}
|\psi(x, \delta x; F) - m(F) \delta x| \leq & x^{\frac{1}{2}} \left|\sum_{n \geq 1} a(1 / 2 + i t(n; F)) x^{i t(n; F)}\right| \nonumber \\ 
& + \left|\sum_{1 \leq j \leq r} \lambda_{j} \sum_{n \geq 1} \frac{(x + \delta x)^{- n / (\max_{1 \leq k \leq r} \lambda_k)}}{n} \right|
\end{align}
so that
\begin{align}
& |\psi(x, \delta x; F) - m(F) \delta x|^2 \nonumber \\
= & x \sum_{m, n \geq 1} a(1 / 2 + i t(m; F)) a^*(1 / 2 + i t(n; F)) x^{i (t(m; F) - t(n; F))} \nonumber \\ 
& + O\left(\left(\log \left(1 - (x + \delta x)^{- \frac{1}{\max_{1 \leq k \leq r} \lambda_{k}}}\right)\right)^2\right).
\end{align}
Consideration is now given to
\begin{align}
\int_{X \leq x \leq 2 X} |\psi(x, \delta x; F) - m(F) \delta x|^{2} \ dx & \nonumber \\
\ll \int_{1 \leq \tau \leq 2} \int_{\frac{1}{2} X \tau \leq x \leq 2 X \tau} |\psi(x, \delta x; F) - m(F) \delta x|^2 \, dx d\tau. &
\end{align}
Firstly,
\begin{equation}
a(1 / 2 + i t(m; F)) a^*(1 / 2 + i t(n; F)) = \frac{(1 + \delta)^{1 / 2 + i t(m; F)} - 1}{1 / 2 + i t(m; F)} \frac{(1 + \delta)^{1 / 2 + i t(n; F)} - 1}{1 / 2 + i t(n; F)},
\end{equation}
also
\begin{align}
& \int_{1 \leq \tau \leq 2} \int_{\frac{1}{2} X \tau \leq x \leq 2 X \tau} x x^{i (t(m; F) - t(n; F))} \, dx d\tau \nonumber \\
= & \frac{2^{2 + i (t(m; F) - t(n; F))} - 2^{- 2 - i (t(m; F) - t(n; F))}}{2 + i (t(m; F) - t(n; F))} X^{2 + i (t(m; F) - t(n; F))} \nonumber \\
& \int_{1 \leq \tau \leq 2} \tau^{2 + i (t(m; F) - t(n; F))} \, d\tau \nonumber \\ 
= & \frac{2^{2 + i (t(m; F) - t(n; F))} - 2^{- 2 - i (t(m; F) - t(n; F))}}{2 + i (t(m; F) - t(n; F))} X^{2 + i (t(m; F) - t(n; F))} \nonumber \\
& \frac{2^{3 + i (t(m; F) - t(n; F))} - 1}{3 + i (t(m; F) - t(n; F))},
\end{align}
with
\begin{align}
& \frac{\left|\left(1 / 2 + i t(n; F)\right) \delta + O(\delta^2)\right|^2}{1 / 4 + (t(n; F))^2} \ll \frac{\left|\left(1 / 2 + i t(n; F)\right) \delta + O(\delta^2)\right|^2}{1 / 4 + (t(n; F))^2} \nonumber \\
& \frac{2^4 + 2^{- 4}}{1 + |t(m; F) - t(n; F)|^2} \nonumber \\
\ll & \frac{\delta^2 + \delta^2 (t(n; F))^{- 2}}{(1 + |t(m; F) - t(n; F)|)^2}.
\end{align}
Similarly
\begin{equation}
\left|\frac{(1 + \delta)^{1 / 2 + i t(m; F)} - 1}{1 / 2 + i t(m; F)} \frac{2^{3 + i (t(m; F) - t(n; F))} - 1}{3 + i (t(m; F) - t(n; F))}\right|^2 \ll \frac{\delta^2 + \delta^2 (t(m; F))^{- 2}}{(1 + |t(m; F) - t(n; F)|)^2}.
\end{equation}
Using the bound of $\log t(n; F)$ for the spacing between the $n$th zero ordinate and the next closest zero ordinate gives
\begin{align}
& \int_{X \leq x \leq 2 X} |\psi(x, \delta x; F) - m(F) \delta x|^2 \, dx \nonumber \\
\ll & X^2 \sum_{m, n \geq 1} \min(\delta^2, (t(m; F))^{- 2}) (1 + |t(m; F) - t(n; F)|)^{- 2} \nonumber \\ 
\ll & \delta^2 X^2 \sum_{m \geq 1} \log t(m; F) 1_{0 < t(m; F) \leq \delta^{- 1}} \ll \delta X^2 \left( \log \frac{2}{\delta}\right)^2
\end{align}
and the lemma then follows by summing over dyadic intervals.
$\Box$

Applying Lemma C then gives
\begin{align}
& Var^{mul}(X, \delta; F) = \frac{1}{6} \delta X (3 \log X - 4 \log 2) + O(\delta^{1 + c / 2} X^2) \nonumber \\
& + O_{\epsilon}(\delta^{1 + c_1 / 2 - \epsilon} X) + O_{\epsilon}(\delta^{1 + c_2 / 2 - \epsilon} X) + O_{\epsilon}(\delta^{1 - \epsilon} X 4^{- J}),
\end{align}
therefore
\begin{align}
& Var^{mul}(X, \delta; F) = \frac{1}{6} \delta X (3 \log X - 4 \log 2) + O(\delta^{1 + c / 2} X) \nonumber \\
& + O_{\epsilon}(\delta^{1 - \epsilon} X (\delta X^{1 / A_1})^{- 2 A_1 / (4 A_1 + 1)}) + O_{\epsilon}(\delta^{1 - \epsilon} X (\delta X^{1 / A_2})^{1 / 2}),
\end{align}
uniformly for $X^{- B_2} \ll \delta \ll X^{- B_1}$ for some $c > 0$ and any fixed $1 / A_2 < B_1 \leq B_2 < 1 / A_1.$

For the regime $\textrm{deg}(F) < A_1 < A_2$ the analysis is almost identical, as may be seen in \cite{Buionthevarof2016}. Changing the function in Lemma A to $R(u) = \int_{|t| \leq u} f(t) \, dt - u (\log u + A)$ results in
\begin{equation}
\int \left(\frac{\sin \kappa u}{u}\right)^2 f(u) \, du = \frac{\pi}{2} \kappa \left(\log \frac{1}{\kappa} + B\right) + O(\kappa^{1 + c}) + O_{\epsilon}(\kappa^{1 + c_1 - \epsilon}) + O_{\epsilon}(\kappa^{1 + c_2 - \epsilon})
\end{equation}
as $\kappa \to 0^{+},$ where $B = A + 2 - \gamma - \log 2.$ Applying this in place of Lemma A then gives
\begin{align}
& \int \left|a(it)\right|^2 \left|\sum_{n \geq 1} \frac{X^{i t(n; F)}}{1 + (t - t(n; F))^2} 1_{|t(n; F)| \leq Z}\right|^2 \, dt \nonumber \\
= & \pi \kappa \left(\textrm{deg}(F) \log \frac{1}{\kappa} + \log \mathfrak{q}(F) + (1 - \gamma - \log 4 \pi) \textrm{deg}(F)\right) + O \left( \kappa^{1 + c} \right) \nonumber \\
& + O_{\epsilon} \left( \kappa^{1 + c_1 - \epsilon} \right) + O_{\epsilon} \left( \kappa^{1 + c_2 - \epsilon}\right) \nonumber \\
= & \frac{\pi}{2} \delta \left(\textrm{deg}(F) \log \frac{1}{\delta} + \log \mathfrak{q}(F) + (1 - \gamma - \log 2 \pi) \textrm{deg}(F) \right) + O \left( \delta^{1 + c} \right) \nonumber \\
& + O_{\epsilon} \left( \delta^{1 + c_1 - \epsilon} \right) + O_{\epsilon} \left( \delta^{1 + c_2 - \epsilon} \right).
\end{align}
The calculation then proceeds as before giving
\begin{align}
& Var^{mul}(X, \delta; F) \nonumber \\
= & \frac{1}{2} \delta (\textrm{deg}(F) \log \frac{1}{\delta} + \log \mathfrak{q}(F) + (1 - \gamma - \log 2 \pi) \textrm{deg}(F)) \nonumber \\
& + O(\delta^{1 + c / 2} X) + O_{\epsilon}(\delta^{1 + c_1 / 2 - \epsilon} X) + O_{\epsilon}(\delta^{1 + c_2 / 2 - \epsilon} X) + O_{\epsilon}(\delta^{1 - \epsilon} X 4^{- J}).
\end{align}
\section{$Var^{fix}(X, h ; F)$}
It is now shown that if $0 < B_1 < B_2 \leq B_3 < 1 / \textrm{deg}(F)$ and
\begin{align}
Var^{mul}(X, \delta; F) = & \frac{1}{2} \delta X \left(\textrm{deg}(F) \log \frac{1}{\delta} - \textrm{deg}(F) (\log 2 \pi + \gamma - 1) + \log \mathfrak{q}(F)\right) \nonumber \\
& + O \left(\delta^{1+c} X \right)
\end{align}
uniformly for $X^{- B_3} \ll \delta \ll X^{- B_1}$ for some $c > 0,$ then
\begin{align}
Var^{fix}(X, h ; F) = & h \left(\textrm{deg}(F) \log \frac{X}{h} - \textrm{deg}(F) (\log 2 \pi + \gamma) + \log \mathfrak{q}(F)\right) \nonumber \\ 
& + O_\epsilon\left(h X^\epsilon (h/X)^{c/3} \right) + O_\epsilon \left( h X^{1 - (B_2 - B_1)/3 (1 - B_1) + \epsilon}\right)
\end{align}
uniformly for $X^{1 - B_3} \ll h \ll X^{1 - B_2}.$

Firstly setting $h = \delta x$ gives
\begin{align}
& \int_{X \leq x \leq 2 X, H_1 \leq h \leq H_2} |\psi(x, h; F) - m(F) h|^2 \, dh dx \nonumber \\
= & \int_{R} x |\psi(x, \delta x; F) - m(F) \delta x|^2 \, d\delta dx,
\end{align}
where
\begin{equation}
R = \{(\delta, x) : H_1/x \leq \delta \leq H_2/x \textrm{ and } X \leq x \leq 2 X \}.
\end{equation}
Assuming $H_1 < H_2 < 2 H_1,$ exchanging the order of integration gives an integral over $R = R_1 \cup R_2 \cup R_3,$ where
\begin{equation}
R_1 = \{(x, \delta) : H_1/\delta \leq x \leq 2 X \textrm{ and } H_1/2 X \leq \delta \leq H_2/2 X \},
\end{equation}
\begin{equation}
R_2 = \{(x, \delta) : H_1/\delta \leq x \leq H_2/\delta \textrm{ and } H_2/2 X \leq \delta \leq H_1/X \}
\end{equation}
and
\begin{equation}
R_3 = \{(x, \delta) : X \leq x \leq H_2/\delta \textrm{ and } H_1/X \leq \delta \leq H_2/X \}.
\end{equation}
The integrals over $x$ are of the form
\begin{align} \label{x integral (fixed delta)}
& \int_{X_1 \leq x \leq X_2} x |\psi(x, \delta x; F) - m(F) \delta x|^2 \, dx  \nonumber \\
= & X_2^2 Var^{mul}(X_2, \delta; F) - X_1^2 Var^{mul}(X_1, \delta; F) \nonumber \\
& - \int_{X_1 \leq x \leq X_2} x Var^{mul}(x, \delta; F) \, dx \nonumber \\ 
= & \frac{1}{3} \delta (X_2^3 - X_1^3) \left(\textrm{deg}(F) \log \frac{1}{\delta} + \log \mathfrak{q}(F) - \textrm{deg}(F) (\log 2 \pi + \gamma - 1) \right) \nonumber \\
 & + O(\delta^{1 + c} X_2^3).
\end{align}
For $X_1 = H_1 / \delta$ and $X_2 = 2 X$ the terms coming from \eqref{x integral (fixed delta)} are $$\frac{8}{3} X^3 \textrm{deg}(F) \delta \log \frac{1}{\delta}, - \frac{1}{3} H_1^3 \textrm{deg}(F) \delta^{- 2} \log \frac{1}{\delta},$$
$$\frac{8}{3} X^3 (\log \mathfrak{q}(F) - \textrm{deg}(F) (\log 2 \pi + \gamma - 1)) \delta,$$ $$- \frac{1}{3} H_1^3 (\log \mathfrak{q}(F) - \textrm{deg}(F) (\log 2 \pi + \gamma - 1)) \delta^{- 2} \textrm{ and } O(\delta^{1 + c} X^3),$$ for $X_1 = H_1 / \delta$ and $X_2 = H_2 / \delta$ they are $$\frac{1}{3} (H_2^3 - H_1^3) \textrm{deg}(F) \delta^{- 2} \log \frac{1}{\delta}, \frac{1}{3} (H_2^3 - H_1^3) (\log \mathfrak{q}(F) - \textrm{deg}(F) (\log 2 \pi + \gamma - 1)) \delta^{- 2}$$ $$\textrm{ and } O(\delta^{- 2 + c} H_1^3)$$ and for $X_1 = X$ and $X_2 = H_2 / \delta$ they are $$\frac{1}{3} H_2^3 \textrm{deg}(F) \delta^{- 2} \log \frac{1}{\delta}, - \frac{1}{3} X^3 \textrm{deg}(F) \delta \log \frac{1}{\delta},$$ $$\frac{1}{3} H_2^3 (\log \mathfrak{q}(F) - \textrm{deg}(F) (\log 2 \pi + \gamma - 1)) \delta^{- 2},$$ $$- \frac{1}{3} X^3 (\log \mathfrak{q}(F) - \textrm{deg}(F) (\log 2 \pi + \gamma - 1)) \delta \textrm{ and } O(\delta^{- 2 + c} H_1^3).$$
The integrals over $\delta$ are therefore of the form
\begin{equation}
\int_{\Delta_1 \leq \delta \leq \Delta_2} \delta \log \frac{1}{\delta} \; d\delta = - \frac{1}{4} (\Delta_2^2 (2 \log \Delta_2 - 1) - \Delta_1^2 (2 \log \Delta_1 - 1)),
\end{equation}
\begin{equation}
\int_{\Delta_1 \leq \delta \leq \Delta_2} \delta^{- 2} \log \frac{1}{\delta} \; d\delta = \Delta_2^{- 1} (\log \Delta_2 + 1) - \Delta_1^{- 1} (\log \Delta_1 + 1),
\end{equation}
\begin{equation}
\int_{\Delta_1 \leq \delta \leq \Delta_2} \delta \; d\delta = \frac{1}{2}(\Delta_2^2 - \Delta_1^2)
\end{equation}
and
\begin{equation}
\int_{\Delta_1 \leq \delta \leq \Delta_2} \delta^{- 2} \; d\delta = - (\Delta_2^{- 1} - \Delta_1^{- 1}).
\end{equation}
Setting $\Delta_1 = H_1 / 2 X$ and $\Delta_2 = H_2 / 2 X,$ the integral over $R_1$ yields
\begin{align}
& - \frac{2}{3} X^3 \textrm{deg}(F) \left(\frac{H_2}{2 X} \right)^2 \left( 2 \log \frac{H_2}{X} - 2 \log 2 - 1 \right) \nonumber \\
& + \frac{2}{3} X^3 \textrm{deg}(F) \left(\frac{H_1}{2 X} \right)^2 \left( 2 \log \frac{H_1}{X} - 2 \log 2 - 1 \right) \nonumber \\
& - \frac{1}{3} H_1^3 \textrm{deg}(F) \left( \frac{2 X}{H_2} \right) \left(\log \left( \frac{H_2}{X} \right) - \log 2 + 1 \right) \nonumber \\
& + \frac{1}{3} H_2^3 \textrm{deg}(F) \left( \frac{2 X}{H_1} \right) \left(\log \left( \frac{H_1}{X} \right) - \log 2 + 1 \right) \nonumber \\
&+ \frac{4}{3} X^3 (\log \mathfrak{q}(F) - \textrm{deg}(F) (\log 2 \pi + \gamma - 1)) \left( \left( \frac{H_2}{2 X} \right)^2 - \left( \frac{H_1}{2 X} \right)^2 \right) \nonumber \\ 
&+ \frac{1}{3} H_1^3 (\log \mathfrak{q}(F) - \textrm{deg}(F) (\log 2 \pi + \gamma - 1)) \left( \left( \frac{2 X}{H_2} \right) - \left( \frac{2 X}{H_1} \right) \right) \nonumber \\ 
&+ O(H_1^{2 + c} X^{1 - c}),
\end{align} 
setting $\Delta_1 = H_2 / 2 X$ and $\Delta_2 = H_1 / X,$ the integral over $R_2$ yields
\begin{align}
&\frac{1}{3} (H_2^3 - H_1^3) \textrm{deg}(F) \nonumber \\
& \left( \left( \frac{X}{H_1} \right) \left( \log \left( \frac{H_1}{X} \right) + 1 \right) - \left( \frac{2 X}{H_2} \right) \left( \log \left( \frac{H_2}{X} \right) - \log 2 + 1 \right) \right) \nonumber \\
&- \frac{1}{3} (H_2^3 - H_1^3) (\log \mathfrak{q}(F) - \textrm{deg}(F) (\log 2 \pi + \gamma - 1)) \left( \left( \frac{X}{H_1} \right) - \left( \frac{2 X}{H_2} \right) \right) \nonumber \\ 
&+ O(H_1^{2 + c} X^{1 - c})
\end{align}
and setting $\Delta_1 = H_1 / X$ and $\Delta_2 = H_2 / X,$ the integral over $R_3$ yields
\begin{align}
& \frac{1}{3} H_2^3 \textrm{deg}(F)\left( \left( \frac{X}{H_2} \right) \left( \log \left( \frac{H_2}{X} \right) + 1 \right) - \left( \frac{X}{H_1} \right) \left( \log \left( \frac{H_1}{X} \right) + 1 \right) \right) \nonumber \\
& + \frac{1}{12} X^3 \textrm{deg}(F) \left( \left( \frac{H_2}{X} \right)^2 \left( 2 \log \frac{H_2}{X} - 1 \right) - \left( \frac{H_1}{X} \right)^2 \left( 2 \log \frac{H_1}{X} - 1 \right) \right) \nonumber \\
&+ \frac{1}{3} H_2^3 (\log \mathfrak{q}(F) - \textrm{deg}(F) (\log 2 \pi + \gamma - 1)) \left( \left( \frac{X}{H_1} \right) - \left( \frac{X}{H_2} \right) \right) \nonumber \\ 
&- \frac{1}{6} X^3 (\log \mathfrak{q}(F) - \textrm{deg}(F) (\log 2 \pi + \gamma - 1)) \left( \left( \frac{H_2}{X} \right)^2 - \left( \frac{H_1}{X} \right)^2 \right) \nonumber \\ 
&+ O(H_1^{2 + c} X^{1 - c}).
\end{align}
The total contribution is therefore
\begin{align}
&\left( \frac{1}{3} H_1^2 + \frac{2}{3} H_1^3 + \frac{1}{3} (H_2^3 - H_1^3) H_1^{- 1} - \frac{1}{3} H_2^3 H_1^{- 1} - \frac{1}{6} H_1^2 \right) \textrm{deg}(F) X \log \frac{H_1}{X} \nonumber \\
&+\left( - \frac{1}{3} H_2^3 - \frac{2}{3} H_1^3 H_2^{- 1} - \frac{2}{3} (H_2^3 - H_1^3) H_2^{- 1} + \frac{1}{3} H_2^2 + \frac{1}{6} H_2^2 \right) \textrm{deg}(F) X \log \frac{H_2}{X} \nonumber
\end{align}
\begin{align}
&+\left( \frac{1}{6} (H_2^2 - H_1^2) (2 \log 2 + 1) \textrm{deg}(F) + \frac{2}{3} H_1^3 H_2^{- 1} \textrm{deg}(F) (\log 2 - 1) \right. \nonumber \\
&\left. \vphantom( - \frac{2}{3} H_1^2 \textrm{deg}(F) (\log 2 - 1) + \frac{1}{3} (H_2^2 - H_1^2) (\log \mathfrak{q}_F - \textrm{deg}(F) (\log 2 \pi + \gamma - 1)) \right. \nonumber \\
&\left. \vphantom( + \frac{2}{3} H_1^3 H_2^{- 1} (\log \mathfrak{q}_F - d_F (\log 2 \pi + \gamma - 1)) \right. \nonumber \\ 
& \left. \vphantom(- \frac{2}{3} H_1^2 (\log \mathfrak{q}_F - \textrm{deg}(F) (\log 2 \pi + \gamma - 1)) \right. \nonumber \\ 
&\left. \vphantom( + \frac{1}{3} (H_2^3 - H_1^3) H_1^{- 1} \textrm{deg}(F) + \frac{2}{3} (H_2^3 - H_1^3) \textrm{deg}(F) H_2^{- 1} (\log 2 - 1) \right. \nonumber \\
&\left. \vphantom( - \frac{1}{3} (H_2^3 - H_1^3) (\log \mathfrak{q}_F - \textrm{deg}(F) (\log 2 \pi + \gamma - 1)) H_1^{- 1} \right. \nonumber \\ 
&\left. \vphantom( + \frac{2}{3} (H_2^3 - H_1^3) (\log \mathfrak{q}_F - \textrm{deg}(F) (\log 2 \pi + \gamma - 1)) H_2^{- 1} + \frac{1}{3} H_2^2 \textrm{deg}(F) \right. \nonumber \\ 
&\left. \vphantom( - \frac{1}{3} H_2^3 H_1^{- 1} \textrm{deg}(F) - \frac{1}{12} (H_2^2 - H_1^2) \textrm{deg}(F) \right. \nonumber \\
&\left. \vphantom( + \frac{1}{3} H_2^3 H_1^{- 1} (\log \mathfrak{q}_F - \textrm{deg}(F) (\log 2 \pi + \gamma - 1)) \right. \nonumber \\
& \left. \vphantom(  - \frac{1}{3} H_2^2 (\log \mathfrak{q}_F - \textrm{deg}(F) (\log 2 \pi + \gamma - 1)) \right. \nonumber \\
&\left. \vphantom( - \frac{1}{6} (H_2^2 - H_1^2) (\log \mathfrak{q}_F - \textrm{deg}(F) (\log 2 \pi + \gamma - 1)) \right) X \nonumber \\
&+ O( H_1^{2 + c} X^{1 - c}).
\end{align} 
Collecting terms then gives
\begin{align}
&\int_{X \leq x \leq 2 X, H_1 \leq h \leq H_2} |\psi(x, h; F) - m(F) h|^2 \, dh dx \nonumber \\
=&\frac{1}{2} \textrm{deg}(F) \left( H_2^2 \log \frac{X}{H_2} - H_1^2 \log \frac{X}{H_1} \right) \nonumber \\
&+ \left( \frac{1}{6} (2 \log 2 + 1) + \frac{2}{3} (\log 2 - 1) + \frac{1}{3} - \frac{1}{12} \right) \textrm{deg}(F) (H_2^2 - H_1^2) X \nonumber \\
&+ \left( \frac{1}{3} + \frac{2}{3} - \frac{1}{3} - \frac{1}{6} \right) (H_2^2 - H_1^2) X (\log \mathfrak{q}(F) - \textrm{deg}(F) (\log 2 \pi + \gamma - 1)) \nonumber \\
&+ O( H_1^{2 + c} X^{1 - c}) \nonumber
\end{align}
\begin{align}
= &\frac{1}{2} \textrm{deg}(F) X \left( H_2^2 \log \frac{X}{H_2} - H_1^2 \log \frac{X}{H_1} \right) \nonumber \\
&+ \left( \left( \log 2 - \frac{1}{4} \right) \textrm{deg}(F) + \frac{1}{2} (\log \mathfrak{q}(F) - \textrm{deg}(F) (\log 2 \pi + \gamma - 1)) \right) \nonumber \\
& (H_2^2 - H_1^2) X \nonumber \\
& + O( H_1^{2 + c} X^{1 - c}) \nonumber \\
=&\frac{1}{2} \textrm{deg}(F) X \left( H_2^2 \log \frac{X}{H_2} - H_1^2 \log \frac{X}{H_1} \right) \nonumber \\
&+ \frac{1}{4} (2 \log \mathfrak{q}(F) + (2 \log 2 + 1 - 2 \log \pi - 2 \gamma) \textrm{deg}(F)) (H_2^2 - H_1^2) \nonumber \\
&+ O( H_1^{2 + c} X^{1 - c}).
\end{align}
Therefore
\begin{align}
&\sum_{1 \leq k \leq K} \int_{2^{- k} X \leq x \leq 2^{- k + 1} X, H_1 \leq h \leq H_2} |\psi(x, h; F) - m(F) h|^2 \, dh dx \nonumber \\
=&\int_{2^{- K} X \leq x \leq X} \int_{H_1 \leq h \leq H_2} |\psi(x, h; F) - m(F) h|^2 \, dh dx \nonumber \\
=&\frac{1}{2} \textrm{deg}(F) X \sum_{1 \leq k \leq K} 2^{- k} \left( H_2^2 \log \frac{X}{H_2} - H_1^2 \log \frac{X}{H_1} \right) \nonumber \\
&- \frac{1}{2} \textrm{deg}(F) X \sum_{1 \leq k \leq K} k 2^{- k} (H_2^2 - H_1^2) \log 2 \nonumber \\
&+ \frac{1}{4} (2 \log \mathfrak{q}_F + (2 \log 2 + 1 - 2 \log \pi - 2 \gamma) \textrm{deg}(F)) (H_2^2 - H_1^2) X \sum_{1 \leq k \leq K} 2^{- k} \nonumber \\
& + \sum_{1 \leq k \leq K} O(H_1^{2 + c} X^{1 - c} 2^{k(c - 1)}).
\end{align}
The sums are in the form of the geometric series given by
\begin{equation}
\sum_{1 \leq k \leq K} 2^{- k} = 1 - 2^{- K}
\end{equation}
or by the series given by
\begin{align}
&\sum_{1 \leq k \leq K} k 2^{- k} = 2 \sum_{1 \leq k \leq K} k (2^{- k} - 2^{- (k + 1)}) = 1 + 2 \sum_{2 \leq k \leq K} 2^{- k} - K 2^{- K} \nonumber \\
=&2 \sum_{1 \leq k \leq K} 2^{- k} - K 2^{- K} = 2 - (K + 2) 2^{- K}.
\end{align}
Thus
\begin{align}
&\int_{H_1 \leq h \leq H_2, 2^{- K} X \leq x \leq X} |\psi(x, h; F) - m(F) h|^2 \, dh dx \nonumber \\
=& \frac{1 - 2^{- K}}{2} \textrm{deg}(F) \left( H_2^2 \log \frac{X}{H_2} - H_1^2 \log \frac{X}{H_1} \right) X \nonumber \\
&- \textrm{deg}(F) (H_2^2 - H_1^2) \log 2 X + (K + 2) 2^{- (K + 1)} \textrm{deg}(F) (H_2^2 - H_1^2) \log 2 \nonumber \\
&+ \frac{1 - 2^{- K}}{4} (2 \log \mathfrak{q}(F) + (2 \log 2 + 1 - 2 \log \pi - 2 \gamma) \textrm{deg}(F)) (H_2^2 - H_1^2) X \nonumber \\
&+ O(H_1^{2 + c} X^{1 - c}).
\end{align}
To proceed the work in \cite{MR0480388}, now lemma $6,$ is used to produce a bound similar to that in the previous lemma, but now for a fixed difference.
\subsection{Lemma D}
\begin{equation}
Var^{fix}(X, h; F) \ll h X \left(\log \frac{X}{h}\right)^2
\end{equation}
for $0 < h \leq X.$
\subsubsection{Proof}
Considering $h \leq x / 6$ and supposing that $h \leq \tau - h \leq 2 h$ and $ x \leq u + h \leq 3 x$ then
\begin{align}
& h \int_{x \leq u \leq 2 x} |\psi(u, h; F) - m(F) h|^2 \, du \nonumber \\
\ll & \int_{x \leq u \leq 3 x, h \leq w \leq 3} |\psi(u, w; F) - m(F) w|^2 \, dw du \nonumber \\
\ll & x \int_{x \leq u \leq 3 x, h /  3 x \leq \delta \leq 3 h / x} |\psi(u, \delta u; F) - m(F) \delta u|^2 \, d\delta du \nonumber \\
\ll & x \int_{h / 3 x \leq \delta \leq 3 h / x, x \leq u \leq 3 x} |\psi(u, \delta u; F) - m(F) \delta u|^2 \, du d\delta,
\end{align}
the lemma then follows by applying Lemma C.
$\Box$

Using
\begin{equation}
\int_{H_1 \leq h \leq H_2, 1 \leq x \leq 2^{- K} X} |\psi(x, h; F) - m(F) h|^2 \, dx dh = O(H_2^2 X^{1 + \epsilon} 2^{- K})
\end{equation}
then gives
\begin{align}
&\int_{H_1 \leq h \leq H_2, 1 \leq x \leq X} |\psi(x, h; F) - m(F) h|^2 \, dh dx \nonumber \\
=& \frac{1}{2} \textrm{deg}(F) \left( H_2^2 \log \frac{X}{H_2} - H_1^2 \log \frac{X}{H_1} \right) X \nonumber \\
&+ \frac{1}{4} (2 \log \mathfrak{q}(F) + (- 2 \log 2 + 1 - 2 \log \pi - 2 \gamma) \textrm{deg}(F)) (H_2^2 - H_1^2) X \nonumber \\
&+ O(H_1^{2 + c} X^{1 - c}) + O(H_2^2X^{1 + \epsilon} 2^{- K}).
\end{align}
Setting $H_1 = H$ and $H_2 = (1 + \eta) H$ gives
\begin{align}
&\int_{H \leq h \leq (1 + \eta) H, 1 \leq x \leq X} |\psi(x, h; F) - m(F) h|^2 \, dh dx \nonumber \\
=& \frac{1}{2} \textrm{deg}(F) \left( (1 + 2 \eta + \eta^2) \left(\log \frac{X}{H} - \log (1 + \eta) \right) - \log \frac{X}{H} \right) H^2 X \nonumber \\
&+ \frac{1}{4} (2 \log \mathfrak{q}(F) + (- 2 \log 2 + 1 - 2 \log \pi - 2 \gamma) \textrm{deg}(F)) (\eta^2 + 2 \eta) H^2 X \nonumber \\
&+ O(H^{2 + c} X^{1 - c}) + O((1 + \eta)^2 H^2 X^{1 + \epsilon} 2^{- K}) \nonumber \\
=& \frac{1}{2} \textrm{deg}(F) \left((2 \eta + \eta^2) \log \frac{X}{H} \right) H^2 X \nonumber \\
&+ \frac{1}{4} (2 \log \mathfrak{q}(F) + (- 2 \log 2 + 1 - 2 \log \pi - 2 \gamma) \textrm{deg}(F)) (\eta^2 + 2 \eta) H^2 X \nonumber \\
&+ O(H^{2 + c} X^{1 - c}) + O((1 + \eta)^2 H^2 X^{1 + \epsilon} 2^{- K}) \nonumber \\
=& \eta \left(\textrm{deg}(F) \log \frac{X}{H} + \log \mathfrak{q}(F) - \left(\log 2 \pi + \gamma - \frac{1}{2} \right) \textrm{deg}(F)\right) H^2 X \nonumber \\
&+ O \left( \eta^2 H^2 X \log \frac{X}{H} \right) + O(H^{2 + c} X^{1 - c}) + O(H^2 X^{1 + \epsilon} 2^{- K}).
\end{align}
Following the same analysis starting with
\begin{equation}
Var^{mul}(X, \delta; F) = \frac{1}{6} \delta X (3 \log X - 4 \log 2) + O(\delta^{1 + c} X)
\end{equation}
gives
\begin{align}
& Var^{fix}(X, h; F) = \frac{1}{6} h (6 \log X - (3 + 8 \log 2)) \nonumber \\
& + O_{\epsilon}(h X^{\epsilon} (h / X)^{c / 3}) + O_{\epsilon}(h X^{\epsilon} (h X^{- (1 - B_1)})^{1 / (3 (1 - B_1))}).
\end{align}
\section{Numerics}
The behaviour of $Var^{fix}(X, h; F)$ is now examined numerically for some examples of $L$-functions. The Riemann zeta function, the $L$-function associated to the Ramanujan tau function and some $L$-functions associated to different elliptic curves are used.

The code used to generate the data is given in Appendix A. For the Ramanujan tau and elliptic curve cases the associated von Mangoldt functions require the coefficients in the Dirichlet series, for example for the Ramanujan tau case
\begin{equation}
\Lambda(n; \Delta) = 2 \cos\left(\frac{\log(n)}{\Lambda(n)} \arccos\left(\frac{\lambda(\exp(\Lambda(n)); \Delta)}{2}\right)\right)\Lambda(n).
\end{equation}
The unnormalised coefficients were generated using Sage, \cite{Sage}. For the Ramanujan tau case a normalisation of $n^{11 / 2}$ was then applied to each coefficient to give the $n$th normalised coefficient, $\lambda(n; \Delta).$ For the elliptic curve case there is a square root normalisation, otherwise the von Mangoldt function is the same except at a finite number of points where primes divide the conductor. The elliptic curve is specified by entering its Weierstrass coefficients, $(a_1, a_2, a_3, a_4, a_6),$ which were found on the database of $L$-functions, modular forms, and related objects, \cite{lmfdb}, $L$-functions in subsequent figures and tables are referred to via their label given on the database, which is called the LMFDB label there.
\begin{figure}[h]
\includegraphics[scale=0.33]{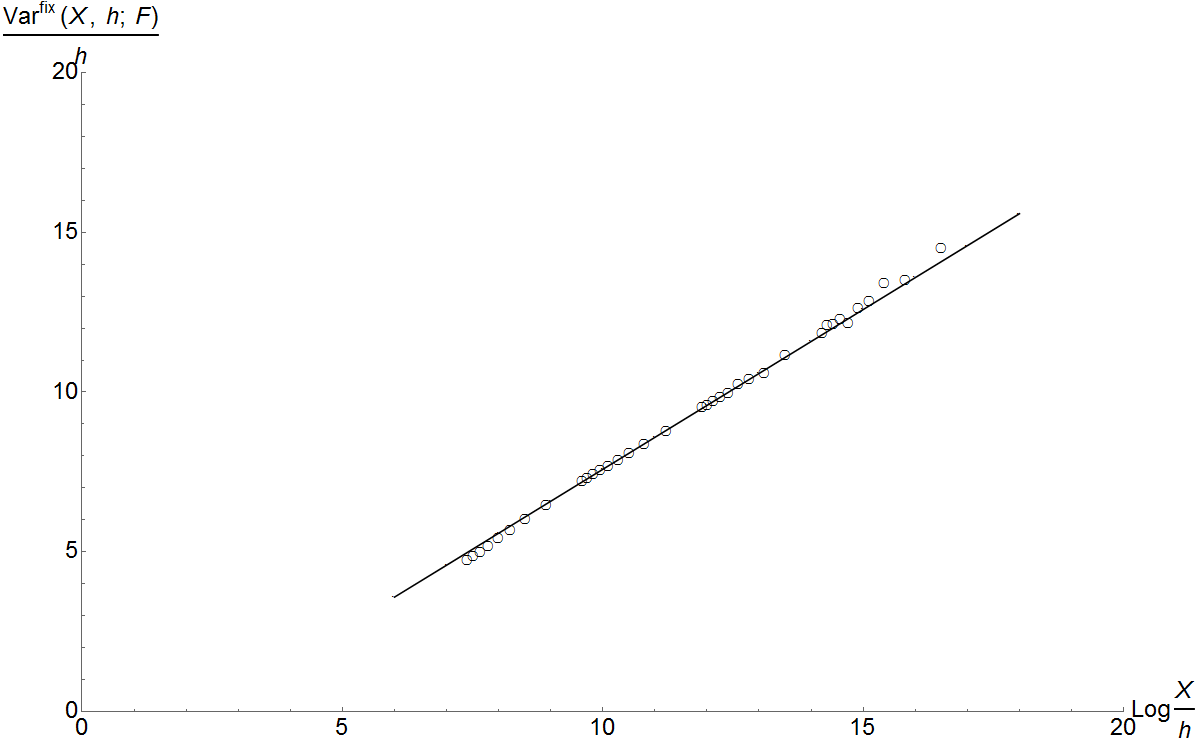} \caption[Variance for the Riemann zeta function case for $X = 15 000 000$]{For the Riemann zeta function case $X$ is taken to be $15 000 000,$ the straight line is given by $\log \frac{X}{h} - (\log 2 \pi + \gamma).$}
\end{figure}
\begin{figure}[h]
\includegraphics[scale=0.33]{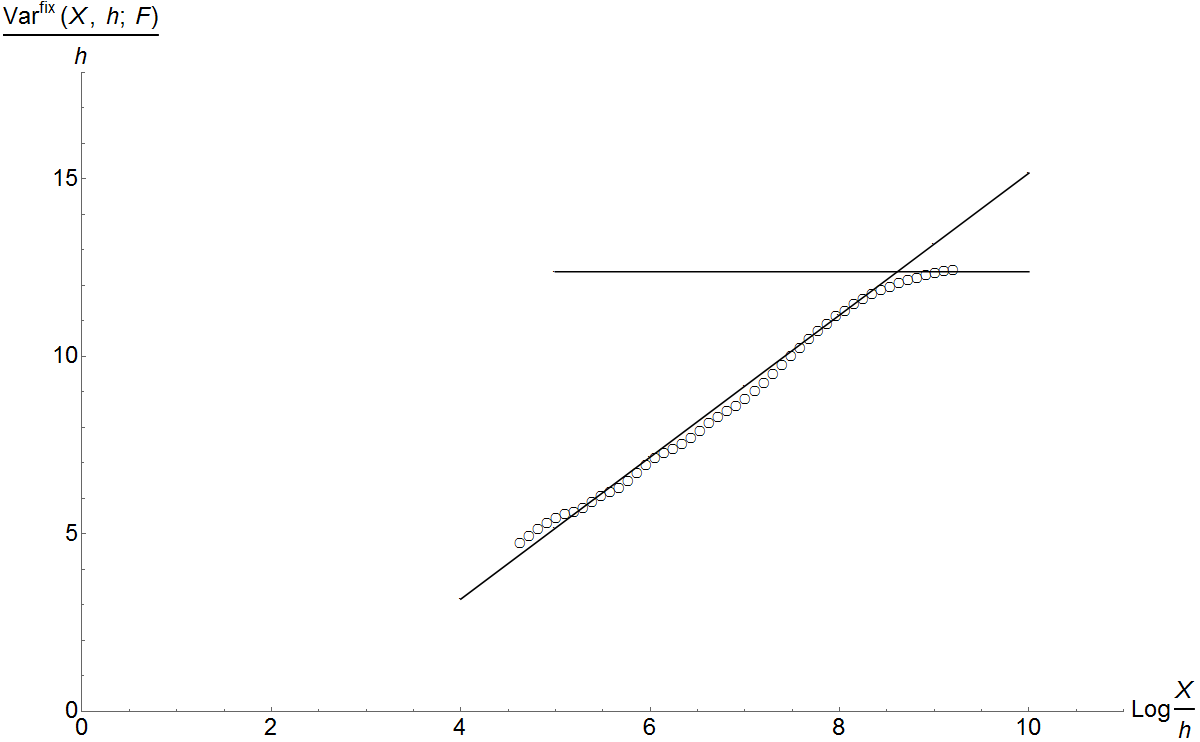} \caption[Variance for the Ramanujan tau function case for $X = 1 000 000$]{This plot shows the appearence of $2$ regimes, illustrated by plotting data for the Ramanujan tau case with $X = 1 000 000.$ The horizontal line is given by $\log 1 000 000 - \frac{1}{6} (3 + 8 \log 2),$ the slanted line is given by $2 \log \frac{X}{h} - 2 (\log 2 \pi + \gamma).$}
\end{figure}
\begin{figure}[h]
\includegraphics[scale=0.33]{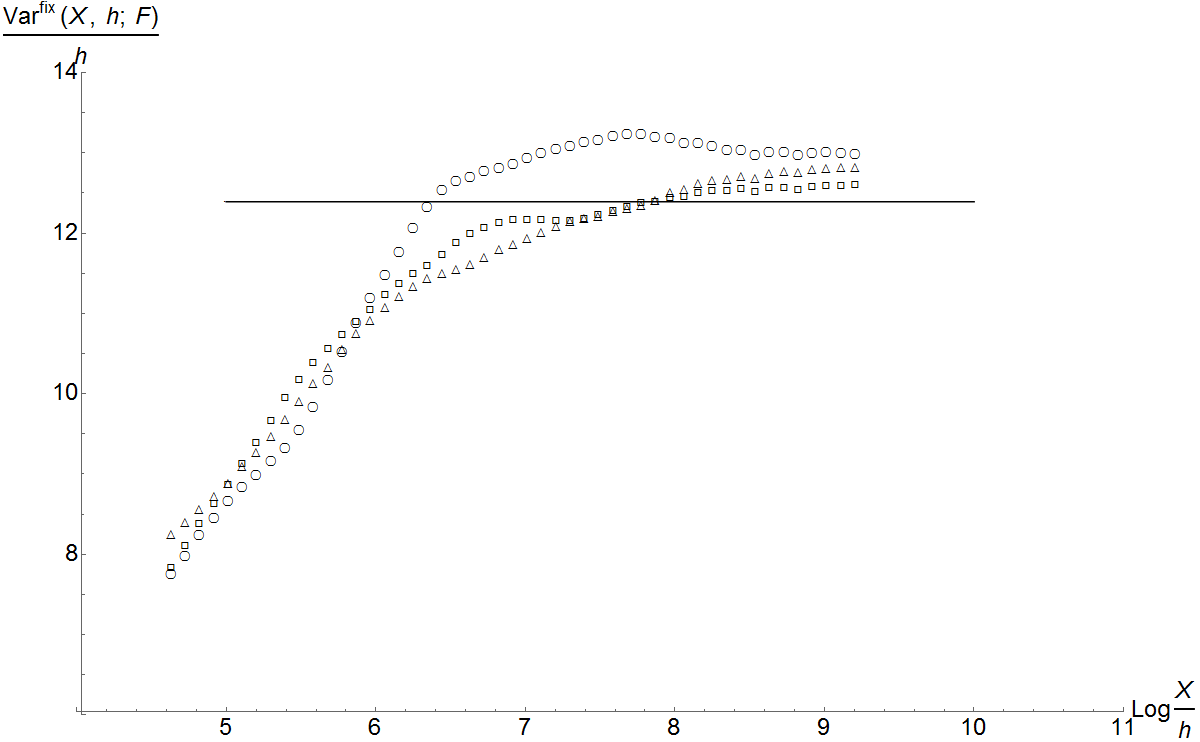} \caption[Variances for the elliptic curve cases E27.a1, E37.a1 and E37.b2 for $X = 1 000 000$]{The universal behaviour in the regime given by the horizontal black line is further illustrated by plotting data for $3$ elliptic curves, $\bigcirc$ are for an elliptic curve of conductor $27$, E27.a1, $\vartriangle$ are for an elliptic curve of conductor $37,$ E37.a1, and the $\square$ points correspond to another conductor $37$ elliptic curve, E37.b2. In each case X = 1 000 000.}
\end{figure}

\chapter{Outlook}
\section{Random matrix theory and \\ higher order correlations}
The universal forms of the $2$-point correlation statistic in both the quantum chaology and analytic number theory setting are both conjectural, however they coincide with the corresponding statistic for random unitary matrices with normally distributed matrix elements, a theorem in this field of study. Random matrix theory is a rapidly developing area of study in both mathematics and physics, with its origins in the study of the multiple interactions present in atomic nuclei. A thorough treatment of the subject from a physical perspective may be found in the textbook \cite{MR2129906}, for more mathematical emphasis \cite{MR2760897} is recommended.

The connection between quantum chaology, analytic number theory and random matrix theory extends beyond $2$-point statistics. In \cite{MR2540409} the rescaled $n$-point correlation statistic was shown to coincide with the corresponding statistic in random matrix theory, given by the determinant
\begin{equation} \label{Random matrix theory correlations}
\textrm{det}_{n \textrm{ x } n} S(x_i - x_j), \textrm{ with } i, j = 1, 2, \ldots, n,
\end{equation}
where
\begin{equation}
S(x) = \frac{\sin \pi x}{\pi x}.
\end{equation}
Following Montgomery's study of the pair correlation of zeros of the Riemann zeta function through $\mathcal{F}(X, T),$ in \cite{MR1283025} Hejhal studied triple correlations of zeros of the Riemann zeta function and subsequently Rudnick and Sarnak considered the general case of $n$-point correlations for $L$-functions associated to automorphic forms in \cite{MR1395406}. From a physical perspective, in a pair of papers Bogomolny and Keating first treated the case of triple and quadruple correlations of zeros of the Riemann zeta function around the time of the work of Hejhal in \cite{MR1363402} and then the $n$-point case in \cite{MR1399479}. In all cases there is a match with the random matrix theory expression \eqref{Random matrix theory correlations}.

The results in this thesis concern $2$-point statistics, however assuming the form \eqref{Random matrix theory correlations}, it would be interesting to determine moments beyond the variance determined in this thesis, generalising the work of Montgomery and Soundararajan in \cite{MR2104891}, taking zero correlations rather than prime correlations as the starting point.
\section{Hardy-Littlewood type conjectures}
In this thesis $Var^{fix}(X, h; F)$ was calculated by starting with a conjecture for the distribution of non-trivial zeros of $F,$ however the variance may be related to the generalised correlation function, 
\begin{equation} \label{Hardy-Littlewood type generalised correlation function}
C(r, N; F) = \frac{1}{N} \sum_{n \leq N} \Lambda(n + r; F) \Lambda^*(n; F),
\end{equation}
and was determined from the Hardy-Littlewood conjecture by Montgomery and Soundararajan for the case of the Riemann zeta function in \cite{MR2104891} in this way. Writing $Var^{fix}(X, h; F)$ as a sum,
\begin{equation}
Var^{fix}(X, h; F) = \frac{1}{X} \sum_{1 \leq x \leq X} \left|\sum_{x < n \leq x + h} \Lambda(n; F) - m(F)h\right|^2,
\end{equation}
gives
\begin{align}
\frac{1}{X} \sum_{1 \leq x \leq X} \sum_{x < n \leq x + h} \left|\Lambda(n; F)\right|^2 + \frac{1}{X} \sum_{1 \leq x \leq X} \sum_{x < m \neq n \leq x + h} \Lambda(m; F) \Lambda^*(n; F) \nonumber \\
- \frac{m(F) h}{X} \sum_{1 \leq x \leq X} \left(\sum_{x < m \leq x + h} \Lambda(m; F) + \sum_{x < n \leq x + h} \Lambda^*(n; F)\right) + (m(F) h)^2,
\end{align}
and so an understanding of the Hardy-Littlewood type generalised correlation function, \eqref{Hardy-Littlewood type generalised correlation function} over short intervals of length $h,$ averaged over $X$ is desired.

The inversion argument in chapter $2$ is well suited to the Riemann zeta function case, where the products involved in the off-diagonal contribution to the $2$-point correlation statistic may be written in terms of Ramanujan sums. If such a decomposition of the products occuring in the general expression for the $2$-point statistic may be found then the method could be replicated for the general case.

Plots for the Riemann zeta function case, the Ramanujan tau function case and an elliptic curve case are given for comparison.
\begin{figure}[h]
\includegraphics[scale=0.33]{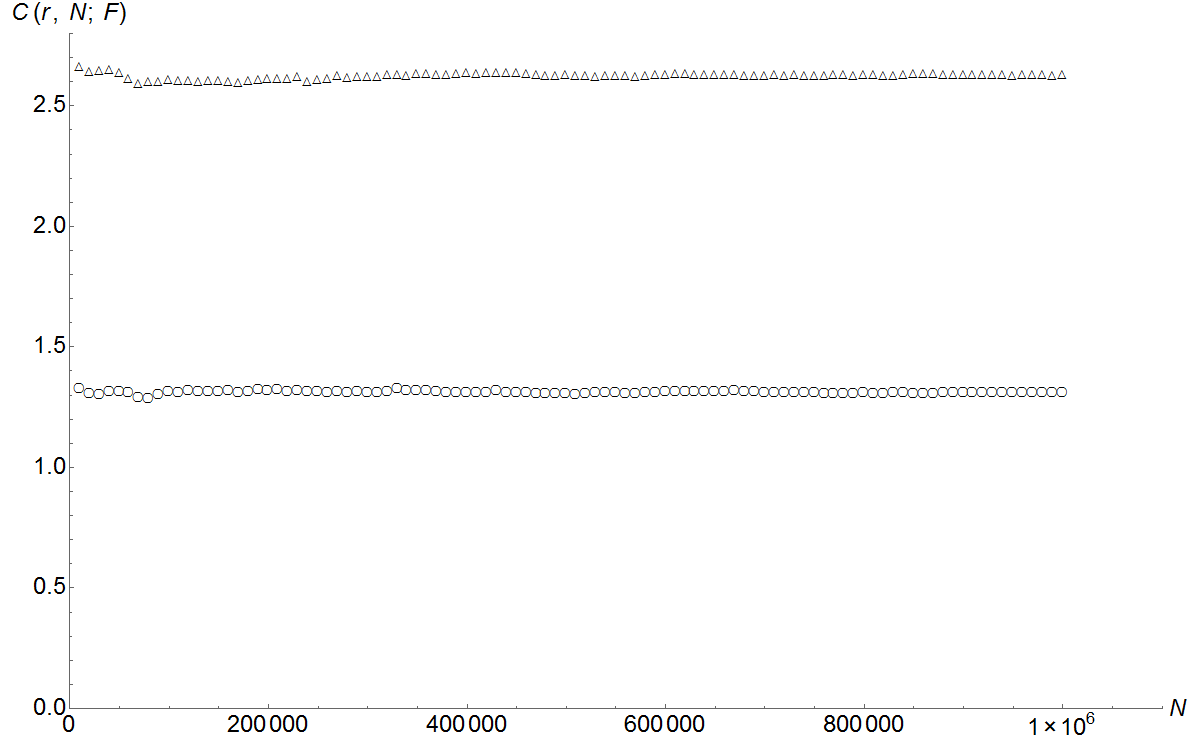} \caption[Arithmetic correlations for the Riemann zeta function case]{Correlations for the Riemann zeta function case, $\bigcirc$ are for $r = 2$ and $\vartriangle$ are for $r = 6,$ in agreement with the Hardy-Littlewood conjecture.}
\end{figure}
\begin{figure}
\includegraphics[scale=0.33]{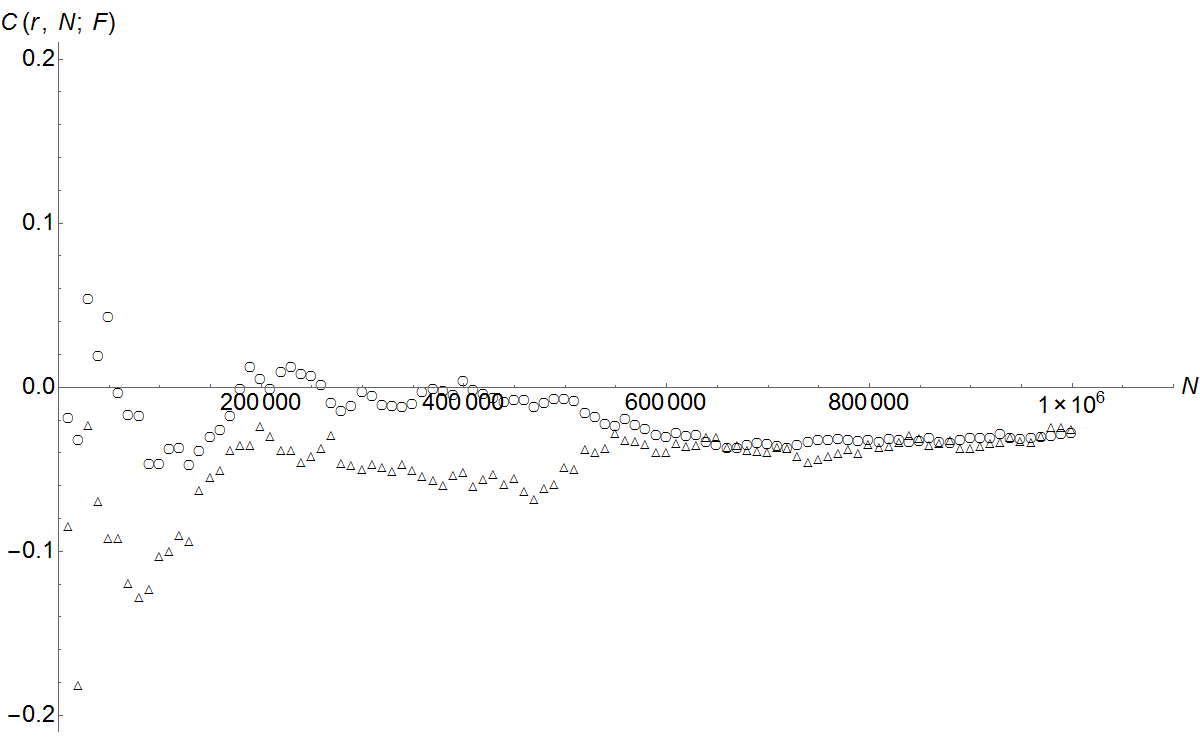} \caption[Arithmetic correlations for the Ramanujan tau function case]{Correlations for the Ramanujan tau function case, $\bigcirc$ are for $r = 2$ and $\vartriangle$ are for $r = 6.$ A negative bias is apparent with the appearence of convergence to a small negative value. This may change for larger values of $N,$ however the code used would take an unfeasible length of time to extend $N$ considerably.}
\end{figure}

\begin{figure}
\includegraphics[scale=0.33]{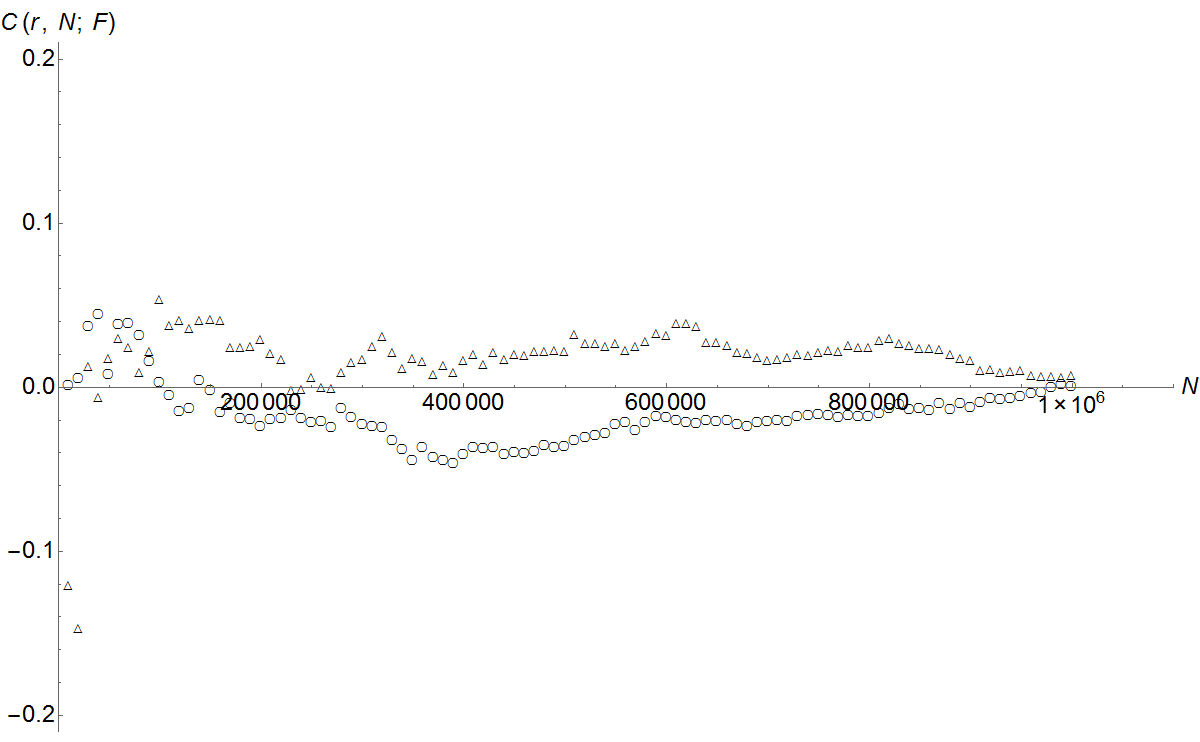} \caption[Arithmetic correlations for the elliptic curve case E11.a3]{Correlations for the elliptic curve E11.a3, $\bigcirc$ are for $r = 2$ and $\vartriangle$ are for $r = 6.$ No correlations are apparent in this case.}
\end{figure}

\appendix

\chapter{Code for numerics}
\section{Python code for Sage}
\subsection{Generating Ramanujan tau values}
sage: d = CuspForms(Gamma0(1), 12).basis()[0] \\
sage: t = [0] \\
sage: t = t + d.coefficients(1010000) \\
sage: tp = [] \\
sage: for p in primes(1010001): \\
....:   tp.append(t[p]) \\
....:
\subsection{Writing Ramanujan tau values to a file}
sage: o = open(`File', `w') \\
sage: for i in range(0, len(tp)): \\
....:   o.write(str(tp[i]) + `$\backslash$n') \\
....: \\
sage: o.close()
\subsection{Generating coefficients associated to an elliptic curve}
sage: E = EllipticCurve(($a_1$, $a_2$, $a_3$, $a_4$, $a_6$)) \\
sage: x = E.aplist(1010000)
\subsection{Writing coefficients associated to an elliptic curve to a file}
sage: o = open(`File', `w') \\
sage: for i in range(0, len(x)): \\
....:   o.write(str(x[i]) + `$\backslash$n') \\
....: \\
sage: o.close()
\section{Java code}
\subsection{Class for the Riemann zeta function and \\ Ramanujan tau function}
import java.io.FileOutputStream; \\
import java.io.PrintStream; \\
import java.io.FileReader; \\
import java.io.IOException; \\
import java.util.Scanner; \\ \\
public class Tau \{ \\
public static final int[] Prime = MaximumValue(1010000); \\
private static int[] MaximumValue(int max) \{ \\
boolean[] isComposite = new boolean[max + 1]; \\
for (int i = 2; i * i $<$ = max; i++) \{ \\
if (!isComposite [i]) \{ \\
for (int j = i; i * j $<$ = max; j++) \{ \\
isComposite [i * j] = true; \\
\} \\
\} \\
\} \\
int numPrimes = 0; \\
for (int i = 2; i $<$ = max; i++) \{ \\
if (!isComposite [i]) numPrimes++; \\
\} \\
int index = 0; \\
int [] Prime = new int [numPrimes]; \\
for (int i = 2; i $<$ = max; i++) \{ \\
if (!isComposite [i]) Prime [index++] = i; \\
\} \\
return Prime; \\
\} \\ \\
private static int log (int Base, int Argument) \{ \\
return (int) (Math.log(Argument) / Math.log(Base)); \\
\} \\ \\
private static double vonMangoldtfunction(int n, double[] a) \{ \\
return a[n - 1]; \\
\} \\ \\
private static double nprime(int n, double [] a)\{ \\
return Math.round(Math.pow(Math.E, vonMangoldtfunction(n, a))); \\
\} \\ \\
private static double nprimepower(int n, double [] a) \{ \\
return (Math.log(n))/(vonMangoldtfunction(n, a)); \\
\} \\ \\
private static double Taunormalisedcoefficient(int n, double[] a) \{  \\
if (a[n - 1] != 0) \{ \\
return (a[n - 1])/(Math.pow(n, 5.5)); \\
\} \\
else return 0; \\
\} \\ \\
private static double TauvonMangoldtfunction(int n, double[] a, double[] b) \{ \\
if (vonMangoldtfunction(n, a) != 0) \{ \\
return 2*(Math.cos((nprimepower(n, a)) \\
*(Math.acos((Taunormalisedcoefficient((int) nprime(n, a), b)) / 2)))) \\
*(vonMangoldtfunction(n, a)); \\
\} \\
else return 0; \\
\} \\ \\
private static double intervalsum(double x, double y, double[] a, double [] b) \{ \\
double value = 0.0; \\
for (int i = (int)(x + 1); i $<$ = x + y; i++) \{ \\
value = value + TauvonMangoldtfunction(i, a, b); \\
\} \\
return value; \\
\} \\ \\
private static double squaremodulusofintervalsum(double x, double y, double[] a, double [] b) \{ \\
return (intervalsum(x, y, a, b))*(intervalsum(x, y, a, b)); \\
\} \\ \\
private static double variance(double x, double y, double[] a, double [] b) \{ \\
double value = 0.0; \\
for(int i = 1 ; i $<$ = x; i++)\{ \\
value = value + squaremodulusofintervalsum(i, y, a, b); \\
\} \\
return value / (x * y); \\
\} \\ \\
private static double correlationfunction(int r, int n, double[] a, double[] b) \{ \\
double sum = 0; \\
for(int i = 1; i $<$ =  n; i++)\{ \\
sum = sum + (vonMangoldtfunction(i + r, a)) * (vonMangoldtfunction(i, a)); \\
\} \\
return sum / n; \\
\} \\ \\
// ----- Main Method ----- // \\ \\
public static void main(String[] args) throws IOException \{ \\
double[] vonMangoldt = new double[1010000]; \\ \\
for (int i = 0; i $<$ Prime.length; i++)\{ \\
for (int j = 1; j $<$ = log(Prime[i], 1010000); j++)\{ \\
vonMangoldt[(int)(Math.pow(Prime[i], j) - 1)] = Math.log(Prime[i]); \\
\} \\
\} \\ \\
FileReader file = new FileReader(``File directory"); \\ \\
double[] Tauatprimearray = new double[79251]; \\ \\
int j = 0; \\
try \{\\
Scanner input = new Scanner(file); \\
while(input.hasNext()) \\
\{ \\
Tauatprimearray[j] = input.nextDouble(); \\
j++; \\
\} \\
input.close(); \\
\} \\ \\
catch(Exception e) \\
\{ \\
e.printStackTrace(); \\
\} \\ \\
double[] Tauprimearray = new double[1010000]; \\
for(int i = 0; i $<$ Prime.length; i++) \{ \\
Tauprimearray[Prime[i] - 1] = Tauatprimearray[i]; \\
\} \\
\} \\	
\}
\subsection{Writing variance data and correlation data to a file}
PrintStream out = new PrintStream(new FileOutputStream(``File.txt")); \\
System.setOut(out); \\ \\ 
for(double h = 50.0; h $<$ = 10000.0; h = 1.1 * h) \{ \\
System.out.println( \\
Math.log(1000000 / h) + `` " + variance(1000000, h, vonMangoldtfunction, Tauprimearray) \\
); \\
\} \\ \\
for(int n = 1; n $<$ = 1000000; n++) \{ \\
if(n \% 10000 == 0)System.out.println( \\
n + `` " + (correlationfunction(r, vonMangoldtfunction, Tauprimearray)) \\
); \\
\}

\chapter{Variance data}
\begin{table}[h]
\begin{center}
\begin{tabular}{| c | c |} \hline
$\log X / h$&$Var^{fix}(X, h; F) / h$ \\ \hline
16.523560759066484&14.520951259388584 \\
15.830413578506539&13.521030326116232 \\
15.424948470398375&13.403677020009642 \\
15.137266397946593&12.845041532187183 \\
14.914122846632385&12.6383584434893 \\
14.73180128983843&12.167262959003123 \\
14.57765061001117&12.29675543048168 \\
14.444119217386648&12.143894781372447 \\ 
14.326336181730264&12.096567426173788 \\ 
14.22097566607244&11.858721924525783 \\ 
13.527828485512494&11.170738706704546 \\ 
13.122363377404328&10.592076902033908 \\ 
12.834681304952548&10.421986805840097 \\
12.611537753638338&10.237567907703243 \\ 
12.429216196844383&9.969909185902116 \\
12.275065517017126&9.848096278758504 \\ \hline
\end{tabular} \caption{Variance data for the Riemann zeta function case for $X = 15 000 000$}
\end{center}
\end{table}
\begin{table}[h]
\begin{center}
\begin{tabular}{| c | c |} \hline
$\log X / h$&$Var^{fix}(X, h; F) / h$ \\ \hline
12.141534124392603&9.721962958385818 \\
12.02375108873622&9.578118995789781 \\ 
11.918390573078392&9.531382786577403 \\
11.225243392518447&8.792654834154554 \\
10.819778284410283&8.382459283873876 \\
10.532096211958502&8.091839805669087 \\
10.308952660644293&7.875314741751184 \\
10.126631103850338&7.692603652059396 \\
9.972480424023079&7.552222846376847 \\
9.838949031398556&7.427696970639272 \\
9.721165995742174&7.312241880462279 \\
9.615805480084347&7.211462585733127 \\
8.922658299524402&6.461382646993843 \\
8.517193191416238&6.016651372253434 \\
8.229511118964457&5.691277480520871 \\
8.006367567650246&5.42447154063884 \\
7.824046010856292&5.190747452086452 \\
7.669895331029034&5.005145706015339 \\
7.536363938404511&4.8543217618304055 \\
7.418580902748128&4.734743947320245 \\ \hline
\end{tabular} \caption{Variance data for the Riemann zeta function case for $X = 15 000 000$}
\end{center}
\end{table}

\begin{table}[h]
\begin{center}
\begin{tabular}{| c | c |} \hline
$\log X / h$&$Var^{fix}(X, h; F) / h$ \\ \hline
9.210340371976184&12.450145500729484 \\ 
9.115030192171858&12.400661530125248 \\
9.019720012367532&12.34708451164469 \\ 
8.924409832563208&12.284762216995691 \\ \hline
\end{tabular} \caption{Variance data for the Ramanujan tau function case for $X = 1 000 000$}
\end{center}
\end{table}
\begin{table}[h]
\begin{center}
\begin{tabular}{| c | c |} \hline
$\log X / h$&$Var^{fix}(X, h; F) / h$ \\ \hline
8.829099652758883&12.197874687361988 \\
8.733789472954557&12.164159618941683 \\
8.638479293150233&12.082616628550099 \\
8.543169113345908&11.958999153859036 \\
8.447858933541584&11.885484133822498 \\
8.352548753737258&11.751711045091765 \\
8.257238573932932&11.628142799666147 \\
8.161928394128608&11.476177630280777 \\
8.066618214324283&11.292865565851903 \\
7.971308034519958&11.13438685471049 \\
7.875997854715633&10.919082193463444 \\
7.7806876749113085&10.710582242779699 \\
7.685377495106984&10.485899640658461 \\
7.590067315302659&10.252602524429552 \\ 
7.494757135498334&10.004228907396628 \\
7.399446955694009&9.768845004850796 \\
7.304136775889684&9.51720291528726 \\
7.208826596085359&9.26350592739765 \\
7.113516416281034&9.040880153862338 \\ 
7.0182062364767095&8.812677793267488 \\
6.922896056672384&8.609000170727947 \\
6.827585876868059&8.451093480225197 \\
6.732275697063734&8.288420442037078 \\ 
6.6369655172594095&8.11362657142882 \\ 
6.541655337455085&7.906395737347605 \\
6.446345157650759&7.6999278055609786 \\ 
6.351034977846434&7.532567145504188 \\
6.25572479804211&7.397706971472282 \\ \hline
\end{tabular} \caption{Variance data for the Ramanujan tau function case for $X = 1 000 000$}
\end{center}
\end{table}
\begin{table}[h]
\begin{center}
\begin{tabular}{| c | c |} \hline
$\log X / h$&$Var^{fix}(X, h; F) / h$ \\ \hline
6.160414618237785&7.27547302341981 \\
6.06510443843346&7.132153930099307 \\
5.969794258629135&6.941217742936219 \\
5.87448407882481&6.709297941076495 \\
5.779173899020485&6.496911012764781 \\
5.68386371921616&6.303311826157115 \\
5.588553539411835&6.174019590689406 \\
5.4932433596075105&6.072318960136859 \\
5.397933179803185&5.901837491254716 \\
5.30262299999886&5.73333600100444 \\
5.207312820194535&5.62867540353775 \\
5.112002640390211&5.549964378657342 \\
5.016692460585886&5.441854148589055 \\
4.92138228078156&5.3061800338739395 \\
4.826072100977235&5.127159813303661 \\
4.730761921172911&4.95365755484686 \\
4.635451741368586&4.740360864001797 \\ \hline
\end{tabular} \caption{Variance data for the Ramanujan tau function case for $X = 1 000 000$}
\end{center}
\end{table}

\begin{table}[h]
\begin{center}
\begin{tabular}{| c | c |} \hline
$\log X / h$&$Var^{fix}(X, h; F) / h$ \\ \hline
9.210340371976184&12.988785064767487 \\
9.115030192171858&12.997764032878326 \\
9.019720012367532&13.004586724380916 \\
8.924409832563208&12.996769241593741 \\
8.829099652758883&12.969625627496745 \\ 
8.733789472954557&13.00376857091496 \\
8.638479293150233&13.007987156980438 \\ \hline
\end{tabular} \caption{Variance data for the elliptic curve case E27.a1 for $X = 1 000 000$}
\end{center}
\end{table}
\begin{table}[h]
\begin{center}
\begin{tabular}{| c | c |} \hline
$\log X / h$&$Var^{fix}(X, h; F) / h$ \\ \hline
8.543169113345908&12.975635494444667 \\
8.447858933541584&13.033284261301656 \\
8.352548753737258&13.033043192354626 \\
8.257238573932932&13.08258479534229 \\
8.161928394128608&13.115676961833989 \\
8.066618214324283&13.122741336640091 \\
7.971308034519958&13.18214683047017 \\
7.875997854715633&13.20063183388228 \\
7.7806876749113085&13.228130275666208 \\
7.685377495106984&13.229483400369974 \\
7.590067315302659&13.205322374563357 \\
7.494757135498334&13.157477154742308 \\
7.399446955694009&13.136853514795833 \\
7.304136775889684&13.083083808143206 \\
7.208826596085359&13.051777633357279 \\
7.113516416281034&12.995012253080189 \\
7.0182062364767095&12.929632471752942 \\
6.922896056672384&12.857957169486875 \\
6.827585876868059&12.815959430094711 \\
6.732275697063734&12.773318932555258 \\
6.6369655172594095&12.702737753973945 \\ 
6.541655337455085&12.642776486984614 \\
6.446345157650759&12.531318431782148 \\
6.351034977846434&12.321691986578813 \\ 
6.25572479804211&12.057288145267053 \\ 
6.160414618237785&11.767987859842052 \\
6.06510443843346&11.47497952479035 \\ 
5.969794258629135&11.18894511351632 \\ \hline
\end{tabular} \caption{Variance data for the elliptic curve case E27.a1 for $X = 1 000 000$}
\end{center}
\end{table}
\begin{table}[h]
\begin{center}
\begin{tabular}{| c | c |} \hline
$\log X / h$&$Var^{fix}(X, h; F) / h$ \\ \hline
5.87448407882481&10.878437827450739 \\ 
5.779173899020485&10.519282456130739 \\
5.68386371921616&10.165806048054563 \\
5.588553539411835&9.834006161915957 \\
5.4932433596075105&9.547996149220548 \\
5.397933179803185&9.318551329935705 \\
5.30262299999886&9.161481632716423 \\
5.207312820194535&8.991573418986622 \\
5.112002640390211&8.833917397378329 \\
5.016692460585886&8.663093726005181 \\
4.92138228078156&8.450088627208237 \\
4.826072100977235&8.239267016764318 \\
4.730761921172911&7.976739439896056 \\
4.635451741368586&7.753322167924815 \\ \hline
\end{tabular} \caption{Variance data for the elliptic curve case E27.a1 for $X = 1 000 000$}
\end{center}
\end{table}

\begin{table}[h]
\begin{center}
\begin{tabular}{| c | c |} \hline
$\log X / h$&$Var^{fix}(X, h; F) / h$ \\ \hline
9.210340371976184&12.82419062430585 \\
9.115030192171858&12.824246155795045 \\
9.019720012367532&12.81357460024457 \\
8.924409832563208&12.797500787698318 \\
8.829099652758883&12.76300497079633 \\
8.733789472954557&12.776251574611361 \\
8.638479293150233&12.749465174239809 \\
8.543169113345908&12.688138502913288 \\
8.447858933541584&12.708073672744717 \\
8.352548753737258&12.668437222504767 \\ \hline
\end{tabular} \caption{Variance data for the elliptic curve case E37.a1 for $X = 1 000 000$}
\end{center}
\end{table}
\begin{table}[h]
\begin{center}
\begin{tabular}{| c | c |} \hline
$\log X / h$&$Var^{fix}(X, h; F) / h$ \\ \hline
8.257238573932932&12.659516631186529 \\ 
8.161928394128608&12.62036080103507 \\
8.066618214324283&12.55357581884703 \\
7.971308034519958&12.508511156031718 \\
7.875997854715633&12.413449182303953 \\
7.7806876749113085&12.355566368813387 \\
7.685377495106984&12.316274278741968 \\
7.590067315302659&12.27403213804591 \\
7.494757135498334&12.218381134327782 \\
7.399446955694009&12.192096796383874 \\
7.304136775889684&12.148118791048645 \\
7.208826596085359&12.090486660647485 \\
7.113516416281034&12.014232065244181 \\
7.0182062364767095&11.943484322389757 \\
6.922896056672384&11.868757249267126 \\
6.827585876868059&11.8063365856975 \\
6.732275697063734&11.703796111556624 \\
6.6369655172594095&11.611703390592691 \\ 
6.541655337455085&11.554883462186949 \\ 
6.446345157650759&11.501008379721927 \\
6.351034977846434&11.4462613335653 \\
6.25572479804211&11.344293582088966 \\
6.160414618237785&11.220498150395324 \\ 
6.06510443843346&11.076834588986078 \\
5.969794258629135&10.91191153989503 \\
5.87448407882481&10.760513116812575 \\ 
5.779173899020485&10.554304781334274 \\
5.68386371921616&10.336733819880877 \\ \hline
\end{tabular} \caption{Variance data for the elliptic curve case E37.a1 for $X = 1 000 000$}
\end{center}
\end{table}
\begin{table}[h]
\begin{center}
\begin{tabular}{| c | c |} \hline
$\log X / h$&$Var^{fix}(X, h; F) / h$ \\ \hline
5.588553539411835&10.128055692162041 \\
5.4932433596075105&9.902032018560478 \\
5.397933179803185&9.687620146916695 \\
5.30262299999886&9.469430820236642 \\
5.207312820194535&9.27655248794704 \\
5.112002640390211&9.096138728539657 \\
5.016692460585886&8.884950651380095 \\
4.92138228078156&8.719455650136204 \\
4.826072100977235&8.56622087360894 \\
4.730761921172911&8.39938554053581 \\
4.635451741368586&8.25496743981971 \\ \hline
\end{tabular} \caption{Variance data for the elliptic curve case E37.a1 for $X = 1 000 000$}
\end{center}
\end{table}

\begin{table}[h]
\begin{center}
\begin{tabular}{| c | c |} \hline
$\log X / h$&$Var^{fix}(X, h; F) / h$ \\ \hline
9.210340371976184&12.586858854104243 \\
9.115030192171858&12.578704232496865 \\
9.019720012367532&12.571410728076517 \\
8.924409832563208&12.556088129168865 \\
8.829099652758883&12.52576463119155 \\
8.733789472954557&12.551634563137636 \\
8.638479293150233&12.54941961187593 \\
8.543169113345908&12.498020071454752 \\
8.447858933541584&12.531953259695925 \\
8.352548753737258&12.505562414336152 \\
8.257238573932932&12.512442904068271 \\ 
8.161928394128608&12.488338228410914 \\
8.066618214324283&12.437665004774983 \\ \hline
\end{tabular} \caption{Variance data for the elliptic curve case E37.b2 for $X = 1 000 000$}
\end{center}
\end{table}
\begin{table}[h]
\begin{center}
\begin{tabular}{| c | c |} \hline
$\log X / h$&$Var^{fix}(X, h; F) / h$ \\ \hline
7.971308034519958&12.421654549038788 \\
7.875997854715633&12.378682364579962 \\
7.7806876749113085&12.361274149552548 \\
7.685377495106984&12.31829396881534 \\
7.590067315302659&12.263657764514718 \\
7.494757135498334&12.209222759576027 \\
7.399446955694009&12.165479862012774 \\
7.304136775889684&12.140385311829222 \\
7.208826596085359&12.141050370807301 \\
7.113516416281034&12.144890118078353 \\ 
7.0182062364767095&12.153719294927466 \\
6.922896056672384&12.146470212662873 \\
6.827585876868059&12.118148317121461 \\
6.732275697063734&12.054966043010904 \\
6.6369655172594095&11.97288737999676 \\
6.541655337455085&11.859037270047974 \\ 
6.446345157650759&11.717806876316573 \\
6.351034977846434&11.579971995721666 \\
6.25572479804211&11.47250888884662 \\ 
6.160414618237785&11.348409028912807 \\
6.06510443843346&11.220901062521488 \\
5.969794258629135&11.033215057792335 \\
5.87448407882481&10.878399679257777 \\ 
5.779173899020485&10.723064374786576 \\
5.68386371921616&10.538634686653987 \\
5.588553539411835&10.36756358226578 \\
5.4932433596075105&10.154221239636886 \\ 
5.397933179803185&9.929323485780593 \\ \hline
\end{tabular} \caption{Variance data for the elliptic curve case E37.b2 for $X = 1 000 000$}
\end{center}
\end{table}
\begin{table}[h]
\begin{center}
\begin{tabular}{| c | c |} \hline
$\log X / h$&$Var^{fix}(X, h; F) / h$ \\ \hline
5.30262299999886&9.644831554416621 \\
5.207312820194535&9.373321479394649 \\
5.112002640390211&9.106265255379748 \\
5.016692460585886&8.855225363911499 \\
4.92138228078156&8.617284767730554 \\
4.826072100977235&8.359511939303564 \\
4.730761921172911&8.08998661024203 \\
4.635451741368586&7.81696731174195 \\ \hline
\end{tabular} \caption{Variance data for the elliptic curve case E37.b2 for $X = 1 000 000$}
\end{center}
\end{table}

\chapter{Correlation data}

\begin{table}[h]
\begin{center}
\begin{tabular}{| c | c | c |} \hline
$N$&$\frac{1}{N} \sum_{n \leq N} \Lambda(n + 2; F) \Lambda^*(n; F)$&$\frac{1}{N} \sum_{n \leq N} \Lambda(n + 6; F) \Lambda^*(n; F)$ \\ \hline
10000&1.330292434&2.665065983 \\
20000&1.309392821&2.644274684 \\
30000&1.303508484&2.646920413 \\ 
40000&1.317547846&2.65332364 \\
50000&1.318375849&2.638610509 \\
60000&1.311938543&2.616507661 \\
70000&1.290404234&2.593017786 \\
80000&1.289646866&2.601933661 \\
90000&1.303001315&2.602631725 \\
100000&1.315229125&2.612901024 \\
110000&1.312343713&2.608655735 \\ 
120000&1.321271774&2.60528435 \\
130000&1.318706331&2.604051205 \\
140000&1.315203744&2.607447111 \\
150000&1.315126018&2.607333763 \\
160000&1.322178079&2.604090594 \\ \hline
\end{tabular} \caption{Arithmetic correlations for the Riemann zeta function case}
\end{center}
\end{table}
\begin{table}[h]
\begin{center}
\begin{tabular}{| c | c | c |} \hline
$N$&$\frac{1}{N} \sum_{n \leq N} \Lambda(n + 2; F) \Lambda^*(n; F)$&$\frac{1}{N} \sum_{n \leq N} \Lambda(n + 6; F) \Lambda^*(n; F)$ \\ \hline
170000&1.313315233&2.600771486 \\
180000&1.317717831&2.606094009 \\ 
190000&1.324211639&2.611799703 \\
200000&1.321438539&2.616957348 \\
210000&1.325138356&2.614514419 \\
220000&1.317339571&2.615611765 \\
230000&1.321479779&2.624054105 \\
240000&1.31802678&2.604904809 \\
250000&1.317973033&2.61035921 \\
260000&1.311987976&2.613364256 \\
270000&1.315378805&2.626899689 \\
280000&1.31461547&2.619509899 \\
290000&1.316062284&2.624911455 \\
300000&1.311057954&2.624646674 \\ 
310000&1.314718972&2.625023088 \\
320000&1.31721448&2.631384133 \\ 
330000&1.327815562&2.632920876 \\
340000&1.323026563&2.627828728 \\
350000&1.321228085&2.637300662 \\
360000&1.320974327&2.6341375 \\
370000&1.317408426&2.632321027 \\
380000&1.31273752&2.632927121 \\
390000&1.31215607&2.635536588 \\
400000&1.313809077&2.640219199 \\
410000&1.313071724&2.635603435 \\
420000&1.313681062&2.641790665 \\ 
430000&1.319656014&2.6416832 \\
440000&1.314557655&2.640253456 \\ \hline
\end{tabular} \caption{Arithmetic correlations for the Riemann zeta function case}
\end{center}
\end{table}
\begin{table}[h]
\begin{center}
\begin{tabular}{| c | c | c |} \hline
$N$&$\frac{1}{N} \sum_{n \leq N} \Lambda(n + 2; F) \Lambda^*(n; F)$&$\frac{1}{N} \sum_{n \leq N} \Lambda(n + 6; F) \Lambda^*(n; F)$ \\ \hline
450000&1.313531866&2.639461691 \\
460000&1.313570707&2.636315535 \\
470000&1.30917193&2.632942088 \\
480000&1.310726659&2.628272383 \\
490000&1.308466455&2.627803473 \\
500000&1.309116154&2.632674472 \\
510000&1.306101506&2.627889146 \\
520000&1.307762803&2.627914266 \\
530000&1.313036606&2.62448664 \\
540000&1.312406615&2.62678395 \\
550000&1.311065403&2.628826515 \\
560000&1.309988549&2.626247692 \\
570000&1.30962527&2.625719937 \\ 
580000&1.31383141&2.627147271 \\
590000&1.314875333&2.630739868 \\
600000&1.315059206&2.632124451 \\
610000&1.316455007&2.636920329 \\ 
620000&1.314995266&2.637231395 \\
630000&1.315608629&2.632405769 \\
640000&1.316391394&2.633532889 \\
650000&1.31761321&2.634037059 \\
660000&1.316263599&2.630269862 \\
670000&1.319153652&2.630183082 \\
680000&1.318301501&2.627405183 \\ 
690000&1.315417367&2.626615549 \\
700000&1.314202169&2.62761107 \\
710000&1.314910515&2.630445224 \\ 
720000&1.313051006&2.627983697 \\ \hline
\end{tabular} \caption{Arithmetic correlations for the Riemann zeta function case}
\end{center}
\end{table}
\begin{table}[h]
\begin{center}
\begin{tabular}{| c | c | c |} \hline
$N$&$\frac{1}{N} \sum_{n \leq N} \Lambda(n + 2; F) \Lambda^*(n; F)$&$\frac{1}{N} \sum_{n \leq N} \Lambda(n + 6; F) \Lambda^*(n; F)$ \\ \hline
730000&1.313523006&2.631143405 \\
740000&1.312293861&2.628627283 \\
750000&1.312469208&2.629411312 \\
760000&1.310385212&2.628315048 \\
770000&1.308861623&2.630529273 \\
780000&1.310115738&2.630281483 \\
790000&1.310322117&2.627889873 \\
800000&1.311478168&2.631032993 \\
810000&1.309441919&2.630056128 \\
820000&1.30838591&2.629612305 \\
830000&1.310988047&2.628115592 \\ 
840000&1.311759522&2.630026929 \\
850000&1.310568486&2.634801336 \\
860000&1.308997526&2.636703902 \\
870000&1.309846896&2.637436458 \\
880000&1.311129623&2.632355285 \\
890000&1.313463103&2.633220029 \\ 
900000&1.312015897&2.631611954 \\ 
910000&1.311863545&2.632670433 \\
920000&1.312964774&2.631907393 \\
930000&1.314065491&2.633640195 \\
940000&1.312955556&2.630757659 \\
950000&1.313283649&2.629469274 \\
960000&1.312243511&2.632099605 \\ 
970000&1.313004146&2.63236835 \\
980000&1.313382561&2.630539801 \\
990000&1.313197086&2.629649839 \\ 
1000000&1.312844345&2.631198767 \\ \hline
\end{tabular} \caption{Arithmetic correlations for the Riemann zeta function case}
\end{center}
\end{table}

\begin{table}[h]
\begin{center}
\begin{tabular}{| c | c | c |} \hline
$N$&$\frac{1}{N} \sum_{n \leq N} \Lambda(n + 2; F) \Lambda^*(n; F)$&$\frac{1}{N} \sum_{n \leq N} \Lambda(n + 6; F) \Lambda^*(n; F)$ \\ \hline
10000&-0.0189774403&-0.08442782339 \\
20000&-0.03242674706&-0.1814906666 \\
30000&0.05350701315&-0.02274914894 \\ 
40000&0.01899372631&-0.06921723123 \\
50000&0.04266321146&-0.09176519492 \\
60000&-0.003253080071&-0.09175070692 \\
70000&-0.01712846298&-0.1194830594 \\
80000&-0.01744808398&-0.1276901007 \\
90000&-0.04648709848&-0.123258718 \\
100000&-0.04665489689&-0.1027644004 \\
110000&-0.03788852796&-0.1000073675 \\
120000&-0.0371087813&-0.08998876875 \\ 
130000&-0.04734425296&-0.09357298413 \\
140000&-0.03871682266&-0.06263512642 \\
150000&-0.03006101438&-0.05455030481 \\ 
160000&-0.02586546694&-0.05060990349 \\ 
170000&-0.01748479684&-0.03822392085 \\
180000&-0.001181937683&-0.03491899473 \\
190000&0.01225223336&-0.03505783379 \\
200000&0.004885815871&-0.0238265759 \\
210000&-0.00098&-0.02971144814 \\
220000&0.009393308044&-0.0382814219 \\
230000&0.01204287425&-0.03843039904 \\
240000&0.007826859765&-0.04543061557 \\
250000&0.006604528201&-0.04173056499 \\
260000&0.001282254576&-0.03697046738 \\
270000&-0.009535612811&-0.02900872698 \\ 
280000&-0.01475572686&-0.04598304049 \\ \hline
\end{tabular} \caption{Arithmetic correlations for the Ramanujan tau function case}
\end{center}
\end{table}
\begin{table}[h]
\begin{center}
\begin{tabular}{| c | c | c |} \hline
$N$&$\frac{1}{N} \sum_{n \leq N} \Lambda(n + 2; F) \Lambda^*(n; F)$&$\frac{1}{N} \sum_{n \leq N} \Lambda(n + 6; F) \Lambda^*(n; F)$ \\ \hline
290000&-0.0113937803&-0.047239434 \\ 
300000&-0.002842263831&-0.05012316992 \\
310000&-0.005211381643&-0.04677785246 \\
320000&-0.01079038242&-0.04838577936 \\
330000&-0.01130281236&-0.05098659701 \\
340000&-0.01186006184&-0.04707257111 \\
350000&-0.009948712072&-0.05030756301 \\
360000&-0.002680989724&-0.05411621594 \\
370000&-0.001133246745&-0.05630497331 \\
380000&-0.002443019706&-0.05959486012 \\
390000&-0.004979440765&-0.05337992938 \\
400000&0.003722088514&-0.05151054842 \\
410000&-0.001602054897&-0.06002416753 \\
420000&-0.004194878164&-0.05617738769 \\
430000&-0.006473868593&-0.05262061414 \\ 
440000&-0.009209925451&-0.05901353435 \\ 
450000&-0.007689610065&-0.0554534687 \\
460000&-0.007502199808&-0.06306767361 \\
470000&-0.01210060229&-0.06787097275 \\ 
480000&-0.00958419432&-0.06110575304 \\
490000&-0.007056085556&-0.0589300837 \\
500000&-0.007128235852&-0.04879096556 \\
510000&-0.008397314552&-0.0500323917 \\
520000&-0.0158207327&-0.03760700881 \\
530000&-0.01832489405&-0.03923128101 \\
540000&-0.02269515311&-0.03672496854 \\
550000&-0.02376031288&-0.02802674811 \\
560000&-0.01921794053&-0.03224394089 \\ \hline
\end{tabular} \caption{Arithmetic correlations for the Ramanujan tau function case}
\end{center}
\end{table}
\begin{table}[h]
\begin{center}
\begin{tabular}{| c | c | c |} \hline
$N$&$\frac{1}{N} \sum_{n \leq N} \Lambda(n + 2; F) \Lambda^*(n; F)$&$\frac{1}{N} \sum_{n \leq N} \Lambda(n + 6; F) \Lambda^*(n; F)$ \\ \hline
570000&-0.02272588482&-0.03260203354 \\
580000&-0.02517267544&-0.03468682827 \\
590000&-0.02889821002&-0.03919581357 \\
600000&-0.03050035577&-0.03924088897 \\
610000&-0.02758590799&-0.03425085182 \\
620000&-0.02981919213&-0.03576753068 \\
630000&-0.02882483191&-0.0350609601 \\
640000&-0.03345084552&-0.03054331654 \\
650000&-0.03527129109&-0.03059001705 \\
660000&-0.03678650911&-0.03634482904 \\
670000&-0.03715968978&-0.03524981789 \\
680000&-0.03502060296&-0.03816686608 \\
690000&-0.03396350548&-0.03873672407 \\
700000&-0.03470933613&-0.03927460319 \\
710000&-0.03593926547&-0.03669712524 \\
720000&-0.03719309038&-0.03689551867 \\
730000&-0.03504866338&-0.04193547039 \\ 
740000&-0.03360291577&-0.04533385386 \\
750000&-0.03212371828&-0.04402232818 \\
760000&-0.03194360359&-0.04165550128 \\
770000&-0.03124277526&-0.04019531646 \\
780000&-0.0319127759&-0.03793252369 \\
790000&-0.03280035278&-0.04008243438 \\
800000&-0.03220979131&-0.03464642662 \\
810000&-0.03358457667&-0.0361304424 \\
820000&-0.0317043789&-0.03592271016 \\
830000&-0.03218264369&-0.03340236767 \\
840000&-0.03334972742&-0.02938632407 \\ \hline
\end{tabular} \caption{Arithmetic correlations for the Ramanujan tau function case}
\end{center}
\end{table}
\begin{table}[h]
\begin{center}
\begin{tabular}{| c | c | c |} \hline
$N$&$\frac{1}{N} \sum_{n \leq N} \Lambda(n + 2; F) \Lambda^*(n; F)$&$\frac{1}{N} \sum_{n \leq N} \Lambda(n + 6; F) \Lambda^*(n; F)$ \\ \hline
850000&-0.0328039171&-0.03097091229 \\
860000&-0.03121950705&-0.03507015215 \\
870000&-0.03317460071&-0.03403036823 \\
880000&-0.03372010312&-0.03232125634 \\
890000&-0.03203089945&-0.03732169961 \\
900000&-0.03076232625&-0.03688399552 \\
910000&-0.03122602894&-0.03589739301 \\
920000&-0.03094921463&-0.03421212047 \\
930000&-0.02841595948&-0.03361505122 \\
940000&-0.03081887357&-0.03051921045 \\
950000&-0.0315540321&-0.03250261653 \\
960000&-0.0309007989&-0.03324071654 \\ 
970000&-0.03012839675&-0.02979650751 \\
980000&-0.02973947929&-0.02418204924 \\
990000&-0.02835974557&-0.02427138206 \\
1000000&-0.02808424592&-0.02517614702 \\ \hline
\end{tabular} \caption{Arithmetic correlations for the Ramanujan tau function case}
\end{center}
\end{table}

\begin{table}[h]
\begin{center}
\begin{tabular}{| c | c | c |} \hline
$N$&$\frac{1}{N} \sum_{n \leq N} \Lambda(n + 2; F) \Lambda^*(n; F)$&$\frac{1}{N} \sum_{n \leq N} \Lambda(n + 6; F) \Lambda^*(n; F)$ \\ \hline
10000&0.001360693324&-0.1205032398 \\
20000&0.005526608274&-0.147000187 \\
30000&0.03702803805&0.01294441124 \\
40000&0.04487384484&-0.006133230287 \\
50000&0.007950980142&0.01763560136 \\ 
60000&0.03857765089&0.0299684617 \\
70000&0.03894359829&0.02436846415 \\ 
80000&0.03199803913&0.00941667465 \\ \hline
\end{tabular} \caption{Arithmetic correlations for the elliptic curve case E11.a3}
\end{center}
\end{table}
\begin{table}[h]
\begin{center}
\begin{tabular}{| c | c | c |} \hline
$N$&$\frac{1}{N} \sum_{n \leq N} \Lambda(n + 2; F) \Lambda^*(n; F)$&$\frac{1}{N} \sum_{n \leq N} \Lambda(n + 6; F) \Lambda^*(n; F)$ \\ \hline
90000&0.01616671351&0.02193135939 \\
100000&0.003366518021&0.05376886526 \\ 
110000&-0.004820349393&0.03822876044 \\
120000&-0.01435434687&0.04072244619 \\
130000&-0.01240894568&0.03598964429 \\
140000&0.004674537615&0.04101144451 \\
150000&-0.001836541417&0.04142766238 \\ 
160000&-0.01493880436&0.04086650171 \\
170000&-0.009496183087&0.02476938222 \\
180000&-0.01874572763&0.02481055577 \\
190000&-0.01948604441&0.02529319654 \\
200000&-0.02391289435&0.02958925117 \\
210000&-0.01936071972&0.02068020678 \\
220000&-0.01850873577&0.01747038783 \\
230000&-0.01391225145&-0.001849554641 \\
240000&-0.01897856564&-0.000784 \\ 
250000&-0.02128396007&0.006196869743 \\
260000&-0.02056421245&-0.0000998 \\
270000&-0.0239424774&-0.000226 \\
280000&-0.01294407074&0.009335738365 \\ 
290000&-0.01823237133&0.01535360847 \\ 
300000&-0.02239250241&0.01701863976 \\
310000&-0.02356350482&0.02512130429 \\
320000&-0.02416480064&0.03134765443 \\
330000&-0.03219216843&0.0212919432 \\
340000&-0.03758607637&0.01163956319 \\
350000&-0.04440585416&0.01811200314 \\
360000&-0.03644117167&0.01592575687 \\ \hline
\end{tabular} \caption{Arithmetic correlations for the elliptic curve case E11.a3}
\end{center}
\end{table}
\begin{table}[h]
\begin{center}
\begin{tabular}{| c | c | c |} \hline
$N$&$\frac{1}{N} \sum_{n \leq N} \Lambda(n + 2; F) \Lambda^*(n; F)$&$\frac{1}{N} \sum_{n \leq N} \Lambda(n + 6; F) \Lambda^*(n; F)$ \\ \hline
370000&-0.04242407093&0.008217930691 \\
380000&-0.04415965109&0.0135441231 \\
390000&-0.04638556444&0.009492067397 \\
400000&-0.04045402702&0.01687636346 \\
410000&-0.03661715982&0.02029750119 \\
420000&-0.03724763127&0.01413224349 \\
430000&-0.03621173446&0.02124933408 \\
440000&-0.04083197481&0.01741707557 \\
450000&-0.03955737675&0.02015212998 \\
460000&-0.04023054229&0.01990697603 \\
470000&-0.03899391608&0.0218177904 \\
480000&-0.03528814401&0.02215101745 \\
490000&-0.03652044498&0.02251188652 \\
500000&-0.03590633097&0.02217368551 \\
510000&-0.03201310385&0.03222708836 \\
520000&-0.03044086494&0.02709231454 \\
530000&-0.02924381283&0.02683025438 \\
540000&-0.02785780813&0.02541288526 \\
550000&-0.02261602292&0.027261404 \\
560000&-0.0213661269&0.02263497032 \\
570000&-0.02593604881&0.02515762476 \\
580000&-0.02128891545&0.02791486977 \\
590000&-0.01770654279&0.03281587457 \\
600000&-0.01794452524&0.03193556161 \\
610000&-0.01971965974&0.03885585589 \\
620000&-0.02107177847&0.03944204942 \\
630000&-0.02188074493&0.0370770951 \\
640000&-0.0198457308&0.02763544504 \\ \hline
\end{tabular} \caption{Arithmetic correlations for the elliptic curve case E11.a3}
\end{center}
\end{table}
\begin{table}[h]
\begin{center}
\begin{tabular}{| c | c | c |} \hline
$N$&$\frac{1}{N} \sum_{n \leq N} \Lambda(n + 2; F) \Lambda^*(n; F)$&$\frac{1}{N} \sum_{n \leq N} \Lambda(n + 6; F) \Lambda^*(n; F)$ \\ \hline
650000&-0.02081392246&0.02787621807 \\
660000&-0.01977918882&0.02555723798 \\
670000&-0.02239569219&0.0214671847 \\
680000&-0.02335690979&0.02060756263 \\
690000&-0.02137447227&0.0182869915 \\
700000&-0.02076795813&0.01683291091 \\
710000&-0.02001433019&0.0171875706 \\
720000&-0.02049796502&0.01821948635 \\
730000&-0.01760380482&0.02018539935 \\
740000&-0.01660947537&0.01982335679 \\
750000&-0.01657891183&0.02143069221 \\
760000&-0.0171344018&0.02245305394 \\
770000&-0.01784892258&0.02233721615 \\
780000&-0.01720520961&0.02581298687 \\
790000&-0.01752178322&0.02452948772 \\
800000&-0.01778293828&0.02465176887 \\
810000&-0.01552380393&0.02901224375 \\
820000&-0.01288973138&0.03030389132 \\
830000&-0.01003305218&0.02717794237 \\
840000&-0.01350272027&0.02602559229 \\
850000&-0.01282143442&0.02419650131 \\
860000&-0.01386408302&0.02420635351 \\
870000&-0.009836966878&0.02344567239 \\
880000&-0.01335636347&0.02050915349 \\
890000&-0.009755359563&0.01802602114 \\
900000&-0.01182808099&0.0168340123 \\
910000&-0.009186812824&0.01067370944 \\
920000&-0.006378927412&0.01118558399 \\ \hline
\end{tabular} \caption{Arithmetic correlations for the elliptic curve case E11.a3}
\end{center}
\end{table}
\begin{table}[h]
\begin{center}
\begin{tabular}{| c | c | c |} \hline
$N$&$\frac{1}{N} \sum_{n \leq N} \Lambda(n + 2; F) \Lambda^*(n; F)$&$\frac{1}{N} \sum_{n \leq N} \Lambda(n + 6; F) \Lambda^*(n; F)$ \\ \hline
930000&-0.007024807306&0.009276106156 \\
940000&-0.006578997168&0.01012236377 \\
950000&-0.005227612016&0.01074753742 \\
960000&-0.003634826018&0.007321265298 \\
970000&-0.002747957287&0.006915929712 \\
980000&0.000408&0.006821232812 \\
990000&0.002042736259&0.006259319362 \\
1000000&0.00071&0.007244883054 \\ \hline
\end{tabular} \caption{Arithmetic correlations for the elliptic curve case E11.a3}
\end{center}
\end{table}

\chapter{Exponential weight}
Consideration is now given to
\begin{align}
\tilde{I}(X,T) = & \int_{|t| \leq T} \left| \sum_{n \geq 1} X^{i t(n; F)} \exp\left(- 2 (t - t(n; F))^2\right) 1_{|t(n; F)| \leq Z} \right|^2 \, dt \nonumber \\
= & \sum_{m, n \geq 1} X^{i (t(m; F) - t(n; F))} 1_{|t(m; F)|, |t(n; F)| \leq Z} \nonumber \\ 
& \int_{|t| \leq T} \exp\left(- 2 ((t - t(m; F))^2 + (t - t(n; F))^2)\right) \, dt.
\end{align}
Following the approach in chapter $3,$ using the density of zeros to restrict the summation to those with $|t(m; F)| \leq Z, |t(n; F)| \leq Z,$ results in an error of size $\ll (\log T)^2.$ Also, extending the domain of integration to the whole real line gives an error of size $\ll (\log T)^3$. Therefore
\begin{align}
\tilde{I}(X,T) = & \sum_{m , n \geq 1} X^{i (t(m; F) - t(n; F))} 1_{|t(m; F)|, |t(n; F)| \leq Z} \nonumber \\
& \int \exp\left(- 2 ((t - t(m; F))^2 + (t - t(n; F))^2)\right) \, dt + O \left( \left( \log T \right)^3 \right) \nonumber \\
= & \sum_{m, n \geq 1} X^{i (t(m; F) - t(n; F))} \exp\left(- 2((t(m; F))^2 + (t(n; F))^2)\right) \nonumber \\
& 1_{|t(m; F)|, |t(n; F)| \leq Z} \int \exp(- 4 t^2 + 4 (t(m; F) + t(n; F)) t) \, dt \nonumber \\
& + O \left( \left( \log T \right)^3 \right).
\end{align}
For $a \in \mathbb{R}^+$ and $b \in \mathbb{C}$ the integral identity
\begin{equation}
\int \exp(- a t^2 + bt) \, dt = \exp(b^2 / (4 a)) \sqrt{\pi / a}
\end{equation}
may now be applied with $a = 4$ and $b = 4 (t(m; F) + t(n; F))),$ giving
\begin{equation}
\tilde{I}(X, T) = \frac{\sqrt{\pi}}{2} \tilde{\mathcal{F}}(X, T; F) + O \left( \left( \log T \right)^3 \right).
\end{equation}
Considering $\mathcal{F}(X, T; F)$ for example in the regime $\textrm{deg}(F) < A_1 < A_2$ then gives
\begin{align}
& \int \left| a(it) \right|^2 \left| \sum_{n \geq 1} X^{i t(n; F)} \exp\left(- 2 (t - t(n; F))^2\right) 1_{|t(n; F)| \leq Z}\right|^2 \, dt \nonumber \\
= & \sqrt{\pi} \kappa \left(\textrm{deg}(F) \log \frac{1}{\kappa} + \log \mathfrak{q}(F) + (1 - \gamma - \log 4 \pi) \textrm{deg}(F)\right) + O \left( \kappa^{1 + c} \right) \nonumber \\
& + O_{\epsilon} \left( \kappa^{1 + c_1 - \epsilon} \right) + O_{\epsilon} \left( \kappa^{1 + c_2 - \epsilon}\right) \nonumber \\
= & \frac{\sqrt{\pi}}{2} \delta \left(\textrm{deg}(F) \log \frac{1}{\delta} + \log \mathfrak{q}(F) + (1 - \gamma - \log 2 \pi) \textrm{deg}(F) \right) + O \left( \delta^{1 + c} \right) \nonumber \\
& + O_{\epsilon} \left( \delta^{1 + c_1 - \epsilon} \right) + O_{\epsilon} \left( \delta^{1 + c_2 - \epsilon} \right),
\end{align}
where Lemma A has been used for 
\begin{equation}
\delta^{- (1 - c_1)} \ll T \ll \delta^{- (1 + c_2)},
\end{equation}
for some $0 < c_1 < 1 \textrm{ and } 0 < c_2 < 1.$
Therefore
\begin{align}
&\int \left| \sum_{n \geq 1} a(1 / 2 + i t(n; F)) X^{i t(n; F)} \exp\left(- (t - t(n; F))^2\right) 1_{|t(n; F)| \leq Z} \right|^2 \, dt \nonumber \\ 
= &\frac{\sqrt{\pi}}{2} \delta \left(\textrm{deg}(F) \log \frac{1}{\delta} + \log \mathfrak{q}(F) + (1 - \gamma - \log 2 \pi) \textrm{deg}(F)\right) \nonumber \\
&+ O \left( \delta^{1 + c} \right) + O_\epsilon \left( \delta^{ 1 + c_1 - \epsilon} \right) + O_\epsilon \left( \delta^{1 + c_2 - \epsilon} \right).
\end{align}
Making use of the Fourier transform,
\begin{equation}
\int \exp(- (t - b)^2) \exp(- 2 \pi i k t) \, dt = \sqrt{\pi} \exp(- \pi^2 k^2) \exp(- 2 \pi i k b),
\end{equation}
an application of Plancherel's formula then results in
\begin{align}
& \int \left|\sum_{n \geq 1} a(1 / 2 + i t(n; F)) X^{i t(n; F)} \exp(- (t - t(n; F))^2) 1_{|t(n; F)| \leq Z}\right|^2 \, dt \nonumber \\
= & \pi \int \left|\sum_{n \geq 1} a(1 / 2 + i t(n; F)) X^{i t(n; F)} \exp(- \pi^2 u^2) \exp(- 2 \pi i t(n; F) u) \right|^2 \, du \nonumber \\
= & \pi \int \left|\sum_{n \geq 1} a(1 / 2 + i t(n; F)) \exp(i t(n; F) (Y + y))\right|^2 \exp(- y^2 / 2) \, dy,
\end{align}
using the substitution $Y = \log X$ and $y = - 2 \pi u.$

The weight now differs to the one removed by Lemma B in chapter $3.$ To proceed Lemma B requires some modification, given a more general initial statement than provided by Lemma B,
\begin{equation}
f(T) = \int g(T + y) h(y) \, dy = O(\exp(- c_1 T)),
\end{equation}
where $0 < c_1 < 1,$ the choice of $H$ such that
\begin{equation}
F(T) = \int g(T + y) H(y) \, dy = O(\exp(- c_2 T)),
\end{equation}
where $0 < c_2 < 1,$ is sought. $F(T) = O(\exp(- c_2 T))$ if $\mathscr{F}[F](s)$ is bounded and analytic in a strip. By the convolution theorem
\begin{equation}
\mathscr{F}[F](s) = \frac{\mathscr{F}[f](s) \mathscr{F}[H](s)}{\mathscr{F}[h](s)},
\end{equation}
therefore $F(T) = O(\exp(- c_2 T))$ if $\frac{\mathscr{F}[H](s)}{\mathscr{F}[h](s)}$ is bounded which will occur if $\mathscr{F}[H](s)$ decays more rapidly than $\mathscr{F}[h](s).$
		 		 
\bibliography{BibTeX} 

\bibliographystyle{alpha}

\end{document}